\documentclass[final,3p,times,number]{elsarticle}
\biboptions{sort&compress}
\journal{Elsevier}

\usepackage[utf8]{inputenc}
\usepackage{gensymb}
\usepackage{hyphenat}
\usepackage{subcaption}
\usepackage{multirow}
\usepackage{graphicx}
\usepackage{bm}
\usepackage{amsmath,amssymb,amsfonts}
\usepackage{setspace}
\usepackage{tikz}
\usetikzlibrary{arrows, arrows.meta}
\usepackage{algorithm, algpseudocode}
\usepackage[normalem]{ulem}
\usepackage{comment}
\usepackage{booktabs} 
\usepackage{colortbl}

\graphicspath{ {./figs/}  
               {./figs/ConditionNumber_EmbeddedCircle/}
               {./figs/Poisson-ClippedMesh/}
               {./figs/Hole/} 
               {./figs/Bimaterial/}
               {./figs/Agglomeration}
}

\renewcommand{\Re}{\mathbb{R}}

\newcommand{\sref}[1]{Section~\ref{#1}}

\newcommand{\fref}[1]{Fig.~\ref{#1}}
\newcommand{\aref}[1]{Algorithm~\ref{#1}}

\newcommand{\vm}[1]{{\bm{#1}}}
\newcommand{\vx}{\vm{x}}

\usepackage[colorlinks=true]{hyperref}

\usepackage{color}
\newcommand{\revised}[1]{{#1}}

\let\oldsout\sout
\renewcommand{\sout}[1]{\textcolor{red}{\oldsout{#1}}}

\makeatletter
\newenvironment{breakablealgorithm}
{
  \begin{center}
    \refstepcounter{algorithm}
    \hrule height.8pt depth0pt \kern2pt
    \renewcommand{\caption}[2][\relax]{
      {\raggedright\textbf{\fname@algorithm~\thealgorithm} ##2\par}%
      \ifx\relax##1\relax 
      \addcontentsline{loa}{algorithm}{\protect\numberline{\thealgorithm}##2}%
      \else 
      \addcontentsline{loa}{algorithm}{\protect\numberline{\thealgorithm}##1}%
      \fi
      \kern2pt\hrule\kern2pt
    }
  }{
    \kern2pt\hrule\relax 
  \end{center}
}
\makeatother

\algnewcommand{\IfThenElse}[3]{
{<else>}
  \State \algorithmicif\ #1\ \algorithmicthen\ #2\ \algorithmicelse\
  #3}

\algnewcommand{\IfThen}[2]{
  \State \algorithmicif\ #1\ \algorithmicthen\ #2}

\usepackage{eqparbox}
\newdimen{\algindent}
\algnewcommand\LeftComment[2]{\hspace{#1\algindent}$\triangleright$ \eqparbox{COMMENT}{\small{\tt #2}} \hfill }
\eqsetmaxwidth{COMMENT}{\linewidth}
\algnewcommand\comment[1]{\hfill$\triangleright~${\small{\tt #1}}}

\errorcontextlines\maxdimen
\makeatletter
\newcommand*{\algrule}[1][\algorithmicindent]{\makebox[#1][l]{\hspace*{.5em}\vrule height .75\baselineskip depth .25\baselineskip}}%

\newcount\ALG@printindent@tempcnta
\def\ALG@printindent{%
  \ifnum \theALG@nested>0
  \ifx\ALG@text\ALG@x@notext
  \addvspace{-3pt}
  \else
  \unskip
  \ALG@printindent@tempcnta=1
  \loop
  \algrule[\csname ALG@ind@\the\ALG@printindent@tempcnta\endcsname]%
  \advance \ALG@printindent@tempcnta 1
  \ifnum \ALG@printindent@tempcnta<\numexpr\theALG@nested+1\relax
  \repeat
  \fi
  \fi
}%
\usepackage{etoolbox}
\patchcmd{\ALG@doentity}{\noindent\hskip\ALG@tlm}{\ALG@printindent}{}{\errmessage{failed to patch}}
\makeatother

\begin{document}

\doublespacing  

\title{CutVEM: Conforming virtual element method on embedded domains
  with shape-agnostic element agglomeration}
  
\author[1]{Ramsharan Rangarajan\corref{cor}}
\ead{rram@iisc.ac.in}

\author[2]{N.\ Sukumar\corref{cor}}
\ead{nsukumar@ucdavis.edu}

\cortext[cor]{Corresponding author}

\address[1]{Mechanical Engineering, Indian Institute of Science Bangalore, Karnataka 560012, INDIA}

\address[2]{Department of Civil and Environmental Engineering,  One Shields Avenue,
University of California, Davis, CA 95616, USA}

\begin{abstract}
  The virtual element method (VEM) is a stabilized Galerkin method
  that is accurate and robust on general polygonal meshes.  This
  feature makes it an appealing candidate for simulations involving
  meshes with embedded interfaces and evolving geometries. 
  However, similar to the finite element method,
  in such scenarios the VEM can also yield poorly conditioned stiffness matrices  due to meshes
  having cut cells. With the objective of developing an embedded
  domain method, we propose a novel element agglomeration algorithm for the virtual element method to
  address this issue.  The agglomeration algorithm renders the VEM
  robust over planar polygonal meshes, particularly on finite element
  meshes cut by immersed geometries.  The algorithm relies on the
  element stability ratio, which we define using the extreme
  eigenvalues of the element stiffness matrix. The resulting element
  agglomeration criterion is free from nebulous polygon quality
  metrics and is defined independently of polygon shapes. The
  algorithm proceeds iteratively and element-wise to maximize the
  minimum element stability ratio, even at the expense of degrading
  elements with better ratios. Crucially, element agglomeration alters
  the number of elements, not the degree of freedom count. The
  resulting method, which we label as CutVEM, retains node locations 
  of cut elements unchanged, and therefore yields discretizations
  that conform to embedded interfaces. This, in turn, facilitates
  straightforward imposition of boundary conditions and interfacial
  constraints. Through detailed numerical experiments that sample
  varied element-interface intersections, we demonstrate that CutVEM
  enjoys dramatically improved condition numbers of global stiffness
  matrices over the VEM. Furthermore, simulations of prototypical
  heat conduction problems with Dirichlet and Neumann boundary conditions on domains
  with immersed geometries show that element agglomeration does not
  noticeably degrade solution accuracy and that CutVEM retains the
  VEM's optimal convergence rate.
\end{abstract}


\maketitle

\section{Introduction}\label{sec:intro}
There has been significant effort over the past two decades to
  develop and advance embedded domain and embedded interface methods,
  also referred to as immersed boundary~\cite{Peskin:2002:TIB} or
  fictitious domain methods~\cite{Glowinski:1994:AFD}.  A review of
  these methods can be found in~\cite{Mittal:2005:IBM}. They aim to
  simplify mesh generation by embedding complex geometries in
  \revised{nonconforming} meshes---often structured (Cartesian) or
  high-quality Delaunay meshes~\cite{Bishop:2003:RSA}.  Embedding
  geometries this way, however, generates cut-cells along boundaries
  and interfaces. Several
  challenges then arise:
\begin{enumerate}[(i)]
\item Constructing approximation spaces by defining basis functions
  over small cut cells can result in poorly conditioned systems of
  equations~\cite{dePrenter:2017:CNA}. The ghost penalty
  stabilization~\cite{Burman:2010:GP}, for instance, addresses this
  issue.
\item Imposing Dirichlet boundary conditions or jump conditions along
cut-cell boundaries require special techniques. Nitsche's method
is a popular choice to this
end~\cite{Dolbow:2009:AEF,Sanders:2012:ANE,Ruberg:2016:UIM}.
\item Computing integrals over cut-elements, which may even be
nonconvex, demands unconventional yet robust numerical integration
techniques, cf.~\cite{Fries:2017:HOM}.
\end{enumerate}

Numerous methods have been proposed in the literature to address the
challenges highlighted. These include the extended finite element
methods for cracks~\cite{Moes:1999:AFE} and material
interfaces~\cite{Sukumar:2001:MHI}, unfitted finite element
methods~\cite{Hansbo:2002:UFE}, embedded discontinuous Galerkin
methods~\cite{Lew:2008:ADI,Rangarajan:2009:DGB}, finite cell
method~\cite{Parvizian:2007:FCM,Duster:2008:TFC}, embedded boundary
techniques~\cite{Wang:2011:AIT}, carving methods on hexahedral
meshes~\cite{Sohn:2013:FES}, CutFEM~\cite{Burman:2015:CDG},
CutIGA~\cite{Elfverson:2018:CUT}, shifted boundary
method~\cite{Main:2018a:SBM,Main:2018b:SBM} and the extreme mesh
deformation method~\cite{Moes:2023:XME,Quiriny:2024:XME}.  Among
these, the shifted boundary method is notable for its approach to
circumvent issues of poor conditioning caused by cut-cells and
numerical integration over them. This method introduces a \emph{shifted
  boundary} that conforms to the mesh, but incurs the task of
reformulating boundary conditions on the shifted boundary.  A review
on the treatment of small cut-elements in embedded interface methods
can be found in~\cite{dePrenter:2023:SCI}.

The methods listed above notwithstanding, the simplicity and
robustness of domain conforming decompositions are hard to overlook~\cite{Noble:2010:ACD,Roberts:2018:AVC,Fries:2018:HCD}. Such methods
provide critical advantages when simulating embedded boundary and
embedded (multimaterial) interface problems. Two crucial challenges
persist in these methods: ill-conditioning caused by cut-elements and
numerical integration over these cut-elements. In this work, we
introduce an element agglomeration algorithm for the virtual element
method~\cite{Beirao:2013:BPV,Ahmad:2013:EPV,Beirao:2014:HGV,Sukumar:2022:AVE}
to address the former issue and adopt the scaled boundary cubature
scheme of~\cite{Chin:2021:SBC} for the latter. The resulting approach,
henceforth referred to as CutVEM, is a virtual element method that is
accurate and robust to element-boundary intersections. Our main
contribution is a shape-agnostic element agglomeration algorithm for
polygonal meshes based on a novel agglomeration criterion, and a
detailed examination of CutVEM's performance over embedded domains
in two dimensions. 

The conforming cut-cell embedded domain method is based on the
construction of agglomerated virtual elements~\cite{Sukumar:2022:AVE}.
The virtual element method (VEM)~\cite{Beirao:2013:BPV} is a
stabilized Galerkin method that is accurate and robust on general
convex and nonconvex polytopal meshes.  For a recent review on the
VEM, the interested reader can refer to~\cite{Beirao:2023:VEM}. In
particular, VEM is known to be robust on $n$-sided polygons
($n \ge 4$) with short edges.  On cutting a Delaunay or quadrangular
finite element mesh, poorly-shaped triangles (interior angle close to
$0^\circ$ or $180^\circ$) or quadrilaterals can arise. 
In embedded methods, the causes and
consequences of ill-conditioning due to small cut-elements are
elaborated in~\cite{Schillinger:2015:TFC}.\footnote{ In the words of
  Schillinger and Reuss~\cite{Schillinger:2015:TFC}: ``Cut elements
  deteriorate the conditioning of the system matrix. Especially those
  elements that are cut in such a way that only a very small part of
  their element domain is located in the physical domain can lead to
  very ill-conditioned system matrices. The finite cell method
  therefore requires the solution of the system of equations by direct
  solvers to guarantee an accurate and stable solution.''}  To remove
such `bad' elements, we agglomerate them with
their neighbors to create polygons having four or more
edges.

The inspiration for the agglomeration algorithm presented here comes
from the work in~\cite{Sukumar:2022:AVE}.  This study highlighted that
using the maximum element-eigenvalue to set the critical time step for
explicit dynamics is impractical for tetrahedral meshes with
slivers. It was demonstrated that agglomeration dramatically decreases
the maximum element-eigenvalue (square of the maximum natural
frequency of the element). Consequently, a larger (computationally
feasible) critical time step was possible for explicit linear
elastodynamic simulations. More recently, this notion of agglomeration
of sliver tetrahedra has been successfully exploited to render
Lagrangian simulations of fluid dynamics to be computationally
efficient~\cite{Fu:2025:PVE}.  Agglomeration alters the number of
elements and their connectivities, but preserves the number of mesh
vertices and hence the degree of freedom count. In 2D, only the area
of the polygon needs to be computed for first-order virtual elements;
however, the approach can be readily extended for higher-order virtual
elements since accurate numerical integration over general polygonal
elements can be efficiently carried out using the scaled boundary
cubature scheme~\cite{Chin:2021:SBC,Chin:2023:NIV}.

We point out that the idea of agglomeration on polygonal meshes has
received significant attention in the recent literature. Noteworthy in
this regard is the work in~\cite{sorgente2023mesh,Sorgente:2024:MOV},
where agglomeration is used on polygonal meshes (modeling large-scale
discrete fractures) to reduce the number of elements without
compromising accuracy.  When compared to the agglomeration algorithm
in~\cite{sorgente2023mesh,Sorgente:2024:MOV}, our criterion
distinguishes itself in several ways: (1) it is element-specific while
being shape-agnostic; (2) circumvents the need for quality metric of
polygons, which is inherently problematic; (3) requires computing for
the extreme eigenvalues of element-level matrices; (4) improves the
\emph{stability ratio} (inverse of the condition number) of poorer
elements, albeit at the expense of better ones; (5) can be performed
as a preprocessing operation and hence is straightforward to
incorporate into existing codes; (6) improves condition number of
cut-elements without altering the nodal locations so that geometric
features are preserved; and (7) retains solution accuracies for
reasonable stability ratios.

The remainder of this paper is structured as follows.  The essentials
of the first-order VEM are presented in~\sref{sec:VEM}.
In~\sref{sec:3}, we present the element agglomeration algorithm. We
highlight the novel agglomeration criterion based on the condition
number of the element-level stiffness matrix introduced
in~\sref{sec:VEM}. In \sref{sec:3-4}, we examine the algorithm's
performance in improving the condition numbers of global stiffness
matrices in the virtual element method using an extensive set of
numerical experiments over meshes with poorly-shaped elements and
embedded \revised{(stationary and evolving)} geometries. 
Numerical experiments are presented in~\sref{sec:results} to establish the accuracy and robustness of CutVEM. First, the benefits of agglomeration to reduce interpolation and finite element approximation errors
over Delaunay meshes (anisotropic
refinement) of flat triangles) are assessed in~\sref{subsec:agg_tests}. Then in
Sections~\ref{subsec:homogeneous_heat_conduction} 
and~\ref{subsec:inhomogeneneous_heat_conduction}, we
present simulations for
prototypical heat conduction problems with Dirichlet and Neumann boundary conditions
on domains with immersed geometries to examine the 
influence of element agglomeration on solution accuracies 
and convergence rates of the CutVEM. We close with a summary of our main findings
in~\sref{sec:conclusions}.

\section{First-order virtual element method}\label{sec:VEM}
In the interest of completeness and for subsequent use in
  \sref{sec:3} and \ref{sec:results}, we include a concise
summary of the first-order virtual element method in the context of
the canonical Poisson boundary-value problem posed over a planar domain \revised{$\Omega$, whose weak form is}:
\begin{align}
  {\rm Find}~ u\in 
   {\rm V} = {\rm H}^1_0(\Omega)
  ~\text{such that} ~{a}(u,v) =
  (f,v)~\text{for all}~v\in {\rm H}^1_0(\Omega), \label{eq:2-1}
\end{align}
where ${a}(u,v) = \int_\Omega\nabla u \cdot \nabla v\, d \vx$, 
$(f,v) = \int_\Omega fv\,d \vx $ 
for a given
$f:\Omega\rightarrow {\mathbb R}$, and ${\rm H}^1_0(\Omega)$ is the
subspace of ${\rm H}^1(\Omega)$ satisfying the homogeneous Dirichlet
boundary condition on $\partial\Omega$.

For the ensuing discussion, it is necessary to introduce some
terminology and notation. Throughout, we reserve the letter 
${\rm E}$
to denote a virtual element and the letter ${\rm K}$ for the
corresponding planar polygon over which ${\rm E}$ is defined. A
polygon ${\rm K}$ is a nonempty and closed subset of ${\mathbb R}^2$,
and is bounded by a closed polygonal chain.
The ${\rm N_K}$ vertices and edges of ${\rm K}$ are enumerated
counterclockwise and are denoted by $v_i$ and $e_i$ for
$i=1,\ldots,N_{\rm K}$, respectively. The Cartesian coordinates of
$v_i$ is $\vx_i\equiv (x_i,y_i)$ and the normal to $e_i$ is given by
$\vm{n}_{i}$. The area and diameter of ${\rm K}$ are $|{\rm K}|$ and
$h_{\rm K}$, respectively. The location of the geometric centroid of
the element is given by $\vx_{\rm K}$. We refer to the boundary of
${\rm K}$ by $\partial {\rm K}$ and an edge belonging to the boundary
as $e\in \partial{\rm K}$.

\subsection{Virtual element space and projection operator}
Given a polygonal discretization ${\cal T}_h$ of $\Omega$, the virtual
element method introduces an approximation of \eqref{eq:2-1} as:
\begin{align}
  {\rm Find}~u_h\in {\rm V}_h~\text{such that}~{a}^h(u_h,v_h) =
  {\rm L}^h(v_h)~\text{for all}~v_h\in {\rm V}_h, \label{eq:2-2}
\end{align}
where the approximation space ${\rm V}_h\subset {\rm V}$ is
constructed polygon-wise, as
\begin{align}
  {\rm V}_h\big|_{\rm K} \equiv {\rm V}_h({\rm K}) = \Bigl\{ v_h\in {\rm
  H}^1({\rm K})\,:\,\Delta v_h = 0, v_h\big|_e\in {\mathbb
  P}_1(e)~\forall e\in \partial {\rm K}, v_h\big|_{\partial{\rm
  K}}\in {\rm C}^0(\partial {\rm K}) \Bigr\} \label{eq:2-3},
\end{align}
with ${\mathbb P}_1(\omega)$ being the usual set of affine functions
on the set $\omega$.  
It is 
known, for instance, that ${\rm V}_h({\rm K})$ coincides with the
space of harmonic generalized barycentric
coordinates~\cite{Floater:2015:GBC}.

The first-order VEM defines the discrete bilinear and linear forms
${a}^h(\cdot,\cdot)$ and ${\rm L}^h(\cdot)$ appearing in
  \eqref{eq:2-2} in a polygon-wise manner.  Over a polygon ${\rm K}$, 
  the restrictions of the bilinear and linear forms in
  \eqref{eq:2-1} and \eqref{eq:2-2} to ${\rm K}$ are given by
\begin{subequations}\label{eq:2-4}
  \begin{align}
    a^h_{\rm K}(v_h,w_h) &= {a}_{\rm K}(\Pi_1^\nabla v_h,\Pi_1^\nabla w_h) +
                                 \tau (v_h-\Pi_1^\nabla v_h,w_h-\Pi_1^\nabla w_h)_{\rm
                            K},  \label{eq:2-4a}\\
    {\rm L}^h_{\rm K}(v_h) &= (\Pi_0f,\overline{v_h})_{\rm K} ,\label{eq:2-4b}
  \end{align}
\end{subequations}
with $\tau>0$ being a
stabilization parameter chosen to be sufficiently large to render
${a}^h(\cdot,\cdot) \equiv\sum_{{\rm K}\in {\cal T}_h}a^h_{\rm
    K}(\cdot,\cdot)$ coercive on ${\rm V}_h$. 
The discrete
  linear form is similarly
  ${\rm L}^h(\cdot) \equiv\sum_{{\rm K}\in{\cal T}_h} {\rm L}^h_{\rm
    K}(\cdot)$. 
In~\eqref{eq:2-4b}, we have used $\overline{v}_h$ to
denote the nodal average
$(1/{\rm N_K})\sum_{j=1}^{\rm N_K}v_h({\bf x}_j)$ of
$v_h\in {\rm V}_h({\rm K})$. The elliptic projection operator
$\Pi_1^\nabla :{\rm V}_h({\rm K})\rightarrow {\mathbb P}_1({\rm K})$
and the $L^2$-orthogonal projection
$\Pi_0:{\rm V}_h({\rm K})\rightarrow {\mathbb R}$ play crucial roles
in the method. Their choices enable the VEM to bypass explicitly
determining a basis $\{\varphi_1,\ldots,\varphi_{\rm N_K}\}$ spanning
${\rm V}_h({\rm K})$; instead, the projections of these basis functions
onto the affine space spanned by the scaled monomials
$m_1=1, m_2 = (x-x_{\rm K})/h_{\rm K}$ and
$m_3=(y-y_{\rm K})/h_{\rm K}$ suffice. The elliptic projection
operator $\Pi_1^\nabla$ satisfies
\begin{subequations}\label{eq:2-5}
  \begin{align}
    a_{\rm K}(v_h-\Pi_1^\nabla v_h, m) &= 0 ~\text{for all}~m\in {\rm M}_{\rm
                                 K}= {\textrm{Span}}\{m_1,m_2,m_3\}, \label{eq:2-5a}\\
    \overline{\Pi_1^\nabla v_h} &= \overline{v_h} , \label{eq:2-5b}
  \end{align}
\end{subequations}
while the projection $\Pi_0$ onto constants is defined such that
$(\Pi_0v_h-v_h,1)=0$.

\subsection{Element stiffness matrix and element load vector}
In the interest of brevity, we directly provide expressions for
computing the element stiffness matrix ${\bf K}_{\rm E}$ and the load
vector ${\bf f}_{\rm E}$. Omitting details of the steps involved in
deriving these from~\eqref{eq:2-4}
and~\eqref{eq:2-5}, we have
\begin{align} 
  {\bf K}_{\rm E}
  &= \underbrace{\vm{\Pi}_*^\top \tilde{\bf G}
    \vm{\Pi}_*}_{{\bf K}_{\rm E}^{\rm consis}} +
    \underbrace{\tau\,({\bf I}-\vm{\Pi})^\top({\bf I}-\vm{\Pi})}_{{\bf K}_{\rm
                    E}^{\rm stab}}, \label{eq:2-6}
\end{align}
where ${\bf I}$ represents the identity matrix of size
${\rm N_E}\times {\rm N_E}$, and the splitting into a consistent term
and a stabilization term replicates that in \eqref{eq:2-4a}. The
matrices appearing in \eqref{eq:2-6} are given by
$\vm{\Pi}_*={\bf G}^{-1}{\bf B}$ and $\vm{\Pi}={\bf D}\vm{\Pi}_*$,
where, in turn
\begin{align*}
    {\bf G}_{3\times 3} &= 
             \begin{bmatrix}
               1 & 0 & 0 \\ \\
               0 & \frac{|{\rm E}|}{h_{\rm E}^2} & 0 \\ \\
               0 & 0 & \frac{|{\rm E}|}{h_{\rm E}^2} 
             \end{bmatrix}, \quad
                       {\bf B}_{3\times {\rm N_E}} = 
              \begin{bmatrix}
                \frac{1}{{\rm N_E}} & \ldots & \frac{1}{\rm N_E} \\ \\
                \{ \frac{1}{h_{\rm E}} \ \ 0 \} \cdot \vm{a}_1 
              & \ldots &
              \{ \frac{1}{h_{\rm E}} \ \ 0 \} \cdot \vm{a}_{\rm N_E} 
              \\ \\
              \{ 0 \ \ \frac{1}{h_{\rm E}} \} \cdot \vm{a}_1 
              & \ldots &
              \{ 0 \ \ \frac{1}{h_{\rm E}} \} \cdot \vm{a}_{\rm N_E} 
            \end{bmatrix}, \quad
                         {\bf D}_{\rm N_E\times 3} = 
                         \begin{bmatrix}
                          m_1({\bf x}_1) & m_2({\bf x}_1) & m_3({\bf
                            x}_1) \\
                          \vdots & \vdots & \vdots & \\
                          m_1({\bf x}_{\rm N_E}) & m_2({\bf x}_{\rm
                            N_E}) & m_3({\bf x}_{\rm N_E})
                        \end{bmatrix}.
\end{align*}
The matrix $\tilde{\bf G}$ equals ${\bf G}$ but with zeros inserted
into the first row. The vectors $\vm{a}_i$, for $i=1$ to ${\rm N_E}$,
appearing in the expression for ${\bf B}$ are given by
\begin{align}
  \vm{a}_i &= \int_{\partial E} \varphi_i \vm{n} \, ds  
             = \frac{1}{2} \vm{n}_{{i-1}} |e_{i-1}| + 
             \frac{1}{2} \vm{n}_{{i-1}} |e_{i}|  
             = \frac{1}{2} \begin{bmatrix}
               y_{i+1} - y_{i-1} \\    
               x_{i+1} - x_{i-1}
             \end{bmatrix} , \notag
\end{align}
where cyclic vertex indexing is assumed, so that $\vm{x}_{\rm N_E+1}=\vm{x}_1$
and $\vm{x}_{0} = \vm{x}_{\rm N_E}$. Finally, the element load vector
follows from \eqref{eq:2-4b} as
\begin{align}
  {\bf f}_{\rm E} = |{\rm E}|\,f({\bf x}_{\rm E})\,
  \begin{bmatrix}
    1\\ \vdots \\ 1
  \end{bmatrix}_{{\rm N_E}\times 1}. \notag
\end{align}

The global stiffness matrix ${\bf K}$ and load vector ${\bf F}$ are
computed by assembling contributions from the element matrices and
vectors in the usual manner. The degrees of freedom conjugate to the
approximation space in the method are, conveniently, function values
of $u_h$ at the nodes of the mesh ${\cal T}_h$.  

The expressions in \eqref{eq:2-6}, while directly useful for numerical
computations~\cite{Sutton:2017:VEM}, necessarily hide important
insights on the method. For instance, the matrix
${\bf K}_{\rm E}^{\rm consis}$ is rank-deficient (${\rm rank}=2$),
except on triangles. The matrix $\vm{\Pi}_*$ is the representation of
the projection operator $\Pi_1$ in the scaled monomial basis. It is
necessary that ${\bf G}={\bf B}{\bf D}$, which is a useful check in
numerical implementations. Discussions of these and further aspects of
the derivations leading to the expressions for ${\bf K}_{\rm E}$ can
be found in~\cite{Beirao:2013:BPV}.  The stabilization term included
in the stiffness matrix is not essential in the VEM.  
\revised{We refer 
to~\cite{Berrone:2021:LOS} for a recent stabilization-free extension of
the VEM for the Poisson problem.} Similarly,
\cite{Chen:2023:SFV,Chen:2023:SFS,Xu:2024:SFV} define stiffness
matrices for the linear elasticity problem without the need for
stabilization.

\section{Element agglomeration}
\label{sec:3}
Element agglomeration, essentially consisting of merging adjacent
elements along a common edge, falls under the category of topological
operations for mesh improvement.  It is a familiar operation in the
literature on improving triangle and quadrilateral meshes and is used
routinely in practice~\cite{daniels2009localized,
  kinney1997cleanup}. Recent works examine its application to meshes
of polygons~\cite{Sorgente:2024:MOV,sorgente2023mesh}. Here, we
introduce a robust and purposeful generalization of the operation to
improve the performance of virtual elements.  Specifically, we propose
employing element agglomeration to improve the condition number of the
global stiffness matrix in the virtual element method.

Recall that a mesh ${\cal T}_h$,
  being a polygonal partition of a domain, satisfies the condition
that if ${\rm K}_1$ and ${\rm K}_2$ are distinct polygons in the
partition, then ${\rm K}_1\cap {\rm K}_2$ is either empty or consists
of a collection of their common vertices and edges.  If
${\rm K}_1\cap {\rm K}_2$ contains a common edge, we refer to the
corresponding virtual elements ${\rm E}_1$ and ${\rm E}_2$ as
\emph{edge-adjacent}.  Denoting the set of vertices of ${\rm K}$ by
${\sf V}({\rm K})$, we say that a pair of edge-adjacent 
virtual
elements ${\rm E}_1$ and ${\rm E}_2$ are \emph{agglomerable} if
${\sf V}({\rm K}) = {\sf V}({\rm K}_1)\cup {\sf V}({\rm K}_2)$, and
denote the set of agglomerable neighbors of ${\rm E}$ in
the mesh by ${\cal A}({\rm E})$. It follows that if
${\rm E}_2\in {\cal A}({\rm E}_1)$, then
${\rm E}_1\in {\cal A}({\rm E}_2)$ as well. This observation warrants
referring to such edge-adjacent elements ${\rm E}_1$ and ${\rm E}_2$
simply as a pair of agglomerable elements.

\begin{figure}[t]
  \subfloat[\label{fig:3-1a}]{\includegraphics[width=0.48\textwidth]{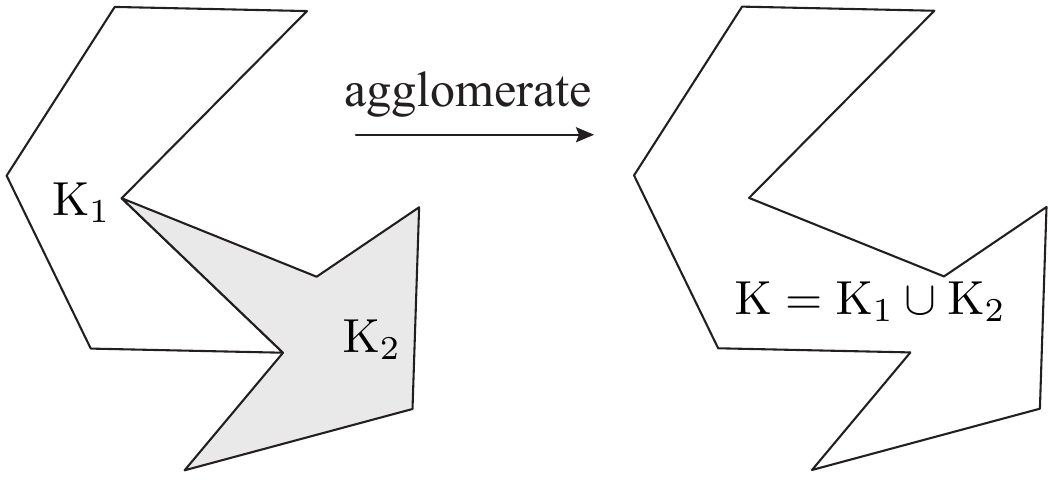}}
  \hfill
  \subfloat[\label{fig:3-1b}]{\includegraphics[width=0.48\textwidth]{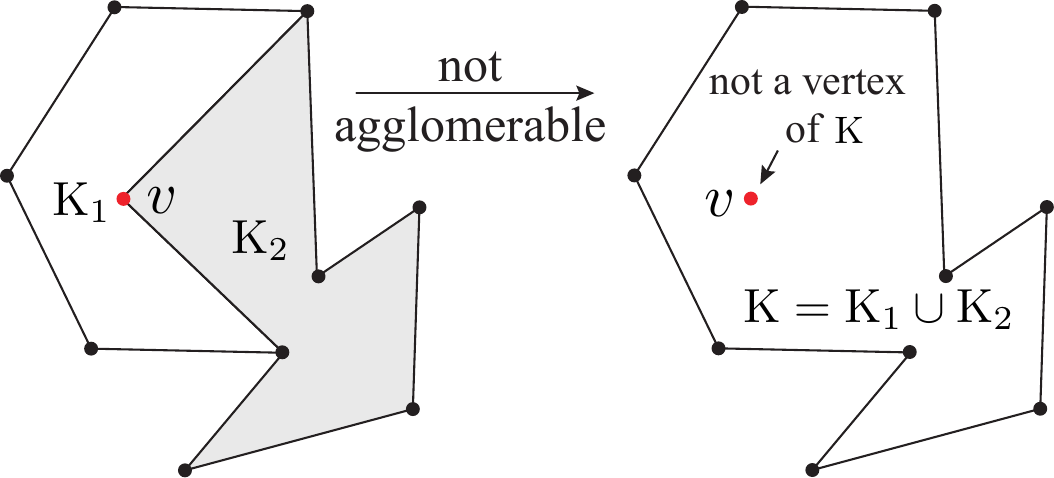}}
  \caption{The image on the left illustrates the agglomeration of
    edge-adjacent neighbors.  The image on the right depicts elements
    that are not agglomerable, since the vertex $v$ common to polygons
    ${\rm K}_1$ and ${\rm K}_2$ is not a vertex of the polygon
    ${\rm K}_1\cup{\rm K}_2$.}
  \label{fig:3-1}
\end{figure}

We can now formally define the element agglomeration operation. Given
a pair of agglomerable elements ${\rm E}_1$ and ${\rm E}_2$, their
agglomeration ${\rm E}_1\uplus {\rm E}_2$ is the virtual element
defined on the polygon ${\rm K}_1\cup{\rm K}_2$. Figure \ref{fig:3-1}
illustrates the definition with an example. Figure \ref{fig:3-1b}
depicts a scenario in which a pair of elements ${\rm E}_1$ and
${\rm E}_2$ are edge-adjacent but not agglomerable because
${\sf V}({\rm K})\subsetneq {\sf V}({\rm K}_1)\cup {\sf V}({\rm
  K}_2)$.  The condition that vertices of the union
  ${\rm K}_1\cup{\rm K}_2$ coincide with the union of vertices of
  ${\rm K}_1$ and ${\rm K}_2$ ensures that element agglomeration
  preserves the set of nodal degrees of freedom of ${\rm E}_1$ and
  ${\rm E}_2$.

\subsection{Agglomeration criterion}
\label{sec:3-1}
Define the \emph{stability ratio} of a virtual element ${\rm E}$
  as
\begin{align}
  \sigma({\rm E}) = \lambda_{\rm min}/\lambda_{\rm max}, ~\text{where}~{\rm
  Spectrum}({\bf K}_{\rm E}) = \{0<\lambda_{\rm min}\leq \ldots \leq
  \lambda_{\rm max}\}, \label{eq:3-1}
\end{align}
where ${\bf K}_{\rm E}$ is the element stiffness matrix given by~\eqref{eq:2-6}.  The zero eigenvalue in the spectrum corresponds to
the eigenvector $[1,\ldots,1]^\top \in {\mathbb R}^{N_E}$. That the remaining
eigenvalues are strictly positive is ensured by choice of the
stabilization parameter, as mentioned in \sref{sec:VEM}. As a ratio of
extreme eigenvalues, it is helpful to interpret $\sigma$ as the
inverse of the condition number of the element stiffness matrix (while
ignoring the zero eigenvalue).

Our criterion for agglomerating elements is based on the stability
ratio. Let $0<\sigma_\epsilon\ll 1$ denote a desired user-defined
lower bound for the stability ratio.  If
$\sigma({\rm E})<\sigma_\epsilon$, we identify the optimal
agglomerable neighbor of ${\rm E}$ as
\begin{align}
{\rm E_{nb}} = \left\{\tilde{\rm E}\in {\cal A}({\rm E})\,:\, \sigma({\rm
  E}\uplus \tilde{\rm E})=\max_{\hat{\rm E}\in {\cal A}({\rm
  E})}\sigma({\rm E}\uplus \hat{\rm E}) \right\}. \label{eq:3-2}
\end{align}
Then, we agglomerate ${\rm E}$ and ${\rm E_{nb}}$ if
$\sigma({\rm E}\uplus {\rm E_{nb}})>\min\{\sigma({\rm
  E}),\sigma({\rm E_{nb}}) \}$.

\begin{figure}[t]
  \centering
  \includegraphics[width=\textwidth]{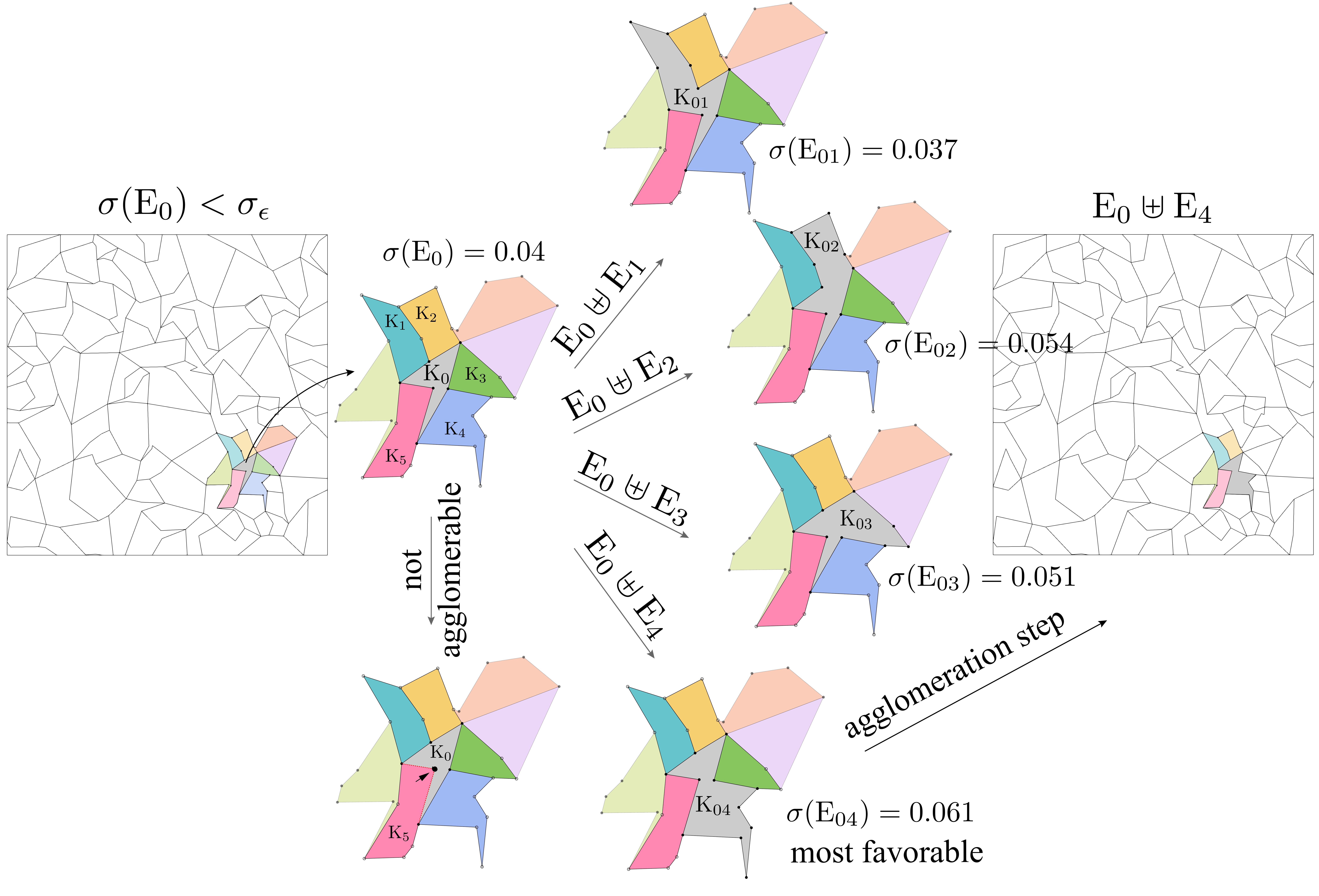}
  \caption{An example illustrating the agglomeration criterion defined
    in \eqref{eq:3-2}. The stability ratio of ${\rm E}_0$ (in gray) is
    improved from $0.04$ to $0.061$ by agglomerating it with
    ${\rm E}_4$.}
  \label{fig:3-2}
\end{figure}

Figure \ref{fig:3-2} illustrates this criterion with an example. The
element ${\rm E}_0$ in gray has a stability ratio of $0.04$. With
$\sigma_\epsilon$ set to $0.05$ (say), ${\rm E}_0$ is a candidate for
agglomeration with one of its neighbors ${\rm E}_1,\ldots,{\rm
  E}_5$. Notice that ${\rm E}_5$ is not agglomerable with ${\rm E}_0$.
Equation~\eqref{eq:3-2} thus stipulates identifying the best candidate
from among ${\rm E}_1,\ldots,{\rm E}_4$. The figure indicates the
stability ratio determined for each of the agglomerated elements
${\rm E}_{0i}\equiv{\rm E}_0\uplus {\rm E}_i$, $i=1$ to $4$. For
instance, agglomerating ${\rm E}_0$ with ${\rm E}_1$ results in an
element with a poorer ratio that ${\rm E}_0$. We find ${\rm E}_4$ to
be the best candidate for agglomeration with ${\rm E}_0$; hence
${\rm E_{nb}}={\rm E}_4$ in \eqref{eq:3-2} and the agglomerated
element's stability ratio of $0.061$ satisfies our requirement that it
be larger than $\min\{\sigma({\rm E}_0),\sigma({\rm E}_4) \}=0.04$.

\begin{breakablealgorithm}
  \caption{Optimal agglomerable neighbor}
  \begin{algorithmic}[1]
    \Function{optimal\_neighbor}{${\rm mesh}={\cal M},{\rm face}={\rm
        K},{\rm threshold}=(\sigma_\epsilon,\beta)$}

    \vspace{0.05in}
    \IfThen{$\sigma({\rm E})> \sigma_\epsilon$}{\Return ${\rm null}$}
    
    \vspace{0.05in}
    \State $({\rm K_{nb}^\star},\sigma_{\rm max})\leftarrow ({\rm
      null},\sigma({\rm E}))$ \Comment{initialization}
    
    \vspace{0.05in}
    \For{${\rm K_{nb}}\in {\cal M}.{\rm neighbors}({\rm K})$}

    \vspace{0.05in}
    \State ${\rm K_{agg}}\leftarrow {\rm K}\cup {\rm K_{nb}}$  \Comment{agglomerated polygon}

    \vspace{0.05in}
    \If{${\sf V}({\rm K_{agg}})=  {\sf V}({\rm K})\cup {\sf V}({\rm K_{nb}})$} \Comment{agglomerability}

    \vspace{0.05in}
    \If{${\cal M}.{\rm domain\_id}({\rm K})={\cal M}.{\rm domain\_id}({\rm K_{nb}})$} \Comment{subdomain}

    \vspace{0.05in}
    \If{$\sigma({\rm E_{agg}}) >\min\{\sigma_\epsilon, \beta
      \sigma({\rm E}),\beta\sigma({\rm E_{nb}})\}$} \Comment{improvement}

    \vspace{0.05in}
    \If{$\sigma({\rm E_{agg}})>\sigma_{\rm max}$} \Comment{optimality}
    \State $({\rm K_{nb}^\star},\sigma_{\rm max})\leftarrow ({\rm K_{nb}}, \sigma({\rm E_{agg}}))$
    \EndIf
    \EndIf
    \EndIf
    \EndIf
    \EndFor
    
    \hspace{-0.1in}
    \Return ${\rm K_{nb}^\star}$ \Comment{possibly ${\rm null}$}
    \EndFunction
  \end{algorithmic}
  \label{algo:3-1}
\end{breakablealgorithm}

\aref{algo:3-1} details the steps involved in implementing a
  generalization of~\eqref{eq:3-2}. Specifically, the algorithm
  incorporates two enhancements of the agglomeration criterion.
  First, in addition to the lower
bound $\sigma_\epsilon$ for the stability ratio, the algorithm
requires a parameter $\beta>1$. The parameter is used in step 8 to
decide if the agglomerated element ${\rm E_{agg}}$ represents an
improvement over the existing elements ${\rm E}$ and ${\rm
  E_{nb}}$. The criterion in \eqref{eq:3-2} only requires
$\sigma({\rm E_{agg}})>\min\{\sigma({\rm E}),\sigma({\rm E_{nb}})\}$,
which ensures that the poorer among the stability ratios of
neighboring elements ${\rm E}$ and ${\rm E_{nb}}$ improves. This,
however, can create polygons having a large number of vertices but
with only a marginally better stability ratio. The parameter $\beta$
provides a way to insist on substantial improvement in the stability
ratio with each agglomeration operation.  Note that the choice of
$\beta$ is only meaningful in scenarios where agglomeration improves
the stability ratio but not sufficiently to exceed
$\sigma_\epsilon$. The second enhancement in the algorithm
  generalizes~\eqref{eq:3-2} to permit domains with embedded
  boundaries and interfaces. Specifically, step 7 in the algorithm refines the set of elements
  agglomerable with ${\rm E}$ to only include edge-adjacent neighbors
  classified as belonging to the same subdomain in the mesh.
  
\subsection{Rationale for agglomeration}
\label{sec:3-2}
\begin{figure}[t]
  \centering
  \subfloat[\label{fig:3-3a}]{\includegraphics[height=0.35\textwidth]{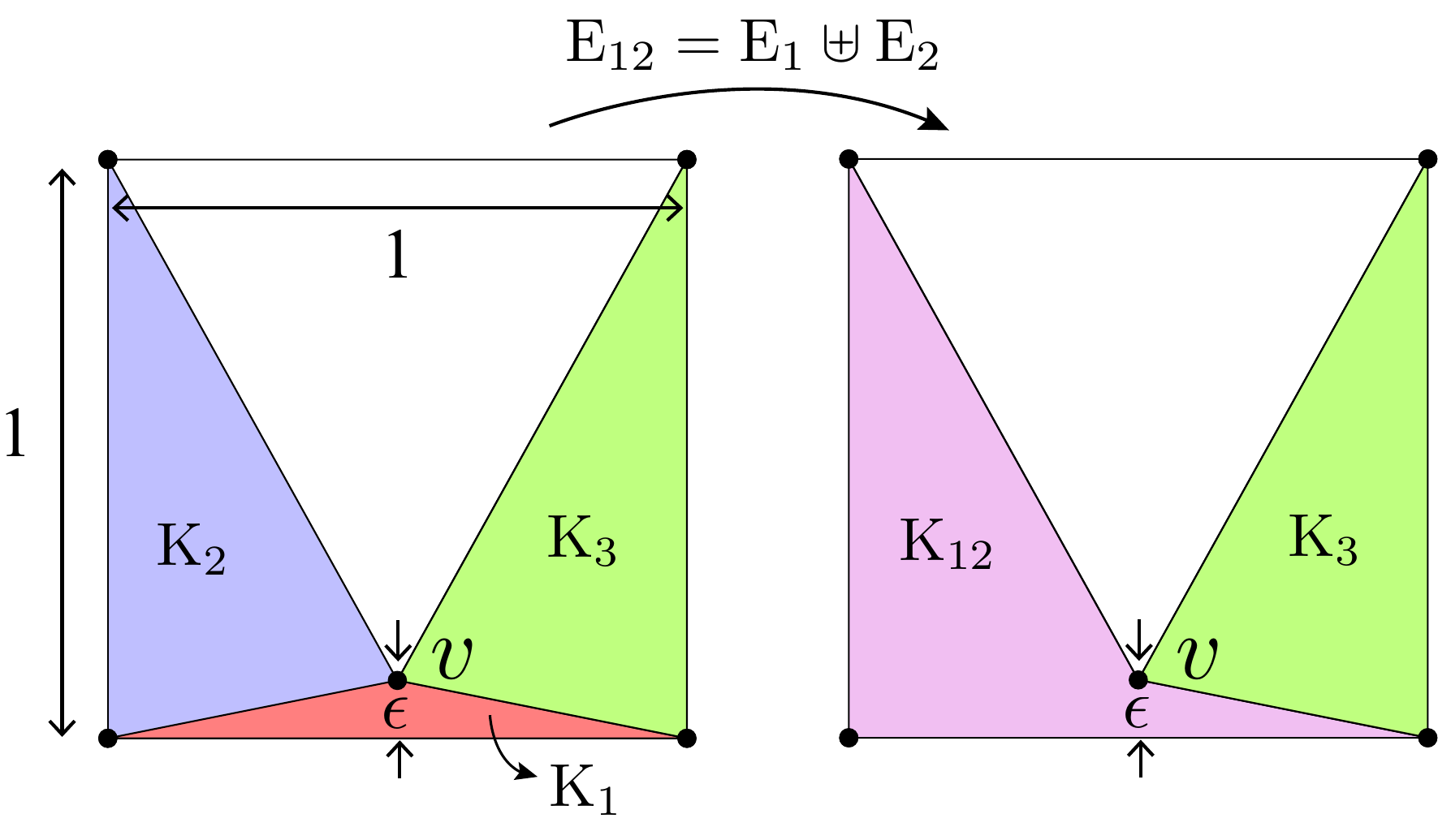}}
  \hfill
  \subfloat[\label{fig:3-3b}]{\includegraphics[height=0.29\textwidth]{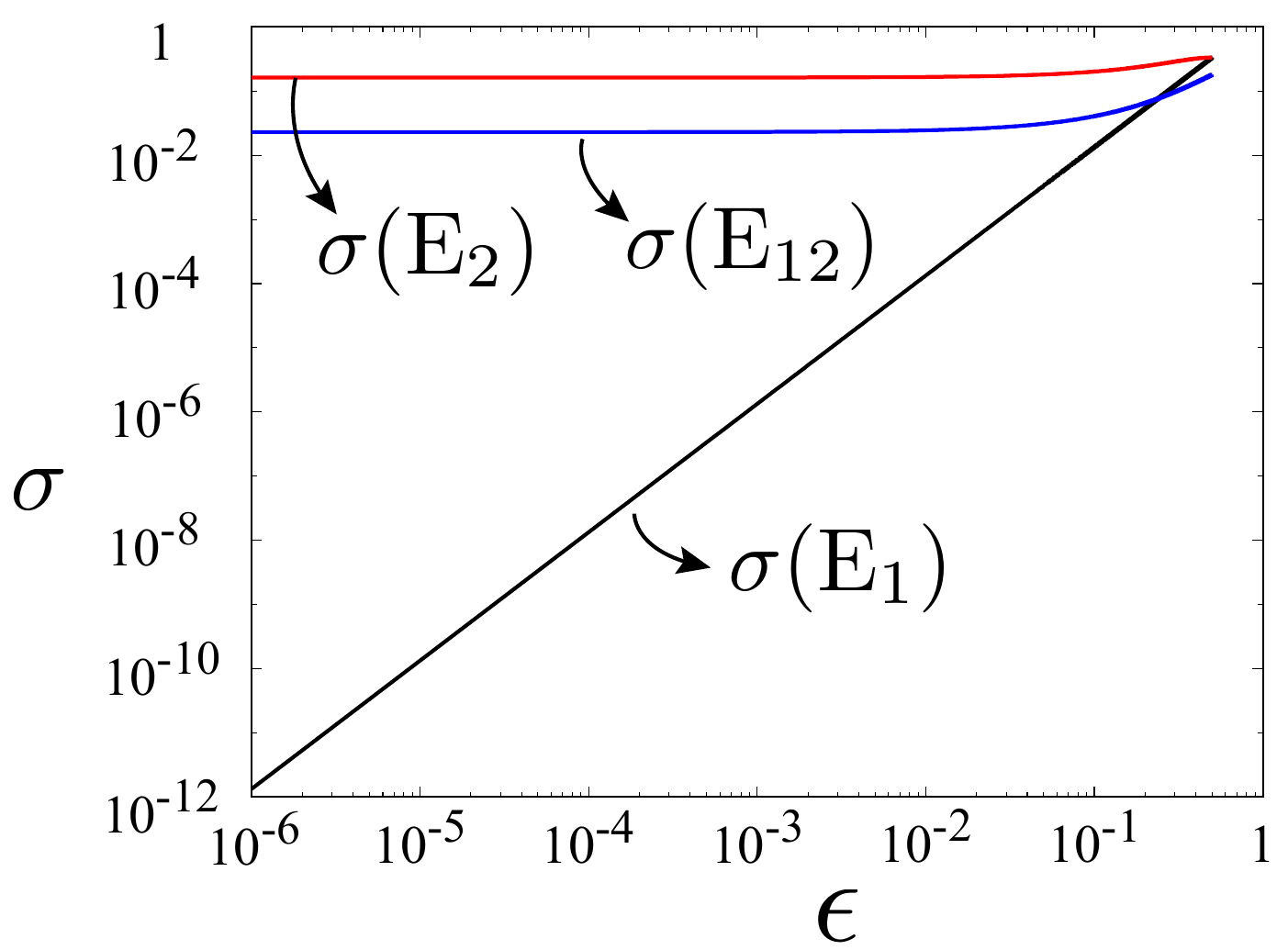}}
  \\[5pt]
  \subfloat[\label{table:3-3c}]{
    \begin{tabular}{c
       >{\columncolor{gray!15}}c
      >{\columncolor{green!10}}c
       >{\columncolor{gray!15}}c
       >{\columncolor{green!10}}c
      >{\columncolor{gray!15}}c
      >{\columncolor{green!10}}c}
      \toprule
      $\epsilon$ &
      $\lambda_{\rm min}^{\rm VEM}$ &
      $\lambda_{\rm min}^{\rm AGG}$ &
      $\lambda_{\rm max}^{\rm VEM}$ &
      $\lambda_{\rm max}^{\rm AGG}$ &
      $\lambda_{\rm max}^{\rm VEM}/\lambda_{\rm min}^{\rm VEM}$ &
      $\lambda_{\rm max}^{\rm AGG}/\lambda_{\rm min}^{\rm AGG}$           \\             
    \midrule
      $10^{-1}$ 
      & $8.7 \times 10^{-1}$  
      & $8.7 \times 10^{-1}$ 
      & $1.1 \times 10^{1}$ 
      & $8.6 \times 10^{0}$
      & $1.3 \times 10^{1}$  & $9.9 \times 10^{0}$ \\
      $10^{-2}$ 
      & $8.4 \times 10^{-1}$ & $8.3 \times 10^{-1}$ & $7.8 \times 10^{1}$ & $1.1 \times 10^{1}$
      & $9.3 \times 10^{1}$ & $1.3 \times 10^{1}$ \\
      $10^{-5}$ 
      & $8.3 \times 10^{-1}$ & $8.2 \times 10^{-1}$ & $7.5 \times 10^{4}$ & $1.1 \times 10^{1}$
      & $9.0 \times 10^{4}$ & $1.3 \times 10^{1}$ \\
      $10^{-8}$ 
      &$8.3 \times 10^{-1}$ & $8.2 \times 10^{-1}$  & $7.5 \times 10^{7}$ & $1.1 \times 10^{1}$
      & $9.0 \times 10^{7}$ & $1.3 \times 10^{1}$ \\
      \bottomrule
  \end{tabular}}
  \caption{An example illustrating the improvement in stability ratio by
  agglomerating neighboring elements. (a) Partition of the unit square
  with triangles, where the location of the vertex $v$ is set by the
  parameter $\epsilon>0$. (b) Plot showing the poor stability ratio of
  element ${\rm E}_1$ as $\epsilon\searrow 0$. In contrast, the
  agglomerated element ${\rm E}_{12}$ enjoys a superior stability
  ratio, seemingly independent of $\epsilon$. (c) Extreme
  eigenvalues of the global stiffness matrix before and after element
  agglomeration.}
  \label{fig:3-3}
\end{figure}

\begin{figure}
  \centering
  \subfloat[\label{fig:3-4a}]{\includegraphics[height=0.29\textwidth]{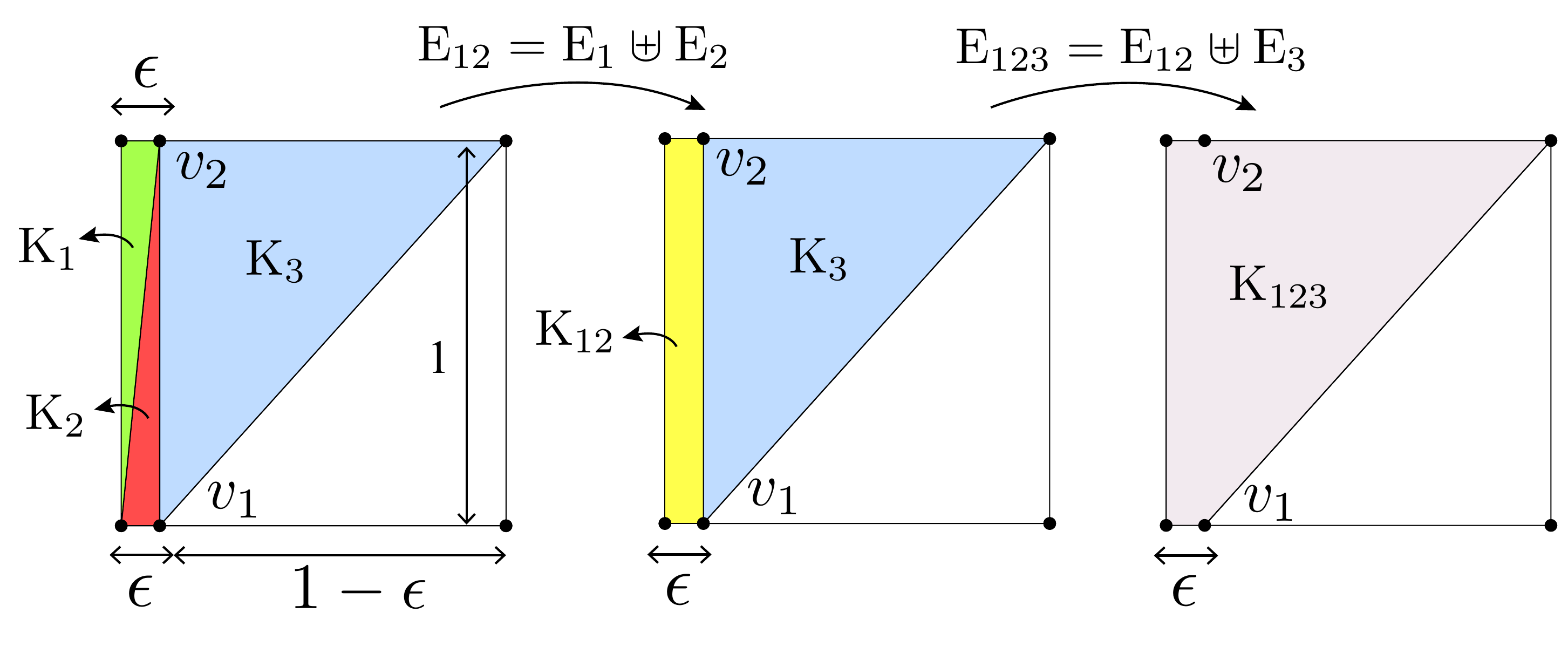}}
  \hfill
  \subfloat[\label{fig:3-4b}]{\includegraphics[height=0.24\textwidth]{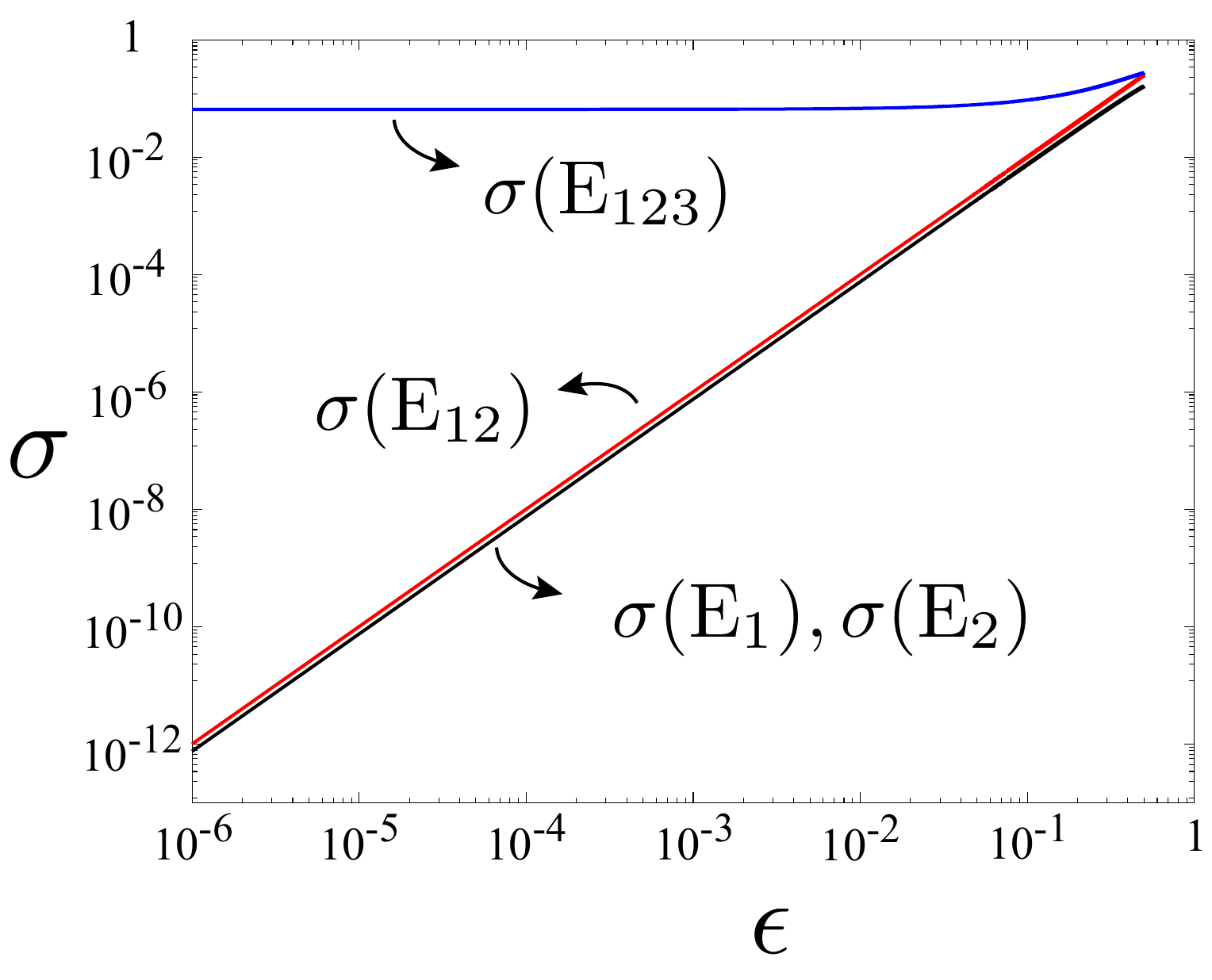}}
  \\[5pt]
  \subfloat[\label{table:3-4c}]{
    \begin{tabular}{c
    >   {\columncolor{gray!15}}c
    >   {\columncolor{green!10}}c
    >   {\columncolor{gray!15}}c
    >   {\columncolor{green!10}}c
    >   {\columncolor{gray!15}}c
    >   {\columncolor{green!10}}c}
      \toprule
      $\epsilon$ &
      $\lambda_{\rm min}^{\rm VEM}$ &
      $\lambda_{\rm min}^{\rm AGG}$ &
      $\lambda_{\rm max}^{\rm VEM}$ &
      $\lambda_{\rm max}^{\rm AGG}$ &
      $\lambda_{\rm max}^{\rm VEM}/\lambda_{\rm min}^{\rm VEM}$ &
      $\lambda_{\rm max}^{\rm AGG}/\lambda_{\rm min}^{\rm AGG}$
        \\
      \midrule
      $10^{-1}$ 
      & $6.1 \times 10^{-1}$  
      & $4.7 \times 10^{-1}$ 
      & $1.1 \times 10^{1}$ 
      & $3.8 \times 10^{0}$
      & $1.8 \times 10^{1}$ & $8.0 \times 10^{0}$ \\
      $10^{-2}$ 
      & $6.1 \times 10^{-1}$ & $4.4 \times 10^{-1}$ & $1.1 \times 10^{2}$ & $4.5 \times 10^{0}$
      & $1.8 \times 10^{2}$ & $1.0 \times 10^{1}$ \\
      $10^{-5}$ 
      &$6.1 \times 10^{-1}$ & $4.4 \times 10^{-1}$  & $1.0 \times 10^{5}$ & $4.6 \times 10^{0}$
      & $1.6 \times 10^5$ & $1.0 \times 10^1$ \\
      $10^{-8}$ 
      &$6.1 \times 10^{-1}$ & $4.4 \times 10^{-1}$  & $1.0 \times 10^{8}$ & $4.6 \times 10^{0}$
      & $1.6 \times 10^8$ & $1.0 \times 10^1$ \\
      \bottomrule
  \end{tabular}}
  \caption{An example revealing the need for iterative element
    agglomeration. Figure (a) shows that creating the rectangle-shaped
    element ${\rm E}_{12}$ by agglomerating ${\rm E}_1$ and
    ${\rm E}_2$ does not improve the stability ratio much. However,
    agglomerating ${\rm E}_{12}$ with ${\rm E}_3$ results in the
    element ${\rm E}_{123}$ having a dramatically better stability
    ratio, as shown by the plot in (b).
    (c) Extreme
    eigenvalues of the global stiffness matrix before and after element
  agglomeration.}
  \label{fig:3-4}
\end{figure}

We illustrate the significance of element agglomeration for improving
stability ratios using the examples depicted in Figs.~\ref{fig:3-3}
and~\ref{fig:3-4}. The former shows a triangular partition of a square
domain where the vertex $v$ is positioned above the bottom edge at a
height $\epsilon>0$. Letting $\epsilon$ approach zero causes the
element ${\rm E}_1$ to collapse. The element's stability ratio plotted
in~\fref{fig:3-3b} amply reflects this fact.  Such scenarios routinely
manifest in meshes with embedded interfaces when an interface splits
facets into distinct subdomains.  Figure~\ref{fig:3-3a} shows the
result of agglomerating ${\rm E}_1$ and ${\rm E}_2$. The resulting
element's stability ratio shown in~\fref{fig:3-3b} remains bounded
away from zero even for very small values of $\epsilon$.  
Figure~\ref{fig:3-3b} highlights another important feature of our
agglomeration criterion---it prioritizes improving poor elements,
even at the expense of better ones. In this example, agglomerating
elements ${\rm E}_1$ and ${\rm E}_2$ is a trade-off. The new element
${\rm E}_{12}$ is better than ${\rm E}_1$ but not ${\rm E}_2$. 

Figure~\ref{fig:3-4} furnishes a second example demonstrating the
utility of element agglomeration. Vertices $v_1$ and $v_2$ of the
triangular partition shown are positioned at a distance $\epsilon>0$
from the left edge of the square domain. Just as in the previous
examples, we find in~\fref{fig:3-4b} that elements ${\rm E_1}$ and
${\rm E}_2$ have poor stability ratios as $\epsilon$ approaches zero.
The rectangle-shaped element ${\rm E}_{12}$ resulting from
agglomerating ${\rm E}_1$ and ${\rm E}_2$ does not fare much better
either. However, agglomerating ${\rm E}_{12}$ with ${\rm E}_3$ results
in element ${\rm E}_{123}$ whose stability ratio remains bounded away
from zero.

The improvements in element stability ratios in the examples in
Figs.~\ref{fig:3-3a} and~\ref{fig:3-3b} achieved through element
agglomeration translate to better conditioned global stiffness
matrices. This is confirmed by the eigenvalues recorded in the
tables listed in Figs.~\ref{table:3-3c} and~\ref{table:3-4c}, respectively.
Noting
that the global stiffness matrix necessarily has one zero eigenvalue,
for the case in~\fref{fig:3-3a}, 
the table in~\fref{table:3-3c} lists the
smallest nonzero eigenvalues of the matrix computed before
$(\lambda_{\rm min}^{\rm VEM})$ and after
$(\lambda_{\rm min}^{\rm AGG})$ element agglomeration for
progressively smaller values of $\epsilon$. Similarly, the largest
eigenvalues are listed as $\lambda_{\rm max}^{\rm VEM}$ and
$\lambda_{\rm max}^{\rm AGG}$. The table in~\fref{table:3-4c}
lists analogous
numbers for the case in Fig.~\ref{fig:3-4a}.  Not surprisingly,
$\lambda_{\rm min}^{\rm VEM}$ and $\lambda_{\rm min}^{\rm AGG}$ are
comparable in magnitude and relatively insensitive to $\epsilon$ in
both tables. However, $\lambda_{\rm max}^{\rm VEM}$ degrades with
$\epsilon$, resulting in the stiffness matrix becoming poorly
conditioned as $\epsilon\rightarrow 0$. On the other hand,
$\lambda_{\rm max}^{\rm AGG}$ remains well-behaved and bounded away
from zero.  Notably, the improvement in $\lambda_{\rm max}^{\rm AGG}$
over $\lambda_{\rm max}^{\rm VEM}$ in both cases is despite the fact
that vertex locations is unaffected by element agglomeration. In
particular, vertex $v$ remains positioned at a distance $\epsilon$
from the base in the agglomerated elements in Fig.~\ref{fig:3-3a}, and
vertices $v_1$ and $v_2$ remain at distance $\epsilon$ from the left
edge in the agglomerated elements in Fig.~\ref{fig:3-4a}. Instead, the
improvement in the conditioning of the stiffness matrices with element
agglomeration is a direct consequence of the agglomeration criterion
seeking to improve element stability ratios.

The examples in Figs.~\ref{fig:3-3} and~\ref{fig:3-4} help highlight a
few additional details.  In both cases, the stability ratio is
improved without altering the position of any vertex.  Vertex $v$ of
${\rm K_{12}}$ remains at a distance $\epsilon$ from the bottom edge
in~\fref{fig:3-3a}; vertices $v_1$ and $v_2$ remain at a distance
$\epsilon$ from the left edge in ${\rm K_{123}}$ 
in~\fref{fig:3-4a}. The fact that elements defined over polygons with
vertices close to each other can enjoy good stability ratios is a
feature of the virtual element method. Element
agglomeration serves as a tool that redefines polygons over which to
construct the elements to exploit this feature of the VEM. The
second point is the need for the second agglomeration operation in the
example in~\fref{fig:3-4}.  This observation, in particular,
hints at the need for iterative element agglomeration.  
Algorithm~\ref{algo:3-2} presented below serves precisely this purpose.
Finally, we draw attention to triangle ${\rm K}_1$ being sliver-shaped
in~\fref{fig:3-3}, and triangles ${\rm K}_{1},{\rm K}_2$ being
needle-shaped in~\fref{fig:3-4}. Such classifications of poorly-shaped
triangles can be commonly found in the meshing 
literature~\cite{Shewchuk:WGL:2002}, and the relationship between numerical
solution accuracies and element shapes have been 
explored~\cite{babuvska1976angle, cao2005error}. For instance, the work in~\cite{babuvska1976angle} highlights slivers to be far worse than
needles in computations with linear triangle finite elements. Here,
the stability ratio naturally classifies both types of triangles
as unfavorable.

\subsection{Iterative element agglomeration}
\label{sec:3-3}
\aref{algo:3-2} below outlines the steps for iterative mesh-wide
element agglomeration.

\begin{breakablealgorithm}
  \caption{Iterative element agglomeration}
  \begin{algorithmic}[1]
    \Function{agglomerate}{${\rm mesh}={\cal M},{\rm params}=(\sigma_\epsilon,\beta),{\rm iterations}={\rm num\_iter}$}

    \For{${\rm iter}=1$ to ${\rm num\_iter}$}
    \vspace{0.05in}
    \State ${\sf pq}\leftarrow \rm{PriorityQueue}\langle\sigma\rangle (\emptyset)$  \Comment{priority queue of polygons sorted by  $\sigma$}

    \vspace{0.05in}
    \For{each polygon ${\rm K}\in {\cal M}$} \Comment{candidate faces for agglomeration}
    \IfThen{$\sigma({\rm E})<\sigma_\epsilon$}{${\sf pq}.{\rm push}({\rm K})$}
    \EndFor

    \vspace{0.05in}
    \While{${\sf pq}$ is not empty}
    \State ${\rm K} \leftarrow {\sf pq}. {\rm pop()}$    \Comment{element with poorest stability ratio}

    \State ${\rm K_{nb}}\leftarrow {\rm  optimal\_neighbor}({\cal
      M},{\rm K},(\sigma_\epsilon,\beta))$ \Comment{invoke \aref{algo:3-1}}

    \vspace{0.05in}
    \If{${\rm K_{nb}}\neq {\rm null}$} \Comment{can agglomerate}

    \State ${\cal M}.{\rm delete\_polygon}({\rm K})$ 

    \State ${\cal M}.{\rm delete\_polygon}({\rm K_{nb}})$

    \IfThen{${\rm K_{nb}}\in{\sf pq}$}{${\sf pq}.{\rm erase}({\rm K_{nb}})$}
    
    \State ${\cal M}.{\rm add\_polygon}({\rm K}\cup{\rm K_{nb}})$ 
    
    \EndIf
    \EndWhile
    \EndFor
    
    \hspace*{-0.1in} \Return ${\cal M}$
    \EndFunction
    \vspace{0.1in}
  \end{algorithmic}
  \label{algo:3-2}
\end{breakablealgorithm}

The key aspect of the algorithm consists of maintaining a list of
candidate elements for agglomeration that prioritizes improving poorer
elements.  To this end, each iteration of the algorithm begins by
populating a priority queue of elements with deficient stability
ratios (steps 3--6) sorted in increasing order of $\sigma$. Our
implementation of the algorithm uses a vector data structure with a
custom comparator for this purpose.  The iteration then examines
elements in decreasing order of priority (steps 7--17). The element
with the poorest stability ratio is popped from the queue, and its
agglomerable neighbor is identified by invoking \aref{algo:3-1}. If a
suitable neighbor is found, the mesh is updated by replacing the
element pair with their agglomeration (steps 11--14).
The numerical implementation of the algorithm is 
carried out in
\texttt{C++}, with use of The Polygon Mesh Processing
Library~\cite{pmp-library}.

Notice that once populated at the start of an iteration, the priority
queue of elements is only pruned thereafter. This is significant in
light of the possibility that the agglomerated element
${\rm E_{agg}}={\rm E}\uplus{\rm E_{nb}}$ replacing neighbors
${\rm E}$ and ${\rm E_{nb}}$ may have a stability ratio that remains
smaller than $\sigma_\epsilon$. In such an eventuality, two
possibilities arise. Either append ${\rm K_{agg}}$ into the queue for
examination during the current iteration or leave it to be examined at
the next. \aref{algo:3-2} does the latter.  In practice, this approach
prevents very poor elements from persisting at the top of the
queue. This has beneficial consequences. It helps bound the number of
edges/vertices of polygons in the mesh and prevents drastic reductions
in the number of elements due to over-agglomeration. Ensuring
mesh-wide inspections before reconsidering an element also increases
the possibility of poorer elements fortuitously benefiting from
element agglomerations in the vicinity.

\section{Element agglomeration and matrix condition numbers}
\label{sec:3-4}
The purported benefit of improving stability ratios of virtual
elements by agglomeration lies in improving the condition number of
the global stiffness matrix.  The data reported in the toy
  examples in Figs.~\ref{fig:3-3} and~\ref{fig:3-4} suggest as
  much. Here, we examine evidence of this through an extensive set of
numerical experiments. To this end, we construct a large variety
  of polygonal partitions 
  \revised{(stationary and evolving geometries)}
  over the unit square. In each case, we
assemble the global stiffness matrix for the Poisson problem with
  homogeneous boundary conditions using expressions provided 
  in~\sref{sec:VEM}, compute its condition number, and examine the range
of condition numbers realized over the ensemble of meshes with and
without element agglomeration.

\subsection{Perturbed triangulations with embedded interfaces}
\label{subsec:pert}
Our strategy to generate polygonal partitions is motivated by problems
involving domains with embedded interfaces and evolving
  geometries. To this end, consider a triangulation ${\cal T}_h$ over
the unit square, consisting of well-shaped triangles, and
characterized by a representative mesh size $h$. Let the zero level
set $\Gamma=\phi^{-1}(\{0\})$ of a level set function
$\phi:{\mathbb R}^2\rightarrow {\mathbb R}$ represent a fictitious
interface contained within the square domain. Given the pair
$({\cal T}_h,\phi)$, we introduce a sequence of meshes
$\{{\cal T}_{h,r}\}_r$, each of which differs from ${\cal T}_h$ only
in the locations of their vertices in the vicinity of $\Gamma$. The
location $\tilde{\bf x}_i$ of the $i$-th vertex in the mesh
${\cal T}_{h,r}$ is related to the corresponding vertex location
${\bf x}_i$ in ${\cal T}_h$ as
\begin{align}
  \tilde{\bf x}_i &= \begin{cases}
    {\bf x}_i + (\delta_1,\delta_2)~\text{if}~|\phi({\bf
      x}_i)|<1.25\,h, \\
    {\bf x}_i~\text{otherwise},
  \end{cases} \label{eq:3-3}
\end{align}
where the coordinate perturbations $\delta_1$ and $\delta_2$ for each
vertex are independently sampled from a uniform distribution on the
interval $[-0.15h,0.15h]$.  In effect, \eqref{eq:3-3} defines
${\cal T}_{h,r}$ as a mesh that is a localized perturbation of
${\cal T}_h$ in the vicinity of the embedded interface $\Gamma$.

Next, for each $r$, we introduce the mesh ${\cal M}_{h,r}$ that is the
result of embedding the interface
$\Gamma_{h,r} = \phi_{h,r}^{-1}(\{0\})$ in ${\cal T}_{h,r}$, where
$\phi_{h,r}$ is the continuous piecewise linear interpolant of $\phi$
defined over the triangulation ${\cal T}_{h,r}$. Unlike the
  triangulation ${\cal T}_{h,r}$, ${\cal M}_{h,r}$ is a polygonal
  mesh. Furthermore, ${\cal M}_{h,r}$ conforms to the piecewise linear
  interpolant $\Gamma_{h,r}$ of $\Gamma$. The sequence of meshes
$\{{\cal M}_{h,r}\}_r$ defined this way helps realize a large variety
of triangle-interface intersections in a manner that is representative
of scenarios that routinely occur in simulations of
moving boundary problems. Meshes with poorly-shaped triangles and
quadrilaterals, having a range of stability ratios, result naturally
in this procedure.  In particular, scenarios such as those
depicted in Figs.~\ref{fig:3-3} and~\ref{fig:3-4} manifest
as well.

The mesh sequence $\{{\cal M}_{h,r}\}_{r}$ serves as our testbed for
evaluating the efficacy of element agglomeration. An important feature
of these meshes is that {\em bad} elements can appear only in the
vicinity of the interface and not elsewhere (or everywhere). Hence,
reasonable choices for $\sigma_\epsilon$ ensure that agglomeration
operations are automatically limited to the vicinity of the interface.
This observation renders our tests and inferences directly relevant to
simulations involving embedded interfaces.

\subsection{Embedding complex interfaces}
In the first set of tests, we consider embedding varied interfaces in
the same triangulation. To this end, we choose a set of six different
level set functions $\phi_1,\ldots,\phi_6$ resulting in interfaces
having a range of feature sizes and curvatures. The negative sublevel
sets of $\phi_1,\ldots,\phi_6$ are shown in gray in \fref{fig:3-5}
alongside the triangulation ${\cal T}_h$. For each pairing of
${\cal T}_h$ with $\phi_\alpha$ ($\alpha=1$ to $6$), we generate a
sequence of ${\rm N}=1000$ realizations of the meshes
$\{{\cal M}_{h,r}^\alpha\}_{r=1}^{\rm N}$. Thus, with six different
interfaces, we generate $6000$ meshes in all. Each mesh
${\cal M}_{h,r}^\alpha$ contains triangles split by the embedded
interface in myriad ways, resulting in elements with stability ratios
expected to sample the entire range of possibilities.

A consequence of how the meshes $\{{\cal M}_{h,r}^\alpha\}_{r,\alpha}$
are generated is that they only contain triangles and
quadrilaterals. Exploiting this fact, for $1\leq r\leq 1000$ and
$1\leq \alpha\leq 6$, we assemble three global stiffness matrices
corresponding to the Poisson problem with each mesh
${\cal M}_{h,r}^\alpha$:
\begin{figure}[t]
  \centering
  \includegraphics[width=0.9\textwidth]{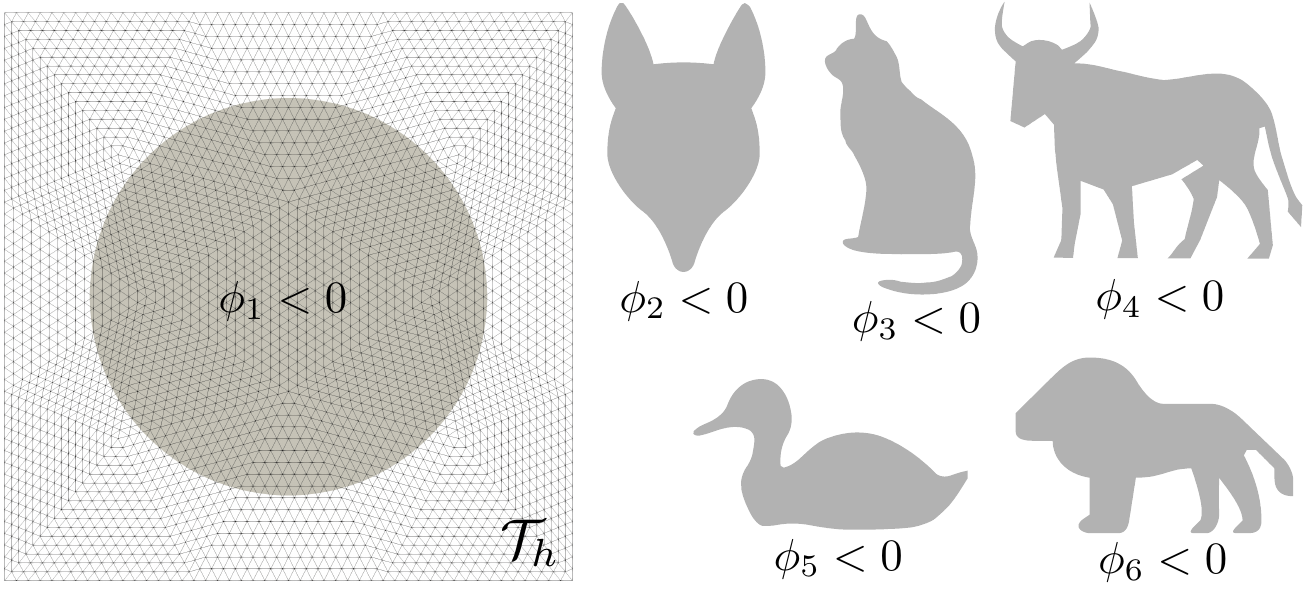}
  \caption{Level set functions corresponding to the six different
    interfaces used in the numerical experiments presented in 
    \sref{sec:3-4}.}
  \label{fig:3-5}
\end{figure}
\begin{enumerate}[(i)]
\item ${\bf K}_{\alpha,r}^{\rm fem}$ computed with isoparametric
  finite elements over each facet,
\item ${\bf K}_{\alpha,r}^{\rm vem}$ resulting from the virtual
  element method, and
\item ${\bf K}_{\alpha,r}^{\rm agg}$ just as in (ii) using the virtual
  element method, but after ${\rm num\_iter}=5$ element agglomeration
  iterations of \aref{algo:3-2} with parameters $\sigma_\epsilon$ and
  $\beta$ set to $0.2$ and $1.2$, respectively. The agglomeration
  criterion in these iterations is constrained by the location of the
  interface. That is, an element can be agglomerated only with an
  edge-adjacent neighbor lying on the same side of the interface. This
  is conveniently enforced in \aref{algo:3-1} by setting distinct
  domain identifiers to elements lying on either side of the
  interface.
\end{enumerate}
Note that ${\bf K}_{\alpha,r}^{\rm fem}$ and
${\bf K}_{\alpha,r}^{\rm vem}$ only involve assembling contributions
from triangles and quadrilaterals in ${\cal M}_{h,r}^\alpha$. Element
agglomeration, however, generates a variety of polygons, which
contribute to ${\bf K}_{\alpha,r}^{\rm agg}$.

For each $r$ and $\alpha$, we compute the condition numbers of the
three stiffness matrices
${\bf K}_{\alpha,r}^{\rm fem}, {\bf K}_{\alpha,r}^{\rm vem}, {\bf
  K}_{\alpha,r}^{\rm agg}$ as the ratios of their largest and smallest
eigenvalues. The results of these computations are compiled 
in~\fref{fig:3-6}. The figure also indicates the condition number
$\kappa_0\approx 1900$ of the finite element stiffness matrix computed
using the triangle mesh ${\cal T}_h$ itself.  It serves as our
benchmark for comparing condition numbers computed over the meshes
with embedded interfaces.  We highlight the key observations revealed by~\fref{fig:3-6}:
\begin{figure}[t]
  \centering
  \includegraphics[width=0.8\textwidth]{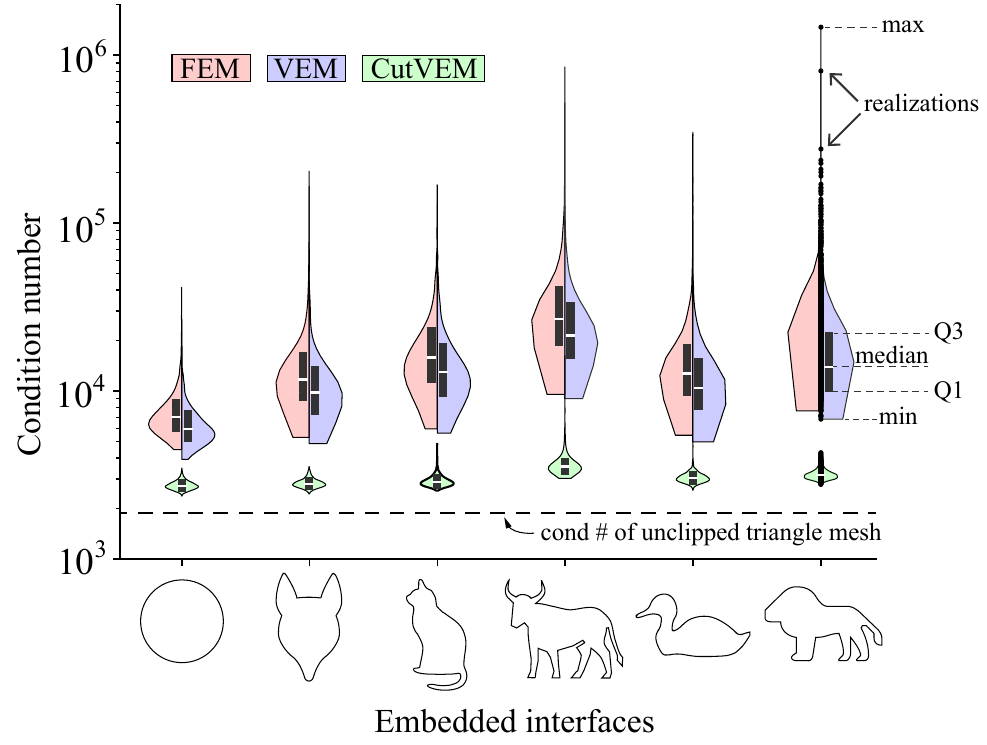}
  \caption{Violin plots summarizing the distribution of condition
    numbers of global stiffness matrices computed with interface
    embedded meshes $\{{\cal M}_{h,r}^\alpha\}_{r,\alpha}$ using the
    finite element method, the virtual element method, and the virtual
    element method after agglomeration iterations. The first and the
    third quartiles (Q1 and Q3) in each distribution are indicated in
    the plot. }
  \label{fig:3-6}
\end{figure}
\begin{itemize}
\item The condition numbers of the finite element stiffness matrices
  are significantly worse compared to $\kappa_0$. This comes as no
  surprise---a large volume of literature on cut-cell and embedded
  interface methods mentioned in \sref{sec:intro} are devoted to
  alleviating precisely this issue.

\item The condition numbers computed with the virtual element method
  fare similarly poorly. While the performance of virtual element
  methods with irregularly shaped elements is 
  well-documented~\cite{sorgente2022role}, the finding here serves as a reminder that simply
  switching from isoparametric finite elements to virtual elements
  does not alleviate the ills of embedding interfaces in meshes.

\item While the median values of the condition numbers for the finite
  and virtual element methods appear roughly an order of magnitude
  worse than $\kappa_0$, the large spread in values is significant in
  practice. The distributions show a long tail favoring larger
  condition numbers. The worst-case scenarios observed here are most
  certainly realizable in simulations involving moving interfaces
  (e.g., shape optimization problems). In particular, the extreme
  values of condition numbers are many orders of magnitude worse than
  $\kappa_0$. This is the case with each of the six interfaces
  considered, and even with the trivial circular interface, confirming
  that the observation is not an artifact of the interface shape.

\item The condition numbers of the matrices
  ${\bf K}_{\alpha,r}^{\rm agg}$ are noticeably better. The median
  values do not exceed $2\kappa_0$, and the spread in the condition
  number is significantly lower. Quite emphatically, the worst
  condition numbers found post-agglomeration are lower than the best
  case scenarios realized for ${\bf K}_{\alpha,r}^{\rm fem}$ and
  ${\bf K}_{\alpha,r}^{\rm vem}$.

\item It is important to note that the dramatic improvements in
    condition numbers indicated in the figure and evident 
    from~\fref{fig:3-6}, are achieved without any vertex 
    adjustments; in
    particular, none of the vertices created by triangle-interface
    intersections are altered during agglomeration. Rather, we
    attribute the improvement observed to the agglomeration criterion
    based on the stability ratio.
\end{itemize}

\begin{figure}[!tbh]
  \centering
  \subfloat[\label{fig:3-7a}]{\includegraphics[width=0.45\textwidth]{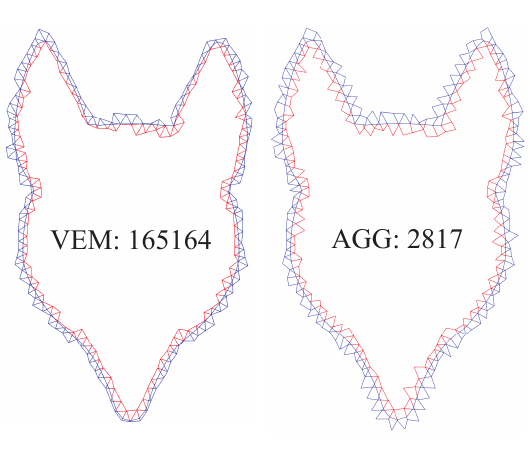}}
  \hfill
  \subfloat[\label{fig:3-7b}]{\includegraphics[width=0.52\textwidth]{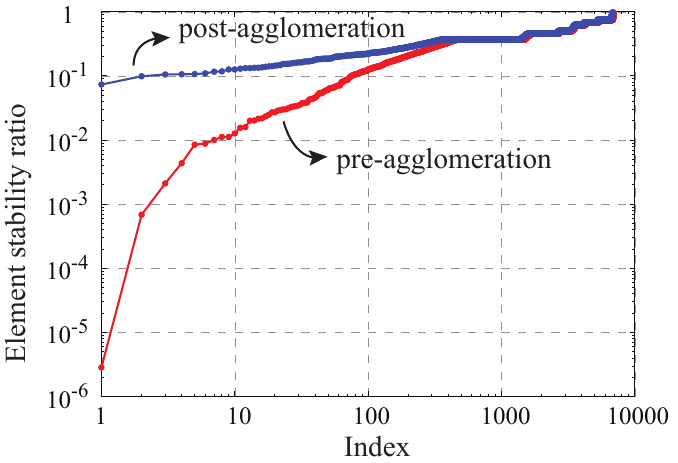}}
  \caption{A comparison of elements along embedded interfaces before
    and after agglomeration iterations is shown in (a). The condition
    numbers of the global stiffness matrix in the virtual element
    method are indicated in each case. The plot in (b) compares the
    stability ratios of elements corresponding to the particular mesh
    realization in (a). Notice the dramatic improvement in the poorest
    few elements due to agglomeration.}
  \label{fig:3-7}
\end{figure}

Figure~\ref{fig:3-7a} contrasts close-up views of elements along the
interface defined by the level set function $\phi_2$ before and after
agglomeration. Elements away from the embedded interface are all
triangles, remain unaffected by agglomeration, and are hence omitted
in the images.  For this specific case, \fref{fig:3-7b} plots the
sorted list of stability ratios of virtual elements defined over the
entire mesh before and after agglomeration. Hence, the plot provides a
direct visualization of the priority queue in~\aref{algo:3-2} at the
start and the end of agglomeration iterations. The logarithmic scale
used helps highlight that just a small fraction of virtual elements
having poor stability ratios in the mesh prior to agglomeration can
cause the global stiffness matrix in the VEM to have a poor condition
number. The plot clearly reveals the improvement in the poorest few
stability ratios resulting from agglomeration iterations. The two
datasets shown differ only over the poorest few elements that lie
along the interface, and agree over elements away from the interface
untouched by agglomeration iterations. Finally, we note that the
poorest $56$ elements after agglomeration iterations still have
stability ratios that are lower than the threshold
$\sigma_\epsilon=0.2$. Further agglomeration iterations do not result
in any improvement, reminding us that element agglomeration is not a
cure-all and does not guarantee arbitrary improvement.

\subsection{Matrix condition numbers and mesh refinement}
\begin{figure}[!tbh]
  \centering \includegraphics[width=0.8\textwidth]{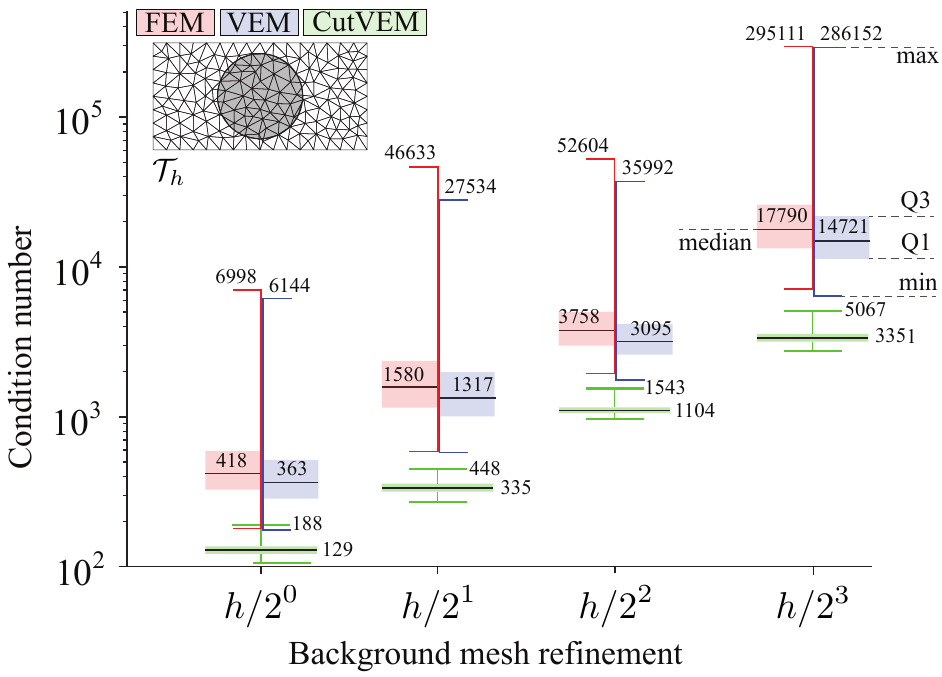}
  \caption{The interplay of element agglomeration and mesh
    refinement. By assembling stiffness matrices using finite elements
    and agglomerated virtual elements for a set of ${\rm N}=1000$
    perturbations of progressively refined meshes with an embedded
    circular interface, the plot shows that the tangible benefit of
    improved condition numbers with agglomeration persists with mesh
    refinement. }
  \label{fig:3-8}
\end{figure}

In a second set of tests, we investigate the influence of mesh
refinement on the condition numbers of stiffness matrices computed
over meshes with embedded interfaces.  We consider a circular
interface $\Gamma$ of radius $R$ embedded in self-similarly refined
triangulations ${\cal T}_{h/2^n}, n=0,1,2,3$, covering a rectangular
domain. The inset in \fref{fig:3-8} shows the mesh ${\cal T}_h$, with
$h\approx R/3$. As done previously, with each mesh ${\cal T}_{h/2^n}$,
we generate a sequence of perturbed triangulations
${\cal T}_{h/2^n,r}$ for $1\leq r\leq {\rm N}$, and in turn, the mesh
sequence ${\cal M}_{h/2^n,r}$ resulting from embedding $\Gamma$ in
${\cal T}_{h/2^n,r}$.

Our specific goal with these tests is to examine the significance of
element agglomeration for improving conditioning numbers of global
stiffness matrices defined over progressively refined meshes. To this
end, we compute ${\bf K}^{\rm fem}_{n,r}$ defined by isoparametric
finite elements, ${\bf K}^{\rm vem}_{n,r}$ by the VEM, and
${\bf K}^{\rm agg}_{n,r}$ by agglomerated virtual elements, over the
mesh ${\cal M}_{h/2^n,r}$ for each $n$ and $r$. For a given $n$, the
spread of condition numbers of
${\bf K}^{\rm fem}_{n,r}, {\bf K}^{\rm vem}_{n,r}$ and
${\bf K}^{\rm agg}_{n,r}$ for $1\leq r\leq {\rm N}$ are depicted 
in~\fref{fig:3-8}. We record the following observations from the 
plot.
\begin{itemize}
\item The median values of the condition numbers of
  ${\bf K}_{n,r}^{\rm fem}, {\bf K}^{\rm vem}_{n,r}$ and
  ${\bf K}_{n,r}^{\rm agg}$ scale roughly quadratically with the mesh
  size $h/2^n$. Hence, the condition number approximately quadruples
  with each subdivision of the mesh. This scaling is consistent with a
  priori estimates~\cite{ern2006evaluation,Mascotto:2018:ICV}. It also
  confirms that element agglomeration does not alter the scaling of
  condition numbers with the mesh size expected in the VEM.

\item Just as in~\fref{fig:3-7}, we find a large spread of condition
  numbers about the median for the stiffness matrices in the finite
  element and virtual element methods. The extreme value are, again,
  orders of magnitude worse than the median. This is a reflection of
  the fact that triangle-interface intersections can generate poor
  finite and virtual elements irrespective of mesh refinement.

\item Condition numbers computed with agglomerated virtual elements
  show marked improvement over their finite element and virtual
  element counterparts.  The median values are significantly lower;
  equally important, the spread of values is significantly smaller. In
  fact, for a given refinement level $n$, the largest condition number
  realized with agglomerated elements is comparable, if not smaller
  than the lowest value realized with the finite element and virtual
  element methods.
\end{itemize}

The observations above attest to the utility of the element
agglomeration operation in problems involving embedded interfaces, independently 
of mesh refinement. Somewhat paradoxically, reducing the
number of elements through agglomeration operations reverses mesh
refinement to improve condition numbers. Of course, it remains to see
whether solution accuracies suffer due to agglomeration. We
investigate this aspect in~\sref{sec:results}.

\subsection{Anisotropic mesh refinement}
\begin{figure}[t]
  \centering
  \subfloat[\label{fig:interpolation_mesh-a}]{\includegraphics[width=0.3\textwidth]{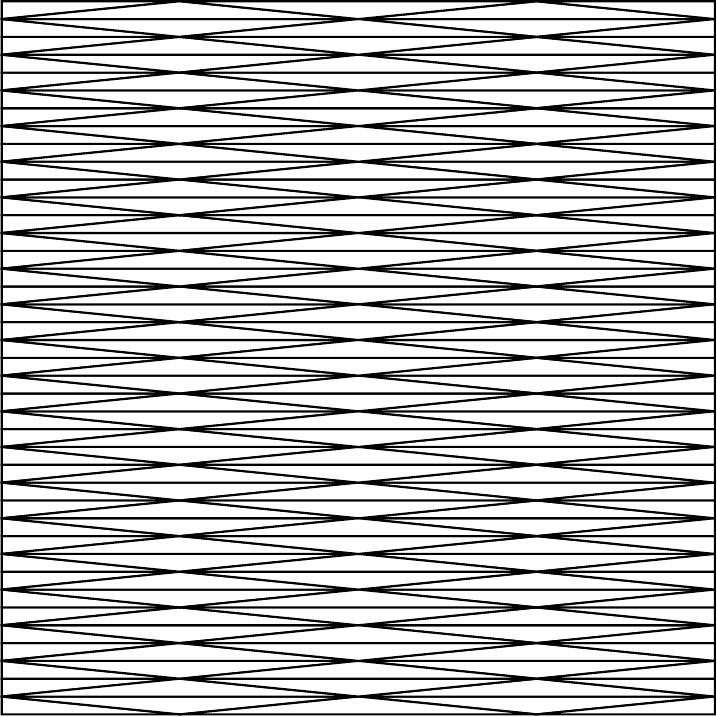}}
  \hspace{0.5in}
  \subfloat[\label{fig:interpolation_mesh-b}]{\includegraphics[width=0.3\textwidth]{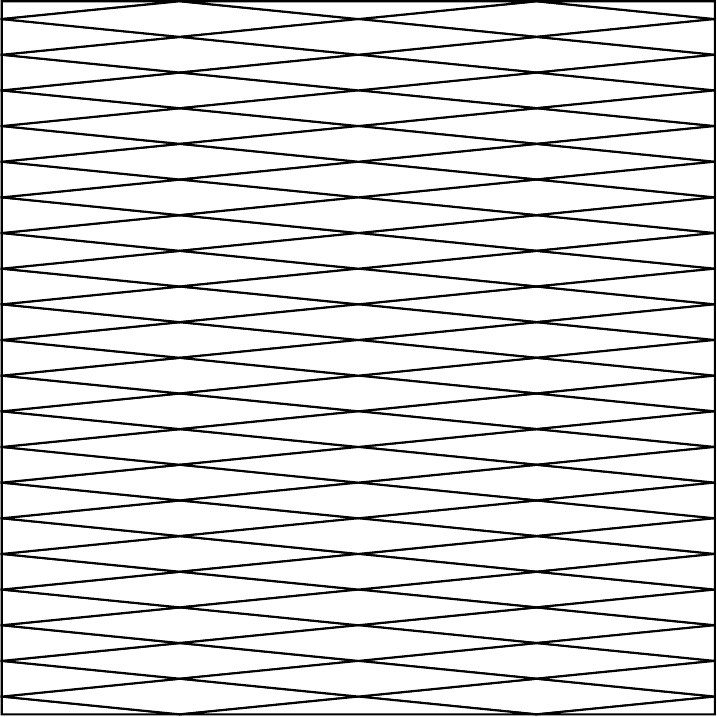}}
  \caption{(a) A triangle mesh of $144$ nodes and $200$ elements
    containing a large set of flattened triangles. (b) Iteratively
    agglomerating elements yields a reduced set of $103$ elements.
    The number of nodes remains unchanged.}
  \label{fig:interpolation_mesh}
\end{figure}

Figure~\ref{fig:interpolation_mesh} shows a generalization of the
  issue of flat triangles highlighted previously in
  Fig.~\ref{fig:3-3}. The mesh in 
  Fig.~\ref{fig:interpolation_mesh-a}
  is the result of tiling a unit square with a stencil of triangles in
  a structured albeit anisotropic manner. The ratio of the number of
  nodes along the vertical and horizontal directions in the mesh shown
  is approximately $10$. Such triangulations have been examined in the
  literature to highlight important distinctions between interpolation
  errors and finite element approximation errors~\cite{Shewchuk:WGL:2002}.

 Despite the seemingly poor aspect ratios of the triangles, the mesh
  in Fig.~\ref{fig:interpolation_mesh-a} is a Delaunay
  triangulation. Hence, the small angles in the mesh are in fact
  maximal for the given set of vertex locations. In particular, this
  implies that the shapes of triangles cannot be improved further by local
  re-triangulation. It is also evident that triangle qualities cannot
  be improved by sequentially perturbing node locations. As a
  consequence of the poor triangle qualities, the smallest nonzero
  eigenvalue of the global stiffness matrix is 
  $0.043$, while the largest
  is $61.5$, implying a condition number of
  $1430$. Figure~\ref{fig:interpolation_mesh-b} shows the result of
  iteratively agglomerating elements, resulting in a mesh consisting
  of quadrilaterals in the interior and triangles along the
  periphery. As a result of element agglomeration, the condition
  number of the stiffness matrix improves to $537$. 
  The number of
  elements reduces from $200$ to $103$, while the number of nodes
  remains unchanged.

\subsection{Stability ratios and element shapes}
\label{sec:3-5}
A distinctive aspect of our agglomeration algorithm is the criterion
for agglomeration based on the stability ratio. As the test cases
above amply demonstrate, improving stability ratios of elements by
agglomeration directly benefits the condition numbers of global
stiffness matrices in the VEM. The choice of the stability ratio over
other routinely used polygon quality metrics is deliberate. The
latter, for instance, serve as the basis for topological and geometric
mesh improvements~\cite{sorgente2023mesh,Sorgente:2024:MOV}. Here, we
highlight that the stability ratio, though possessing some features
desired of quality metrics~\cite{sorgente2023survey}, is not a polygon
quality metric in disguise.

\begin{figure}[t]
  \centering
  \subfloat[\label{fig:3-9a}]{\includegraphics[height=0.32\textwidth]{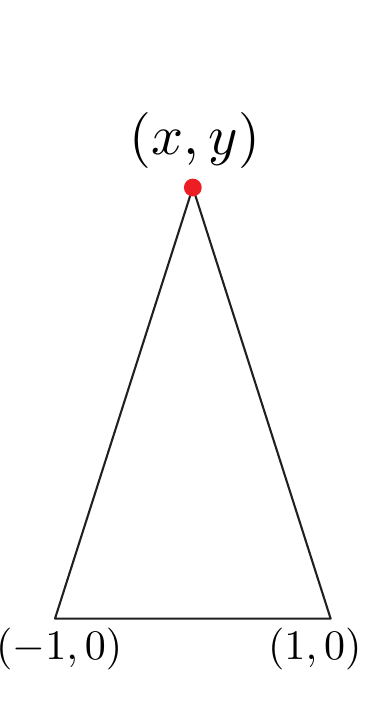}}
  \hfill
  \subfloat[\label{fig:3-9b}]{\includegraphics[height=0.32\textwidth]{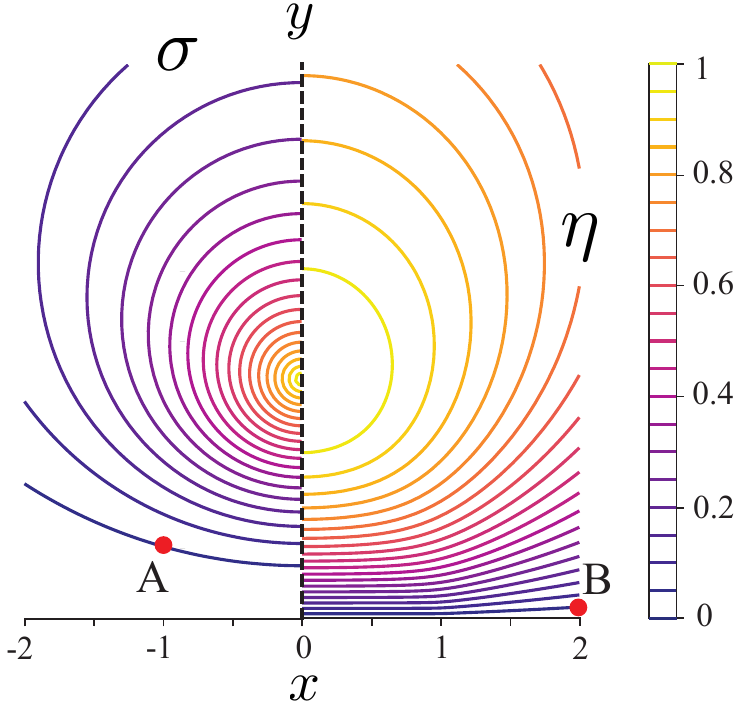}}
  \hfill
  \subfloat[\label{fig:3-9c}]{\includegraphics[height=0.32\textwidth]{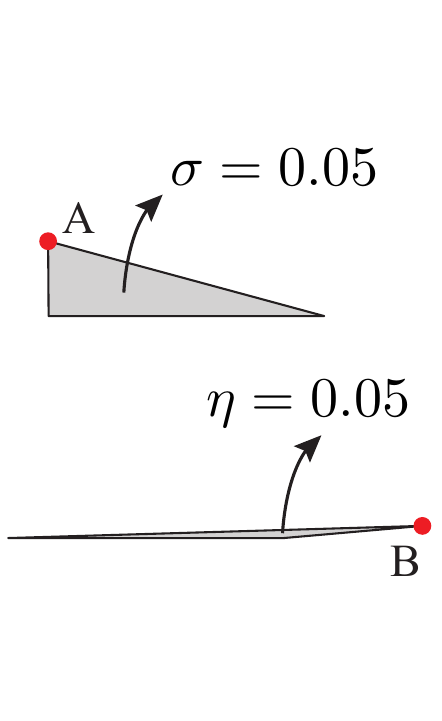}}
  \hfill
  \subfloat[\label{fig:3-9d}]{\includegraphics[height=0.32\textwidth]{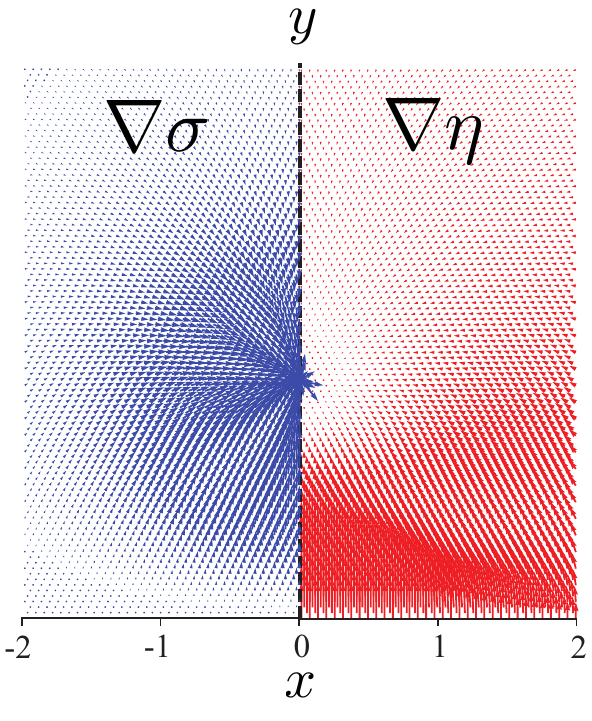}}
  \\
   \subfloat[\label{fig:3-9e}]{\includegraphics[height=0.27\textwidth]{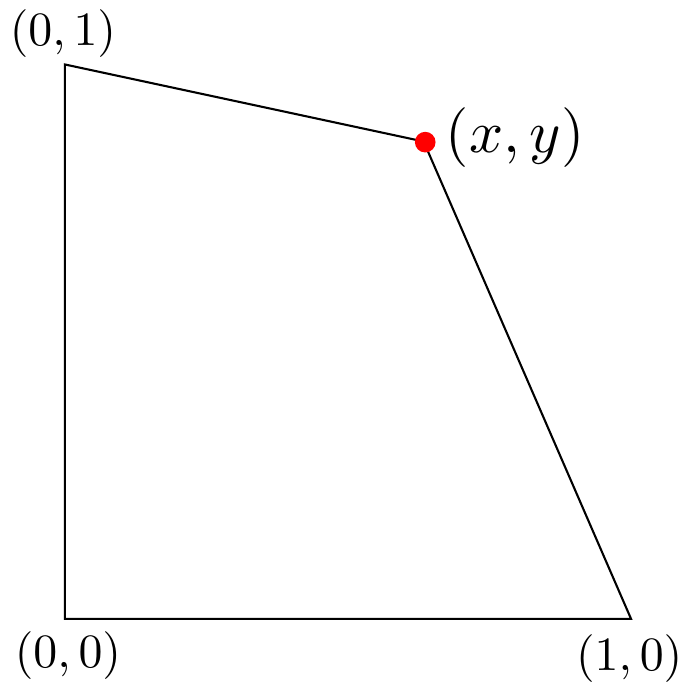}}
  \hfill
  \subfloat[\label{fig:3-9f}]{\includegraphics[height=0.27\textwidth]{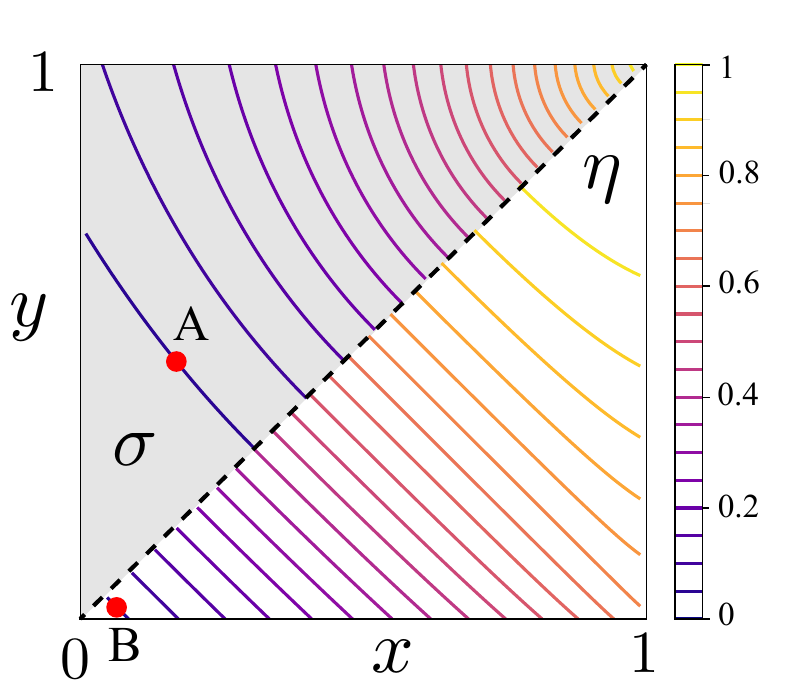}}
  \hfill
  \subfloat[\label{fig:3-9g}]{\includegraphics[height=0.27\textwidth]{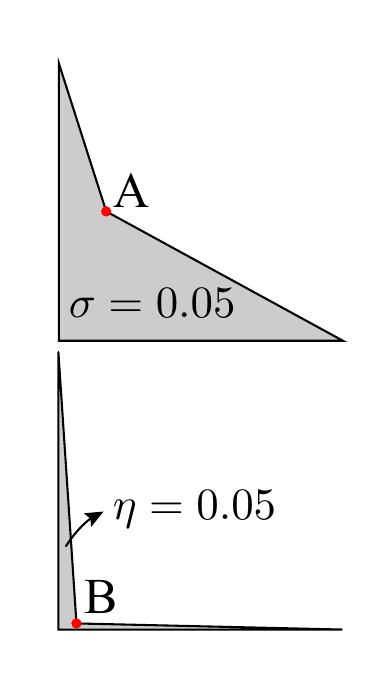}}
  \hfill
  \subfloat[\label{fig:3-9h}]{\includegraphics[height=0.27\textwidth]{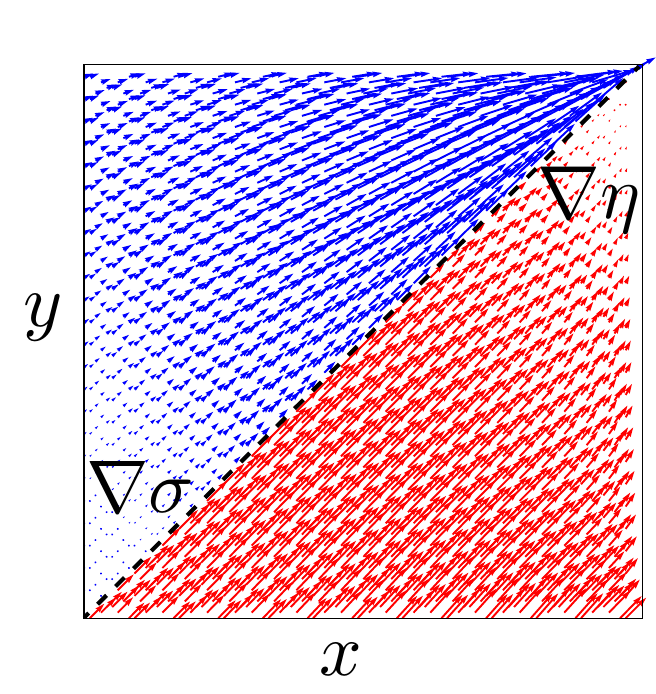}}
  \caption{Contrasting stability ratios and polygon shape qualities
    for a triangle and a quadrilateral.}
  \label{fig:3-9}
\end{figure}

For the sake of definiteness, we consider the commonly used quality
metric $\eta$ given by the ratio of the area and the squared
perimeter. Accordingly, the quality of a $n$-sided polygon ${\rm K}$
is defined as
\begin{align*}
  \eta({\rm K}) = 4n\tan\left(\frac{\pi}{n}\right) \frac{\rm Area(K)}{{\rm
  Perim}^2({\rm K})},
\end{align*}
where ${\rm Area}({\rm K})$ and ${\rm Perim}({\rm K})$ denote the area
and perimeter of ${\rm K}$, respectively. The normalizing factor
appearing in the definition ensures that a regular $n$-sided polygon
is assigned a value of $1$. We examine the distinction between
$\sigma$ and $\eta$ in the special cases of a triangle $(n=3)$ and a
quadrilateral $(n=4)$ element, see \fref{fig:3-9}. 
Figure~\ref{fig:3-9a} shows a triangle having two of its vertices fixed at
$(\pm 1, 0)$ in a Cartesian coordinate system, while its third vertex
is positioned at $(x,y)$, with $y>0$. Permitting a minor abuse of
notation, we denote the stability ratio and the quality of the
resulting triangle by $\sigma(x,y)$ and $\eta(x,y)$,
respectively. Noting both these functions enjoy even symmetry about
$x=0$, we visualize isolevels of $\sigma$ in the region $x\leq 0$ and
$\eta$ in the region $x\geq 0$ in~\fref{fig:3-9b}. In a similar vein, \fref{fig:3-9e} shows a quadrilateral with three of its vertices
located at $(0,0), (1,0)$ and $(0,1)$, and the fourth at $(x,y)$ with
$x,y>0$.  Again, we denote the stability ratio and quality of the
resulting quadrilateral by $\sigma(x,y)$ and $\eta(x,y)$. In this
case, the two maps enjoy a mirror symmetry about the line
$x=y$. Hence, we plot isolevels of the contours of $\sigma$ and 
$\eta$ in the region $x\leq y$ and $x\geq y$, respectively, 
in~\fref{fig:3-9f}.

Examining Figs.~\ref{fig:3-9b} and \ref{fig:3-9f}, the distinction
between $\sigma$ and $\eta$ stands out, especially at small
values. Figures~\ref{fig:3-9c} and~\ref{fig:3-9g} compare element
shapes corresponding to identical values of these functions, from
where we see that seemingly well-shaped elements are assigned small
stability ratios. The distinction between $\sigma$ and $\eta$ is
highlighted more clearly by visualizing their gradients in
Figs.~\ref{fig:3-9d} and \ref{fig:3-9h}. In these figures, the
(uniformly-scaled) arrows in blue and red show the gradients of
$\sigma(x,y)$ and $\eta(x,y)$, respectively. For both the triangle and
the quadrilateral element, $\sigma$ is relatively insensitive to the
location $(x,y)$ of the vertex in regions of low stability ratios and
becomes progressively more sensitive as the stability ratio
increases. In particular, elements with drastically different shapes
can have nearly identical stability ratios at small values of
$\sigma$.  The sensitivity of $\eta$ follows the opposite trend---it
is most sensitive where the quality is poorest. This observation
indirectly suggests that agglomerating elements to improve their
shapes, i.e., based on a quality metric, while seemingly reasonable,
does not correlate well with improving the stability ratio. By
extension, improving element shapes may not generally improve
condition numbers of element stiffness matrices.

\subsection{Comparison with the algorithm of 
Sorgente et al.~\protect\cite{Sorgente:2024:MOV}}
To further draw
distinction between element 
agglomeration criteria based on polygon quality metrics and 
the stability ratio, we compare the performance 
of~\aref{algo:3-2} with that in~\cite{Sorgente:2024:MOV}. 
This reference provides a set of
  five progressively refined triangle meshes, each having a significant
  fraction of poorly-shaped triangles. The agglomerated meshes
  produced by an algorithm improving a specific polygon quality metric
  are also provided.  We refer the reader to~\cite{Sorgente:2024:MOV}
  for details of the quality metric and the algorithm.  Here, we only
  note that a user-defined parameter ${\cal K}$ controls the element
  count in the agglomerated mesh computed by the algorithm proposed
  therein. For each of the five triangle meshes, 
  the agglomerated meshes are
  computed for the choices ${\cal K}=20$ and
  ${\cal K}=40$, which have 
  $20\%$ and $40\%$ of the number of elements in
  the input mesh, respectively.

\begin{figure}[h]
  \centering
  \subfloat[\label{fig:sorgente-a}]{\includegraphics[width=\textwidth]{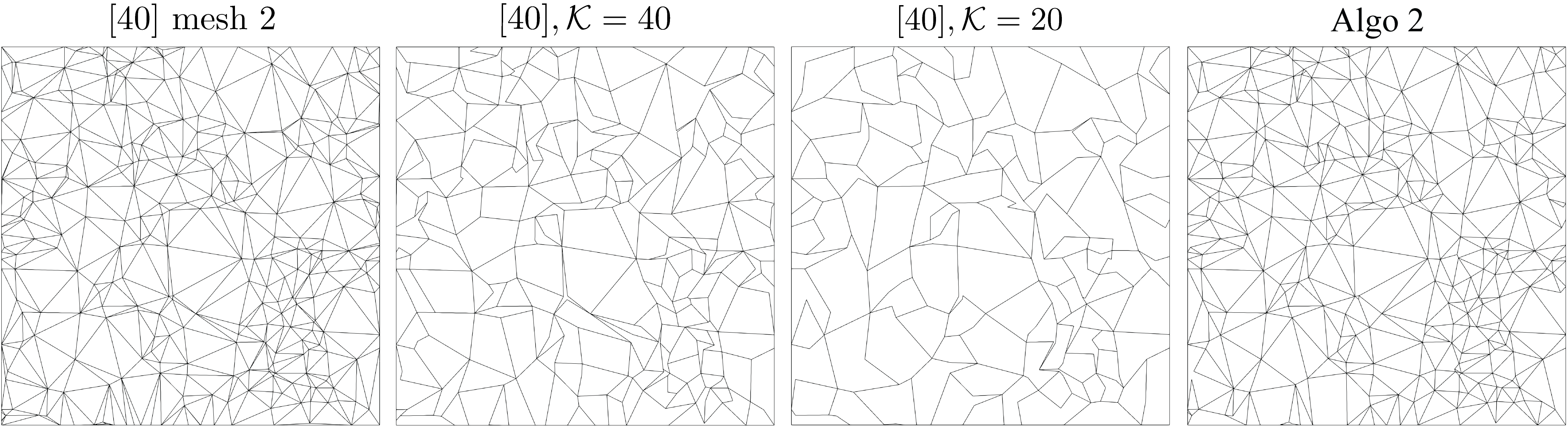}}
  \\
  \subfloat[\label{fig:sorgente-b}]{\includegraphics[height=0.31\textwidth]{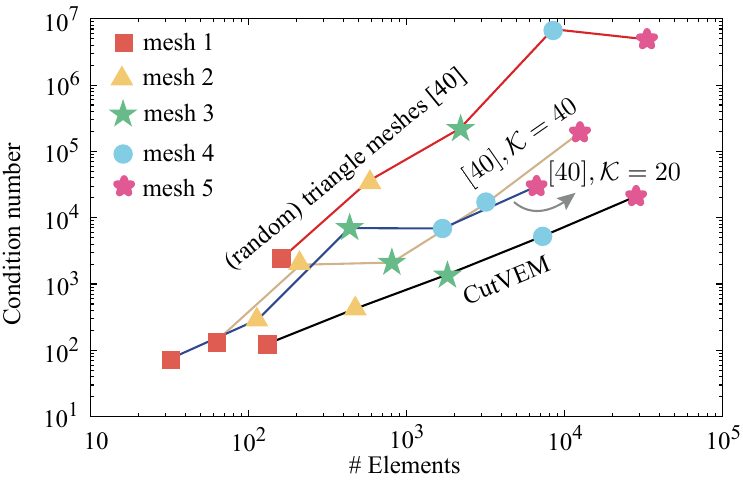}}
  \hfill
    \subfloat[\label{fig:sorgente-c}]{\includegraphics[height=0.31\textwidth]{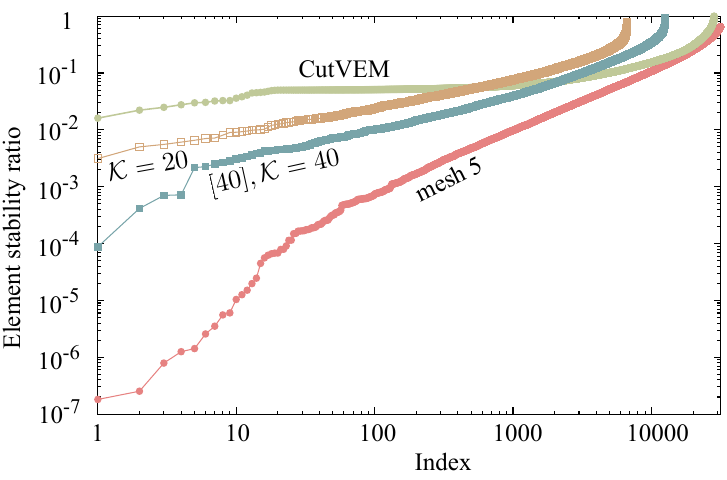}}
    \caption{A comparison of agglomeration algorithms based on polygon
      quality metrics in~\cite{Sorgente:2024:MOV} and the element
      stability ratio in~\aref{algo:3-2}. (a) Input triangle mesh
      labeled `mesh 2,' followed by the pair of meshes computed by the
      algorithm of~\cite{Sorgente:2024:MOV}, and finally the mesh
      computed by~\aref{algo:3-2}.  (b) Comparison of the condition
      numbers of global stiffness matrices. Notice the superior
      performance of \revised{the CutVEM}. This is directly supported
      by examining the element stability ratios shown in (c) for the
      case of `mesh 5.'}
      \label{fig:sorgente}
\end{figure}

Figure~\ref{fig:sorgente-a} shows the second instance in the triangle mesh sequence, which is used to generate the pair of polygonal meshes
that follow with ${\cal K} = 40, 20$~\cite{Sorgente:2024:MOV}. 
The final mesh shown in~\fref{fig:sorgente-a} is the outcome of~\aref{algo:3-2}, computed with the same set of parameters, 
$\sigma_\epsilon=0.2,\, \beta=1.2$ and ${\rm num\_iter}=5$, as employed in
all our experiments discussed thus far. A visual inspection 
of~\fref{fig:sorgente-a} suggests that elements in all three agglomerated
meshes appear better-shaped than the input triangle mesh. The
improvement in shape, however, is more pronounced in the meshes
of~\cite{Sorgente:2024:MOV}. This is to be expected, since the
algorithm employed in~\cite{Sorgente:2024:MOV} 
seeks to render polygons as regular as
possible. As the preceding discussion highlights, this is not the
case with~\aref{algo:3-2} based on the stability ratio. A second
aspect to note in the agglomerated meshes
is that the element count in
the mesh computed using~\aref{algo:3-2} is greater
than both the cases
${\cal K}=40,20$ from~\cite{Sorgente:2024:MOV}. This is a consequence
of the parameter choices in the two algorithms and should not be
considered a shortcoming (or a strength)
of either.

Figure~\ref{fig:sorgente-b} provides a quantitative comparison of the
two algorithms.  This figure plots the condition numbers of global
stiffness matrices in the virtual element method, computed just as we
did in \sref{sec:3-4}, as the ratio of the largest to the smallest
(nonzero) eigenvalues. Using distinct symbols, this plot displays data
for each of the five triangle meshes and their agglomerated
counterparts. The horizontal axis in the plot is the element count in
the mesh, which differs between the two algorithms, and for the cases
${\cal K}=40,20$ in the algorithm of~\cite{Sorgente:2024:MOV}.  In
particular, the data points for the processed meshes appear shifted to
the left of the corresponding triangle meshes since agglomeration
reduces the element count.  We highlight that with each mesh, element
agglomeration using either algorithm improves the condition number
over the input triangle mesh. This observation underlies the growing
interest in exploring agglomeration as a potent strategy to improve
the robustness of the VEM.  At least for this specific example, this
plot reveals that~\aref{algo:3-2} yields markedly better improvement
in condition numbers than the algorithm of~\cite{Sorgente:2024:MOV}.

Figure~\ref{fig:sorgente-c} examines the element stability ratios of
`mesh 5,' the most refined mesh that is considered, and its
agglomerated versions.  The figure plots the sorted list of element
stability ratios for each mesh on a logarithmic scale. We notice that
a significant fraction of elements in the input triangle mesh have
very poor stability ratios. Both agglomeration algorithms improve
these ratios, directly supporting the observation of improved
condition numbers noted in~\fref{fig:sorgente-b}. We also see
that~\aref{algo:3-2} shows more pronounced improvement in the
stability ratios, as intended. Notably, this is despite both cases
${\cal K}=40,20$ agglomerating a greater fraction of elements in the
mesh compared to~\aref{algo:3-2}. Though not the case in this example,
we caution that element stability ratios resulting
from~\aref{algo:3-2} may not be uniformly better than those of the
input mesh.  This is because the algorithm may by design deteriorate
ratios of some elements to improve the ratios of poorer ones.

\revised{
\subsection{Evolving geometry immersed in a fixed background mesh}
\begin{figure}[t]
  \centering
  \subfloat[\label{fig:moving-a}]{\includegraphics[width=\textwidth]{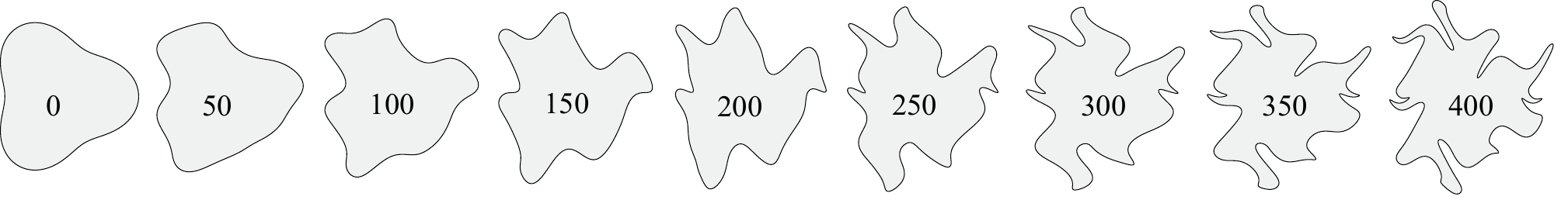}} \\
  \subfloat[\label{fig:moving-b}]{\includegraphics[height=0.24\textwidth]{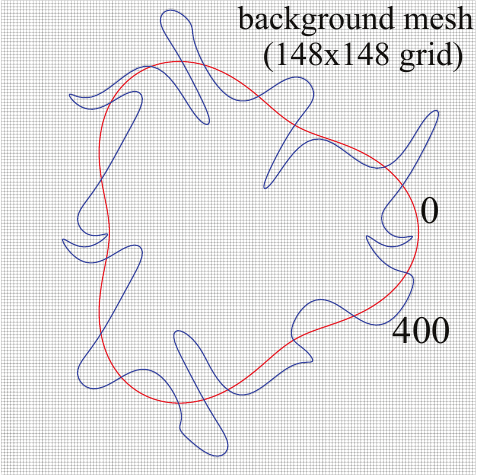}}
  \hfill
  \subfloat[\label{fig:moving-c}]{\includegraphics[height=0.24\textwidth]{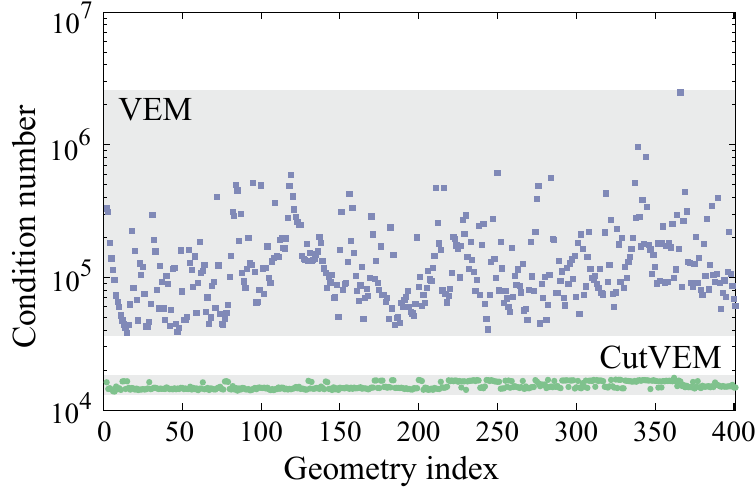}}
  \hfill
  \subfloat[\label{fig:moving-d}]{\includegraphics[height=0.24\textwidth]{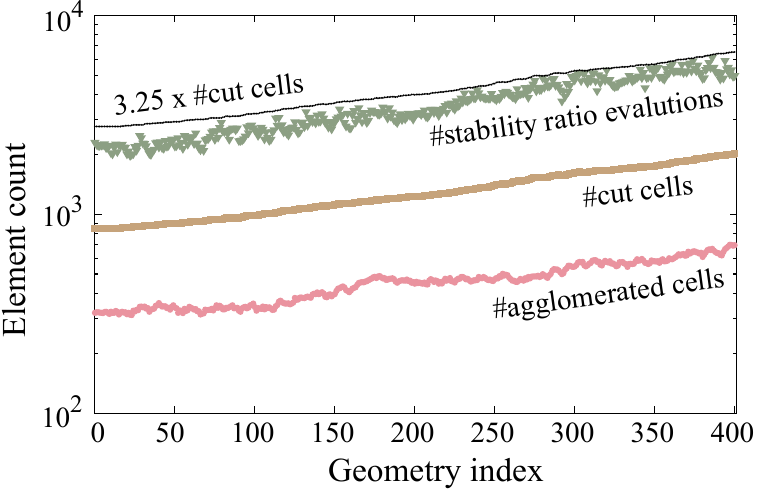}}
  \caption{(a) Shapes of an evolving interface determined by
    integrating~\eqref{eq:geom-evolve}. The step number in the time
    integration scheme, ranging from $0$ at the initial condition to
    $400$ at the final timestep, is used to index the shapes. (b) The
    interfaces computed are embedded in a fixed background mesh; the
    figure shows the mesh, and the first and the last interface
    shapes. (c) Condition numbers of stiffness matrices for the
    Poisson problem in the VEM, computed with the embedded meshes, are
    significantly worse than those computed with the 
    uncut mesh (condition number of
    9,181). We
    also observe a large spread of values. Element agglomeration in the CutVEM
    effectively resolves the issue. (d) Only a small fraction of the
    cut cells are agglomerated for each interface instance. The number
    of stability ratio evaluations, which constitutes the main
    computational effort in agglomeration iterations, is comparable to
    the number of cut cells.}
  \label{fig:moving}
\end{figure}
We conclude this section with an example representative of the
applications that CutVEM is aimed at.  We consider an evolving
interface immersed in a fixed background mesh, replicating a scenario
that arises when simulating problems of shape/topology optimization,
fluid-structure interactions, or phase transformation. To keep the
discussion self-contained, we manufacture a time-evolving interface as
the solution of the differential equation
\begin{align}
  \frac{d}{dt}\begin{bmatrix}x(s,t) \\ y(s,t)\end{bmatrix} = 
  \vm{v}
  \bigl(x(s,t),y(s,t),t \bigr) \ \textrm{for } s \in [0,1] \
  \textrm{and } t>0, \label{eq:geom-evolve}
\end{align}
with the  divergence-free velocity field $\vm{v}$ given by
\begin{align*}
  \vm{v} (x,y,t)
  = \frac{1}{2}\begin{bmatrix} \sin(2\pi y)\cos(\pi t) \\
    \cos(2\pi x) \sin(\pi t)\end{bmatrix},
\end{align*}
and an initial condition defining a three-lobed curve:
\begin{align*}
  \begin{bmatrix}
    x(s,0) \\ y(s,0)
  \end{bmatrix} = r(s) \begin{bmatrix} \cos(2\pi s) \\
    \sin(2\pi s)\end{bmatrix}, \ \textrm{where }
    r(s) = 
    1+0.2\cos(6\pi s). 
\end{align*}
We advect the curve by discretizing the spatial interval $[0,1]$ with
a uniform grid of $601$ points and marching in time using a 4th-order
Runge--Kutta algorithm with a timestep of $\Delta t = 0.00445$ for
$400$ steps. Fig.~\ref{fig:moving-a} shows the interface simulated
this way at a few instants during its evolution. The specific choice
of velocity produces smooth and oscillatory shape deformations. The
sequence of geometries thus imitates a simulation where the interface
progressively develops finer features and larger curvatures. Over the
time duration simulated, the computed piecewise linear curves
$\vx_n (s) \approx \bigl( x(s,n\Delta t),y(s,n\Delta t)\bigr)$
remain simple and connected.

We immerse the sequence of curves $\{\vx_n(s)\}_{n=0}^{400}$ in a
structured square grid ${\cal T}_h$ shown in \fref{fig:moving-b}. The
background mesh is sufficiently refined to accommodate large
curvatures realized in the geometry at later times. In particular,
every geometry instance is immersed in the same grid. Embedding the
simulated interfaces in ${\cal T}_h$ yields the sequence of meshes
$\{{\cal M}_{h,n}\}_n$ indexed by the time step. At each instant,
polygons (triangles, quadrilaterals, pentagons, and hexagons) created
by intersecting square cells in ${\cal T}_h$ with the curve are
limited to the vicinity of the interface; elsewhere, regular
square-shaped elements are retained unchanged from the background
mesh. Fig.~\ref{fig:moving-c} shows the condition numbers of stiffness
matrices for the Poisson problem computed in the VEM with these
embedded meshes. As expected, these values are significantly worse
than the condition number of $9,181$ in the case of the uncut mesh
${\cal T}_h$. Thus, we again find, as we did in the previous examples
discussed in this section, that a small fraction of cut elements in
embedded meshes can drastically worsen matrix condition numbers.  The
considerable spread observed in the data shows that the severity of
the issue is unrelated to mesh refinement. In fact, poor condition
numbers are realized even at initial times when the mesh refinement
far exceeds what is dictated by the interface's curvatures and feature
sizes. Fig.~\ref{fig:moving-c} shows that agglomerating cut elements
effectively remedies the situation. The spread in condition numbers
with the agglomerated meshes is small, and the mean condition number
of approximately $15,200$ in CutVEM compares well with that of the
uncut mesh.

Besides highlighting condition number improvement, this example helps
examine the computational effort involved in mesh
agglomeration. Fig.~\ref{fig:moving-d} plots the number of cut cells
in the embedded meshes. As evident from \fref{fig:moving-a}, the
length of the interface increases with time. The number of cut cells
in ${\cal M}_{h,n}$ therefore increases with $n$. The fraction of
these cut cells that are agglomerated depends crucially on the choice
of the lower bound for the stability ratio $\sigma_\epsilon$, which we
retain unchanged at $0.2$ from our previous examples. The figure shows
the number of successful agglomerations over five iterations in this
example. Notice that it is a small fraction of the number of cut
cells. This is because only cut cells with a stability ratio falling
below $\sigma_\epsilon$ require agglomeration. Furthermore, not all
such cells are successfully agglomerated. Now, the main computation
involved in these element agglomeration operations consists of
evaluating element stability ratios, which in turn, requires
assembling element stiffness matrices and computing their
eigenvalues. These calculations are potentially expensive, at least
when compared to evaluating simple geometric quality measures, such as
the shape quality metric considered in
\sref{sec:3-5}. \fref{fig:moving-d} plots the number of stability
ratio evaluations over five agglomeration iterations required in this
example. The count includes an initial evaluation of stability ratios
of all cut cells to mark those requiring agglomeration. Subsequently,
determining the optimal neighbor to agglomerate each marked cell
requires iterating over its agglomerable neighbors and evaluating
stability ratios of the resulting candidate agglomerated element. The
figure shows that the total number of stability ratio evaluations does
not exceed $3.25$ times the number of cut cells in the embedded
mesh. Of course, this bound is specific to the threshold ratio and the
geometry considered in this example. Nevertheless, it confirms our
expectation that the number of matrix-eigenvalue evaluations are
comparable to the number of cut cells in the mesh for reasonable
choices of $\sigma_\epsilon$.  } 

\section{Numerical experiments}\label{sec:results}
We present numerical experiments to establish the accuracy and
robustness of CutVEM. First, we assess the benefits of agglomeration
in reducing interpolation and finite element approximation errors over
anisotropically refined Delaunay meshes that consist of flat
triangles.  Then we solve a homogeneous heat conduction problem over a
clipped geometry, and use both Dirichlet and mixed (Dirichlet and
Neumann) boundary conditions.  Finally, we solve inhomogeneous heat
conduction problems in immersed (external and internal discs)
geometries: material with a void (hole) and an inclusion (contrasts in
thermal conductivity coefficient) are considered. We show that CutVEM
delivers optimal rates of convergence on agglomerated meshes.

\subsection{Assessing interpolation and finite element approximation errors over meshes with flat triangles}\label{subsec:agg_tests}
\begin{figure}[!bht]
\centering
\begin{subfigure}{0.24\textwidth}
\includegraphics[width=\textwidth]{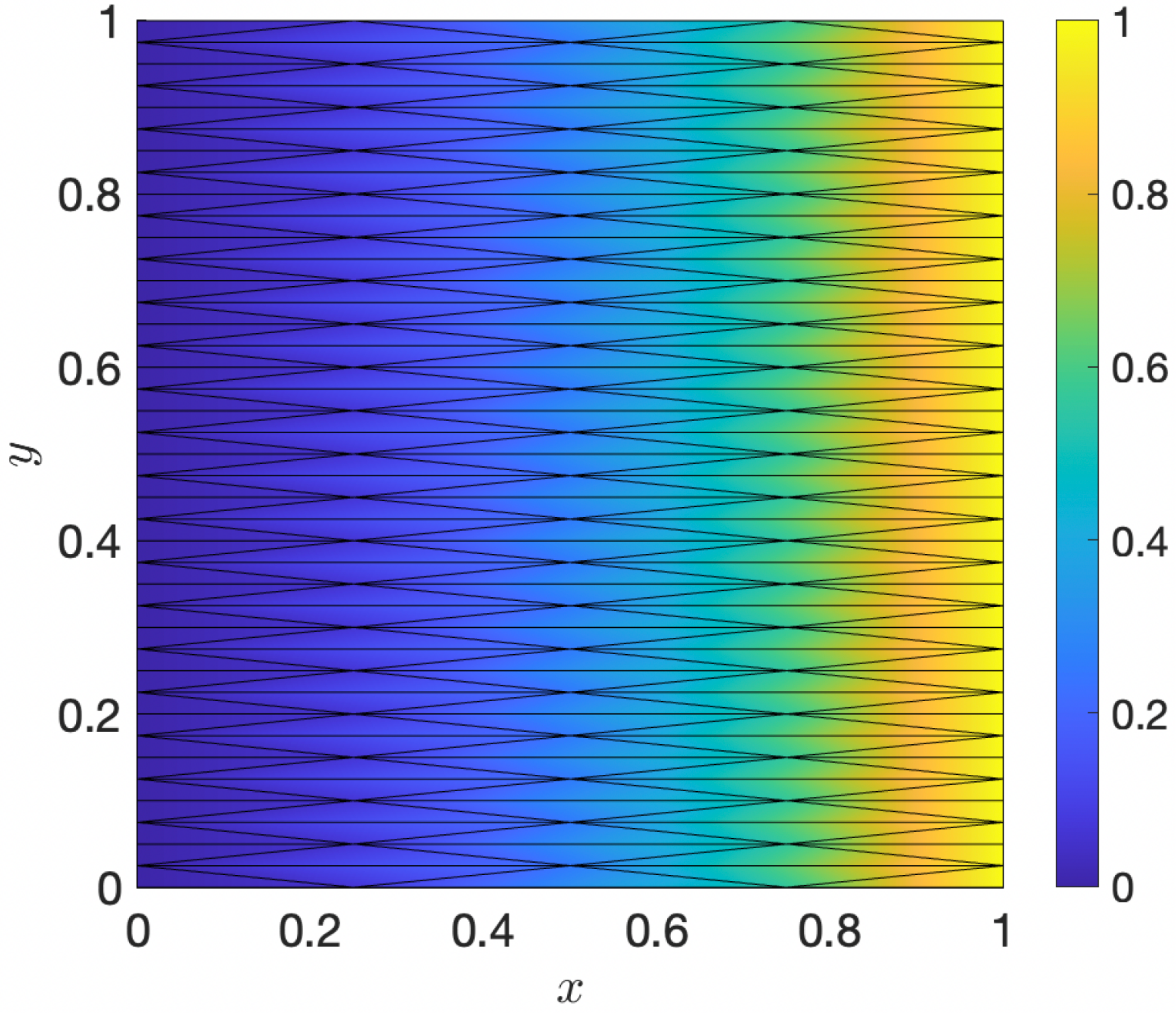}
\subcaption{}\label{fig:udu_interpolation_a}
\end{subfigure} \hfill
\begin{subfigure}{0.235\textwidth}
\includegraphics[width=\textwidth]{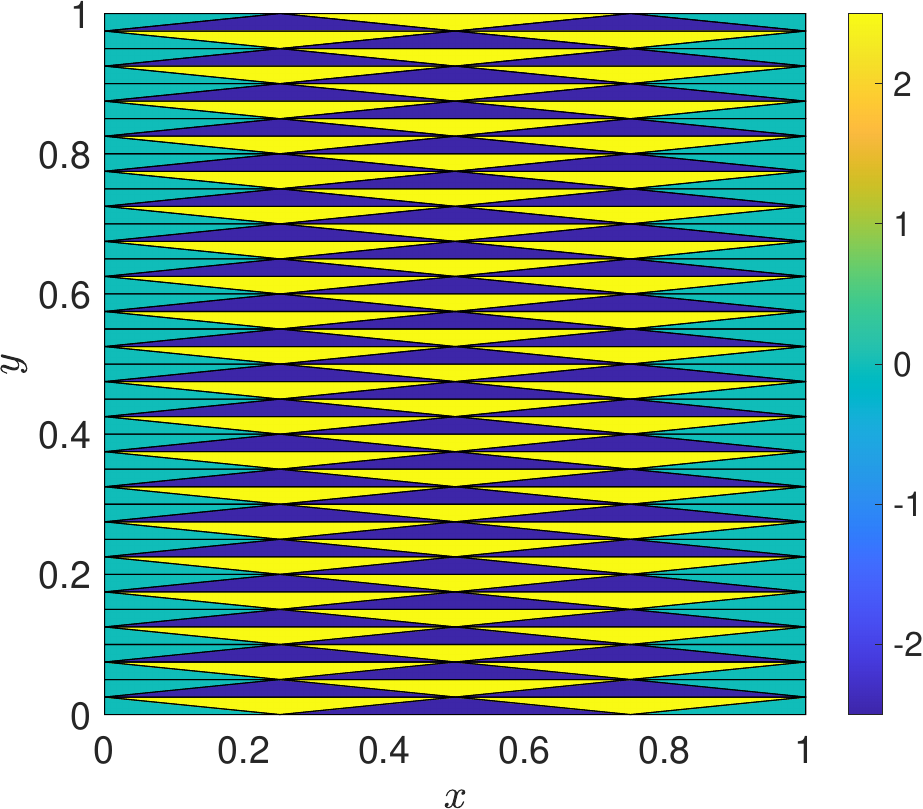}
\subcaption{}\label{fig:udu_interpolation_b}
\end{subfigure} \hfill
\begin{subfigure}{0.24\textwidth}
\includegraphics[width=\textwidth]{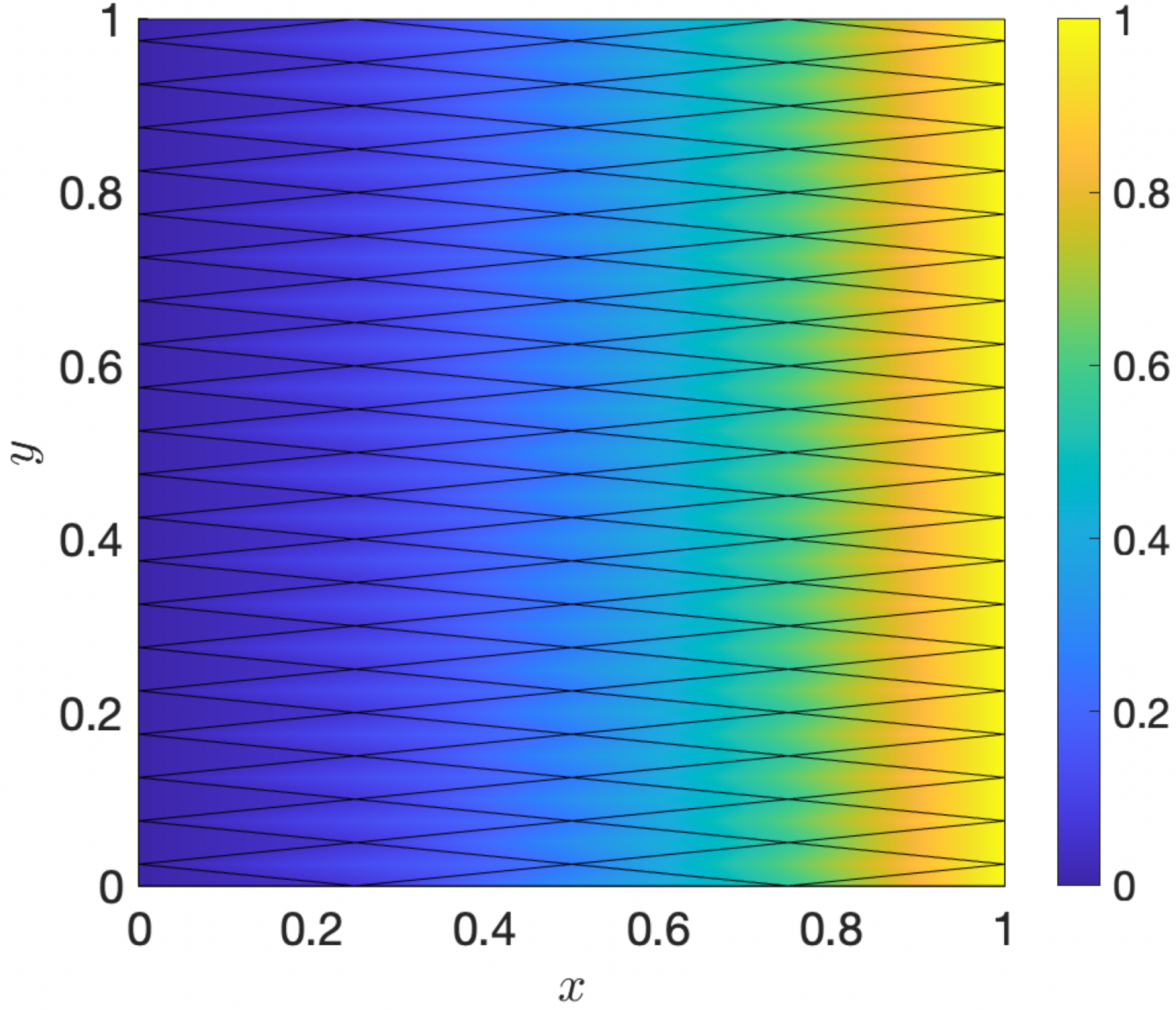}
\subcaption{}\label{fig:udu_interpolation_c}
\end{subfigure} \hfill
\begin{subfigure}{0.24\textwidth}
\includegraphics[width=\textwidth]{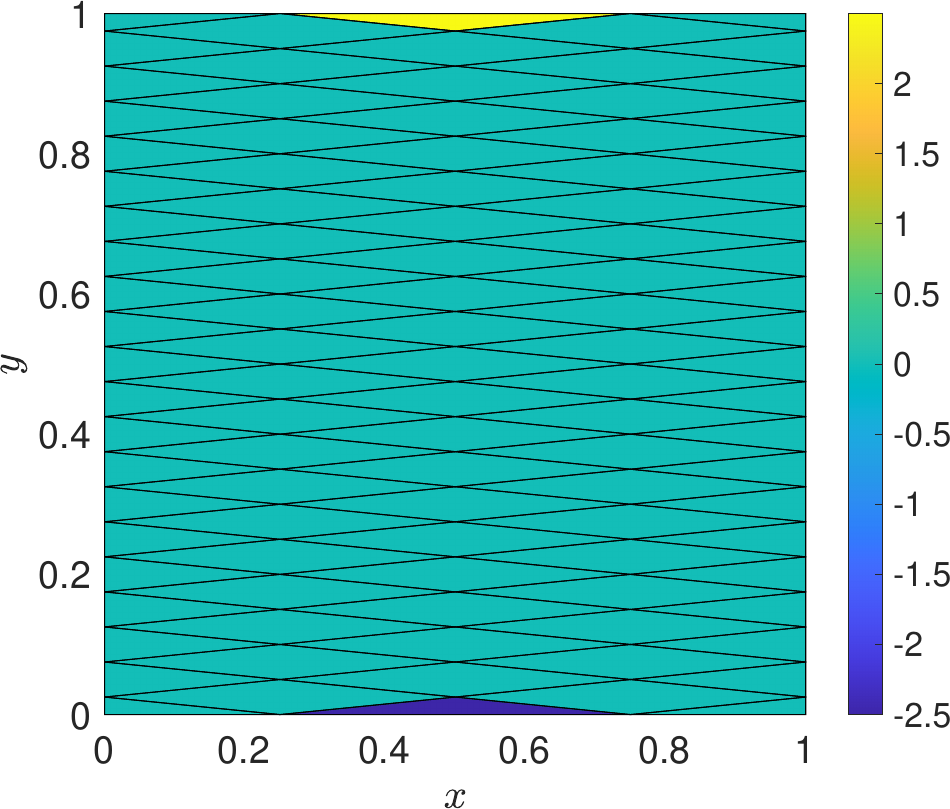}
\subcaption{}\label{fig:udu_interpolation_d}
\end{subfigure}
\caption{Interpolation of the function $u(\vx) = x^2$ on a mesh
    with flat triangles. Contour plots of the FE interpolant over the
    mesh in~\protect\fref{fig:interpolation_mesh-a} are shown for $u$
    and $\partial u/\partial y$ in (a) and (b), respectively.
    Corresponding plots for the interpolant in the virtual element
    method over the agglomerated mesh
    in~\protect\fref{fig:interpolation_mesh-b} are shown in (c) and
    (d). Note the larger errors in $\partial u/\partial y$ (exact
    solution is zero) for the FE interpolant in comparison to the
    virtual element counterpart on the agglomerated
    mesh.}\label{fig:udu_interpolation}
\end{figure}

We reconsider the mesh shown in~\fref{fig:interpolation_mesh}. On such
anisotropically refined Delaunay meshes with flat triangles, it is
known that the interpolation error of the gradient can grow without
bound~\cite{Shewchuk:WGL:2002}. To illustrate this we compute the
interpolated gradient using the FEM on the original Delaunay mesh and the VEM on an agglomerated mesh. The FE mesh
in~\fref{fig:interpolation_mesh-a} has 144 nodes and 200 elements,
whereas the agglomerated mesh in~\fref{fig:interpolation_mesh-b} has
144 nodes and 103 elements.  Contour plots of $u$ using the FEM and
the CutVEM are presented in Figs.~\ref{fig:udu_interpolation_a}
and~\ref{fig:udu_interpolation_c}, respectively, whereas those for
$\partial u/ \partial y$ appear in Figs.~\ref{fig:udu_interpolation_b}
and~\ref{fig:udu_interpolation_d}, respectively.  Since the exact
solution for $\partial u / \partial y$ is zero, we observe that the FE
solution in~\fref{fig:udu_interpolation_b} has large errors that are
concentrated over the flat triangles.
In~\fref{fig:udu_interpolation_d}, the errors in
$\partial u/ \partial y$ for CutVEM are concentrated over two
flat triangles. 

To assess finite element approximation errors and convergence, 
we consider a mesh sequence based on the mesh
shown 
in~\fref{fig:interpolation_mesh-a}, and conduct a convergence study for the following Poisson 
boundary-value problem:
\begin{subequations}
  \label{eq:sinsin}
\begin{align}
-\nabla^2 u &= 2 \pi^2 \sin ( \pi x)  \sin (\pi y) \ \ \textrm{in } \Omega = (0,1)^2 , \\
u &= 0 
    \ \ \textrm{on } \partial \Omega,
\end{align}
\end{subequations}
which has the exact solution $u(\vx) = \sin (\pi x) \sin (\pi y)$.

\begin{figure}[!tbh]
  \centering
  \subfloat[$4\times 4$\label{fig:meshes_poisson_a}]{\includegraphics[width=0.29\textwidth]{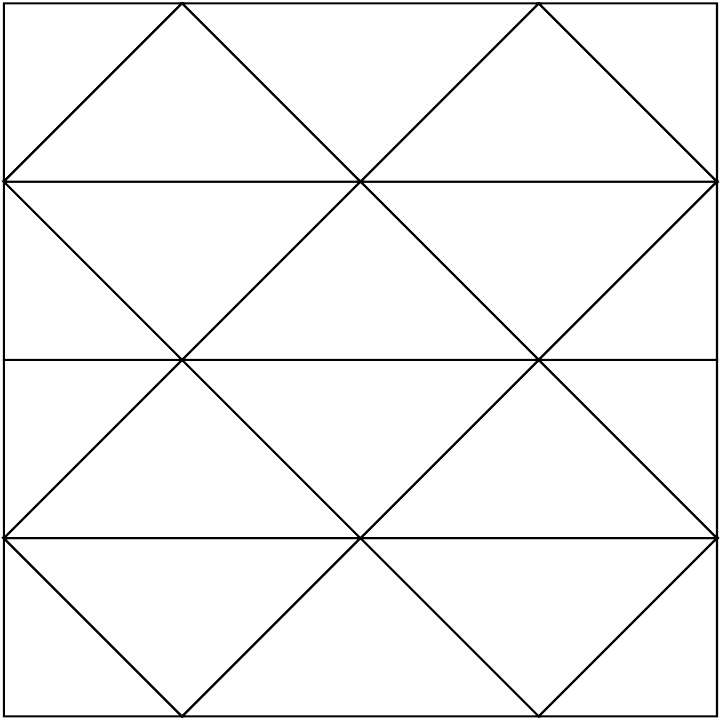}}
  \hfill
  \subfloat[$8\times 8$\label{fig:meshes_poisson_b}]{\includegraphics[width=0.29\textwidth]{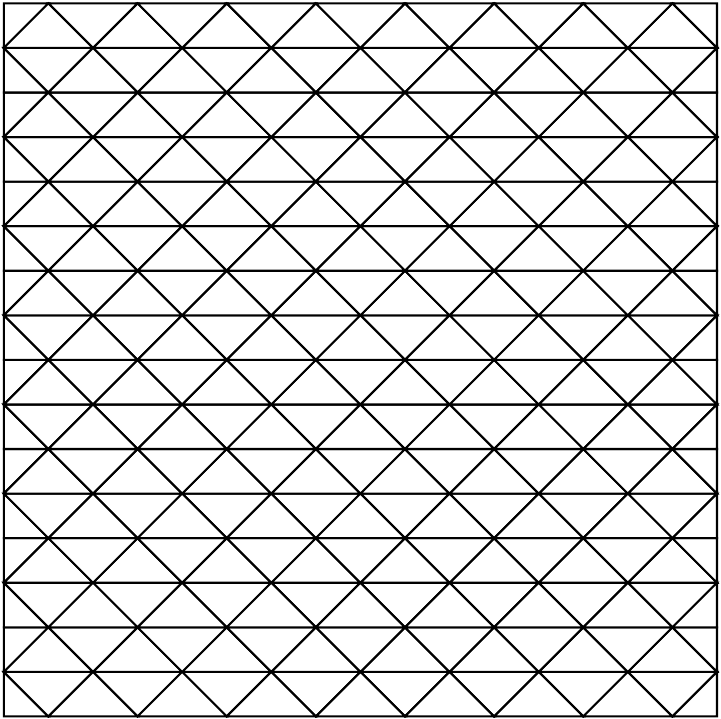}}
  \hfill
  \subfloat[Optimal convergence\label{fig:rates_poisson_fem_a}]{\includegraphics[width=0.38\textwidth]{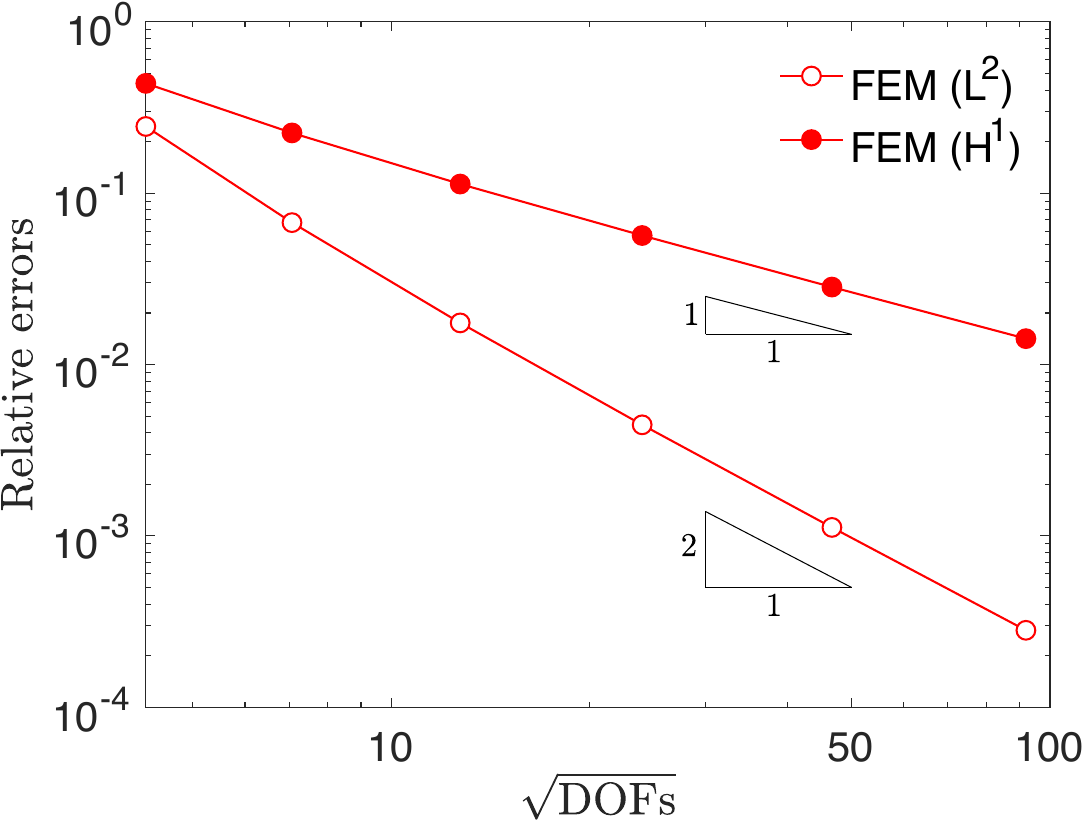}}
  \\[5pt]
  \subfloat[$2\times 10$\label{fig:meshes_poisson_d}]{\includegraphics[width=0.29\textwidth]{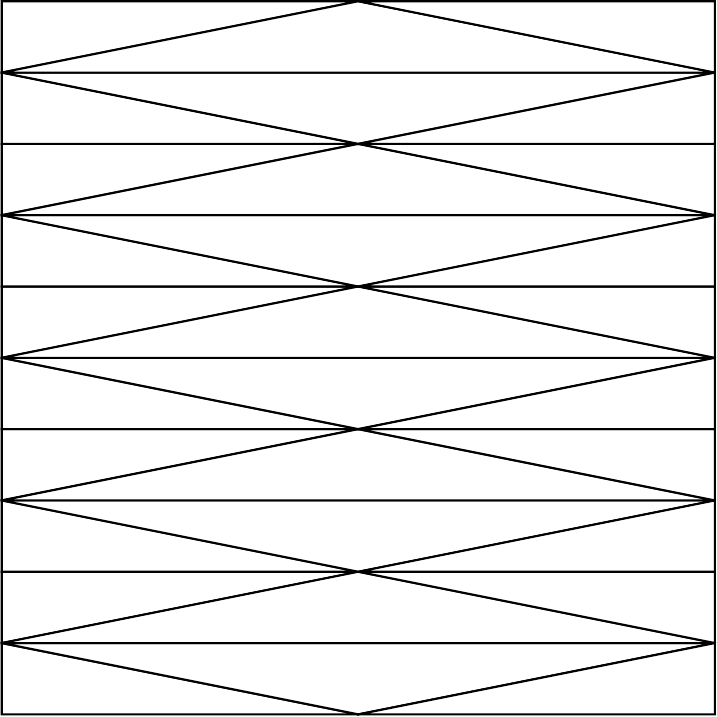}}
  \hfill
  \subfloat[$4\times 40$\label{fig:meshes_poisson_e}]{\includegraphics[width=0.29\textwidth]{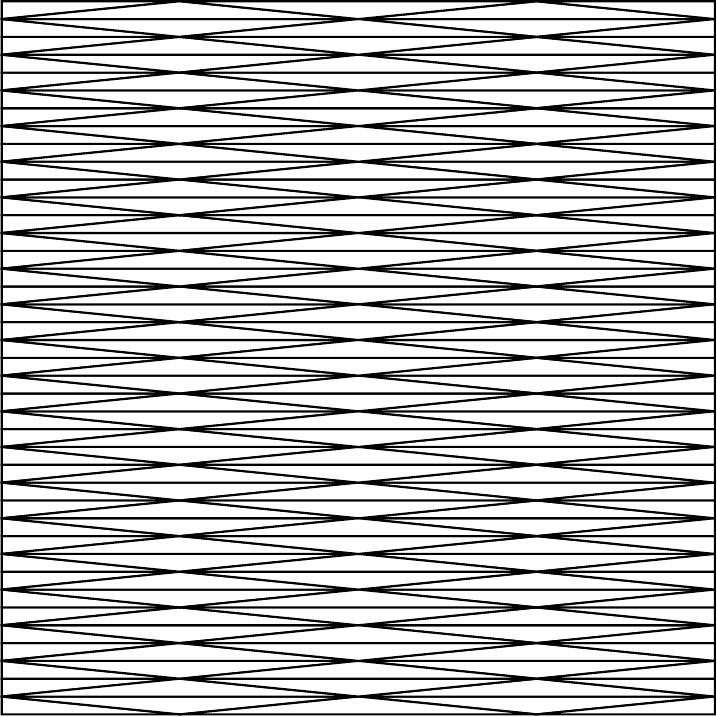}}
  \hfill
  \subfloat[Nonconvergence\label{fig:rates_poisson_fem_b}]{\includegraphics[width=0.38\textwidth]{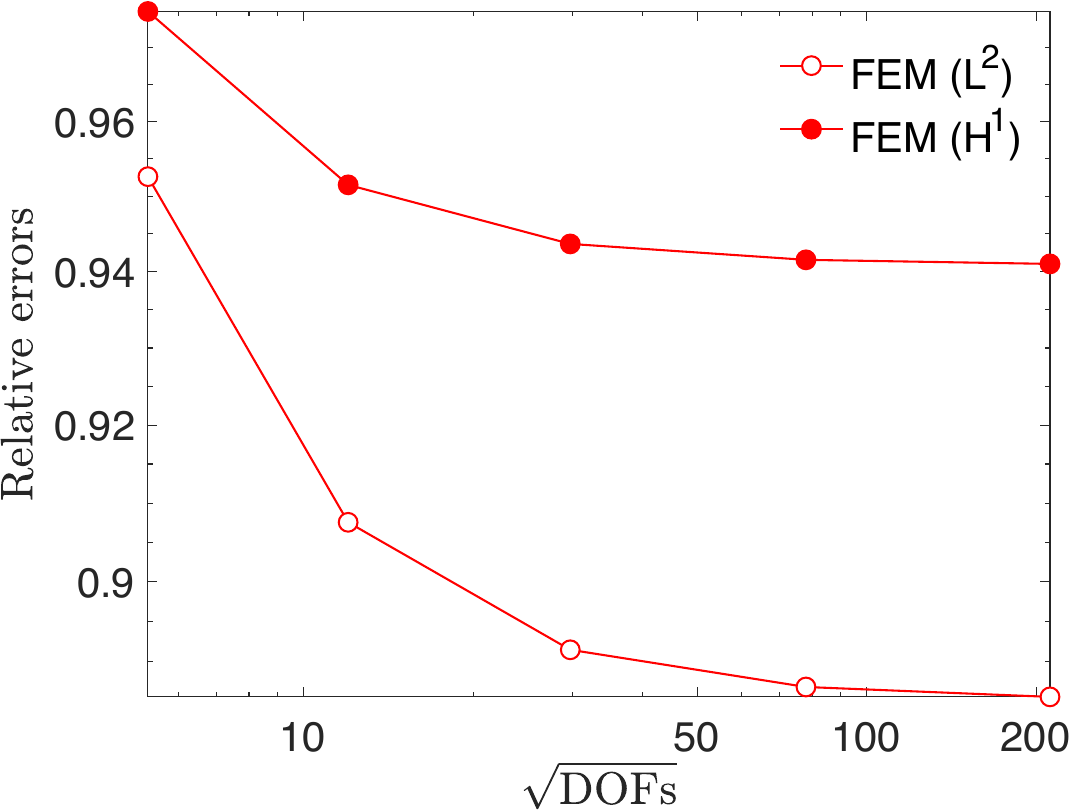}}
  \\[5pt]
  \subfloat[$2\times 10$,
  agglomerated\label{fig:meshes_poisson_g}]{\includegraphics[width=0.29\textwidth]{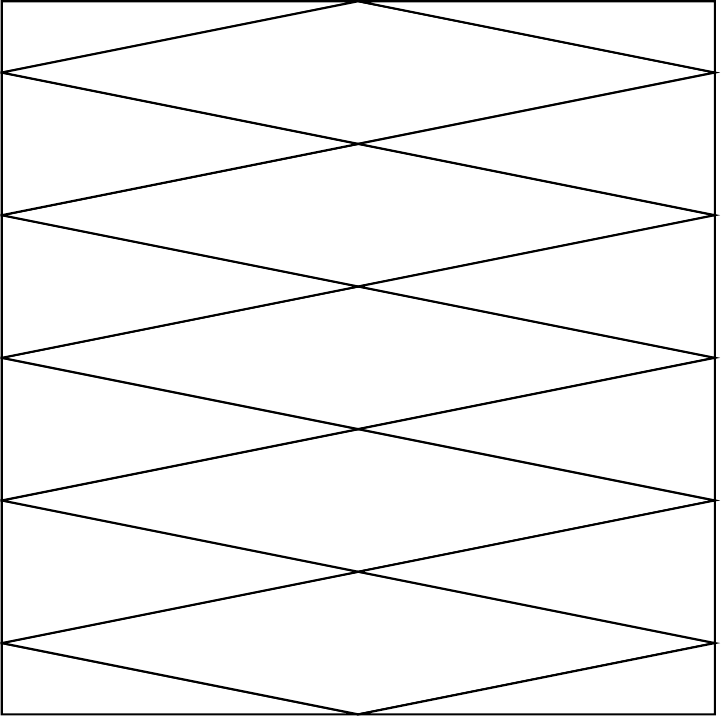}}
  \hfill
  \subfloat[$4\times 40$, agglomerated\label{fig:meshes_poisson_h}]{\includegraphics[width=0.29\textwidth]{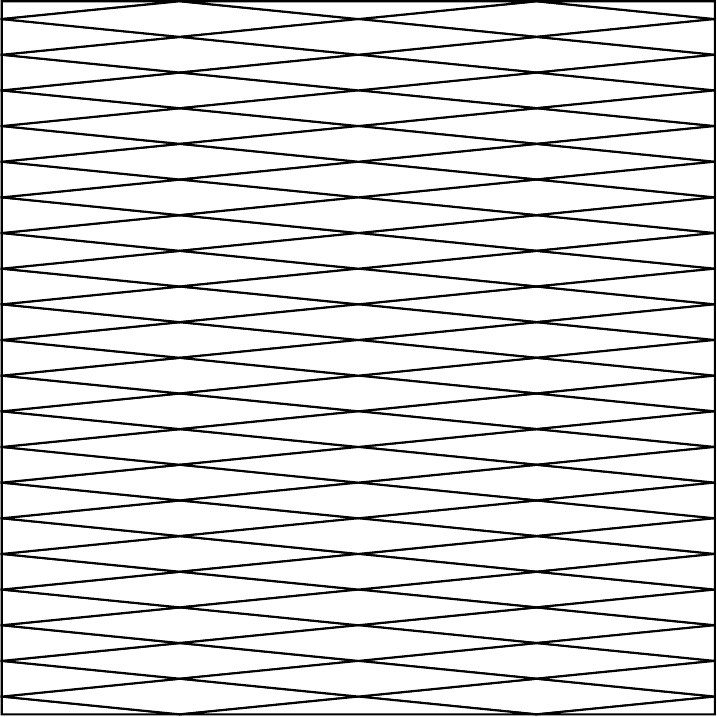}}
  \hfill
\subfloat[Convergence with CutVEM\label{fig:rates_poisson_vem}]{\includegraphics[width=0.38\textwidth]{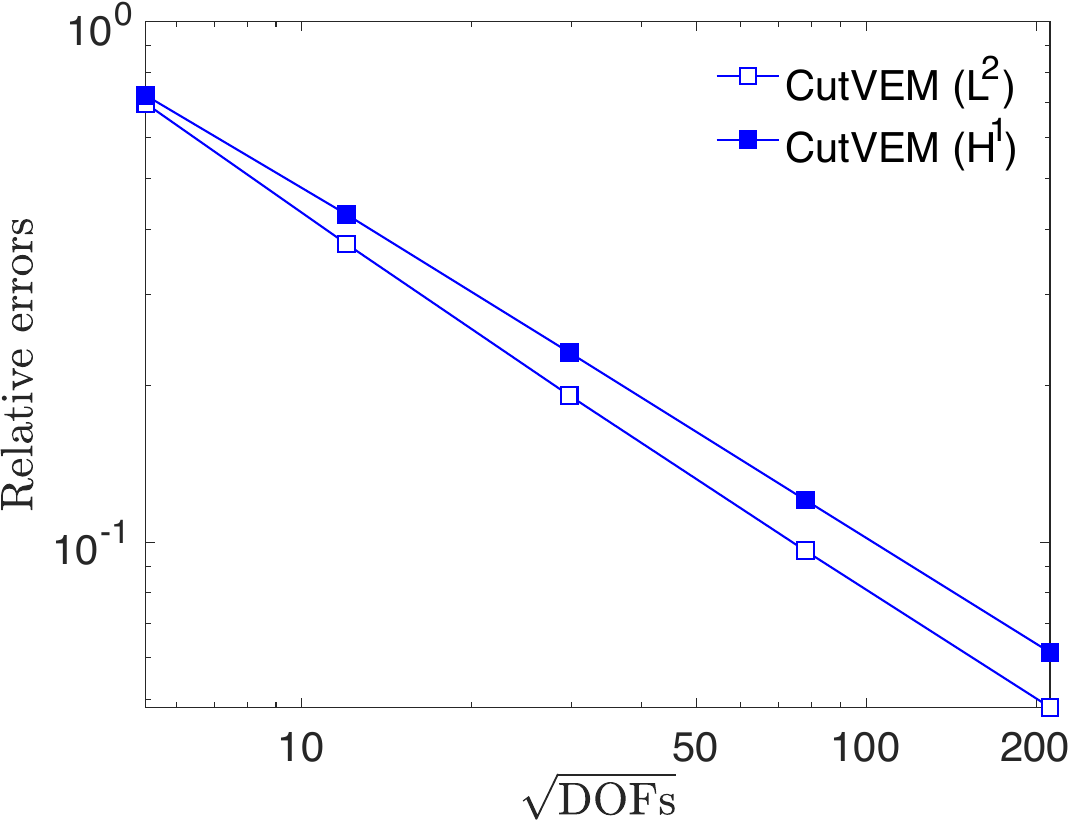}}
\caption{Convergence studies for 
Poisson problem~\eqref{eq:sinsin} using
    uniformly and anisotropically refined mesh sequences. The first
    row shows representative meshes in the former sequence for which
    plot (c) confirms that FEM solutions converge optimally. Over the
    anisotropically refined mesh sequence having flattened triangles as
    shown in the meshes in the second row, plot (f) reveals that FEM
    solutions in fact fail to converge.  With the CutVEM,
    agglomerating elements as shown in the meshes in the last row, we
    regain convergence as seen in plot (i).}
\label{fig:meshes_poisson}
\end{figure}

As a benchmark, we first consider discretizations of well-shaped
  triangles using a sequence of meshes with nodes positioned on
  $4 \times 4$, $8 \times 8$, \dots, $128 \times 128$ grids.
  Representative meshes from this sequence are shown
  in Figs.~\ref{fig:meshes_poisson_a}
  and~\ref{fig:meshes_poisson_b}. Then, we consider a sequence of
  anisotropically refined meshes, starting with a mesh having nodes on
  a $2 \times 10$ grid, and proceeding to $4\times 40, 8\times 160$,
  \dots, $32\times 2560$ grids. 
  Figures~\ref{fig:meshes_poisson_d}
  and~\ref{fig:meshes_poisson_e} 
  show the first pair of meshes in this second sequence. Notice that
  the node spacing in these meshes vary as $h$ and $h^2$ along the
  $x$- and $y$-directions, respectively. These meshes present severe
  test cases because triangle qualities deteriorate with progressive
  refinement. On performing element agglomeration on these
  anisotropically refined meshes, the first pair of agglomerated meshes in this third 
  sequence are shown
  in Figs.~\ref{fig:meshes_poisson_g}
  and~\ref{fig:meshes_poisson_h}. Virtual element method is used
  on these agglomerated meshes.
  
We perform a convergence study on the three sequences of meshes.
  The relative errors of FE solutions in the $L^2$ norm and the
  $H^1$-seminorm are plotted 
  in Figs.~\ref{fig:rates_poisson_fem_a}
  and~\ref{fig:rates_poisson_fem_b}. While
  the convergence rates are optimal for the case of uniform
  refinement, we observe nonconvergence with the mesh sequence having
  flat triangles.  On the other hand,~\fref{fig:rates_poisson_vem}
  shows
  that the errors do converge 
  when using the CutVEM on the agglomerated meshes. 
  The
  convergence rates in the $L^2$ norm and $H^1$-seminorm are 0.69 and
  0.67, respectively. Though suboptimal, which is expected considering
  the anisotropic subdivision procedure, regaining convergence clearly
  highlights the benefits of CutVEM.

\subsection{Heat conduction on clipped geometry}
\label{subsec:homogeneous_heat_conduction}
We consider the Poisson (homogeneous heat conduction with
conductivity $\kappa = 1$) equation,
\begin{subequations}\label{eq:Poisson_clipped}
\begin{align}
-\nabla^2 u &= f \ \ \textrm{in } \Omega = (0,1)^2, 
\intertext{with either Dirichlet or 
mixed (Dirichlet and Neumann) boundary conditions.
The heat source $f \in L^2(\Omega)$ is chosen 
so that the exact solution is:}
u(\vx) &= \sin(4 \pi x) \, \Bigl[ 4 \pi y - \sin(4 \pi y) \Bigr].
\end{align}
\end{subequations}
For the first case, we choose Dirichlet boundary conditions that
are consistent with the above exact solution, and for the
second case, we choose the Dirichlet boundary condition
$u = 0$ on the edges 
$x  = 0$, $x = 1$ and $y = 0$, and
apply zero Neumann conditions on the edge $y = 1$.

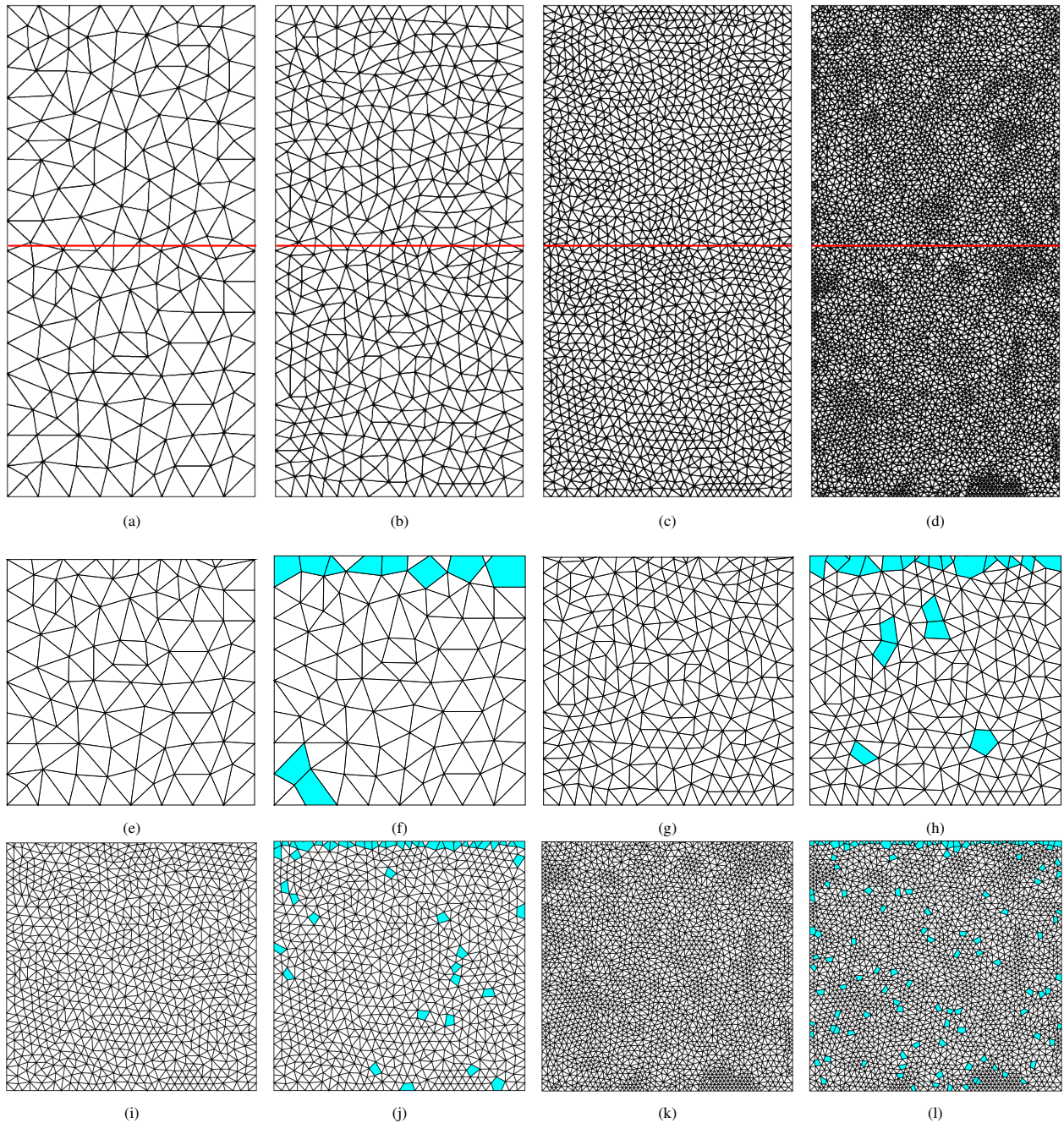
\begin{figure}[!tbh]
\centering
\begin{subfigure}{0.24\textwidth}
\includegraphics[width=\textwidth]{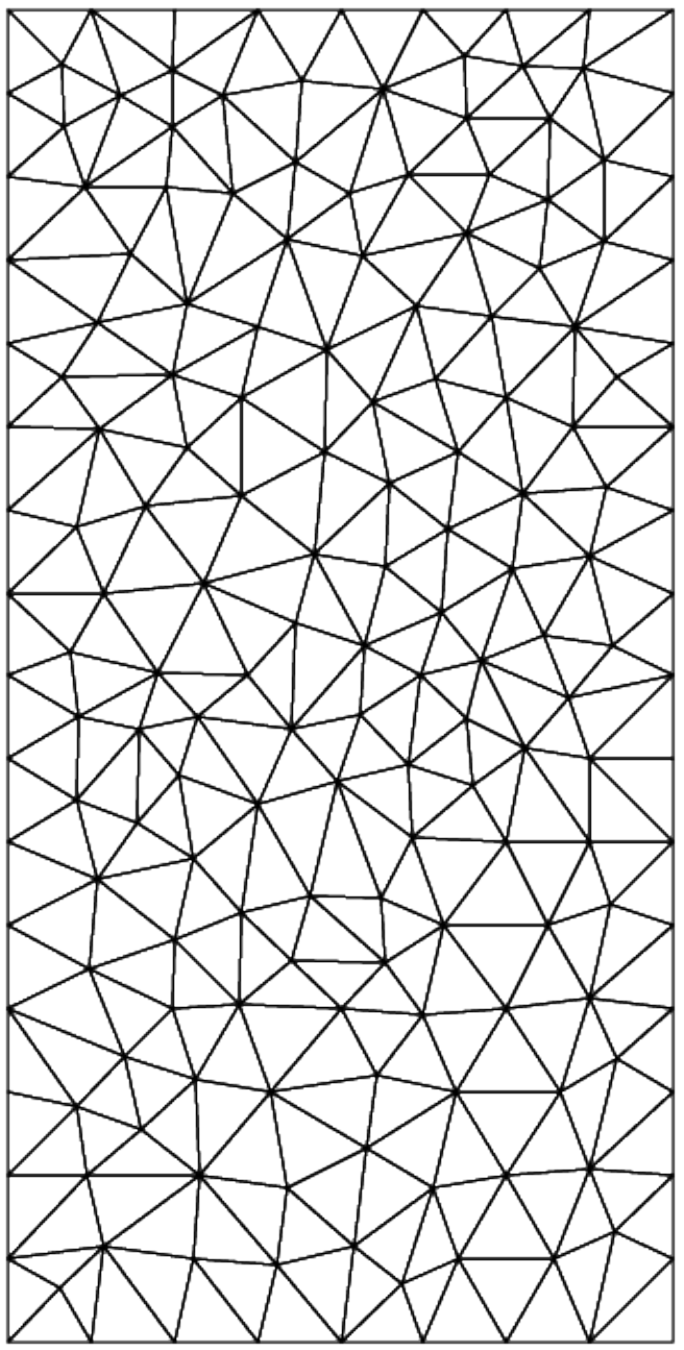}
\subcaption{}\label{fig:poisson_clipped_meshes_a}
\begin{tikzpicture}[overlay]
\draw[red,thick] (0.05,4.81) -- (3.91,4.81);  
\end{tikzpicture}
\end{subfigure} \hfill
\begin{subfigure}{0.24\textwidth}
\includegraphics[width=\textwidth]{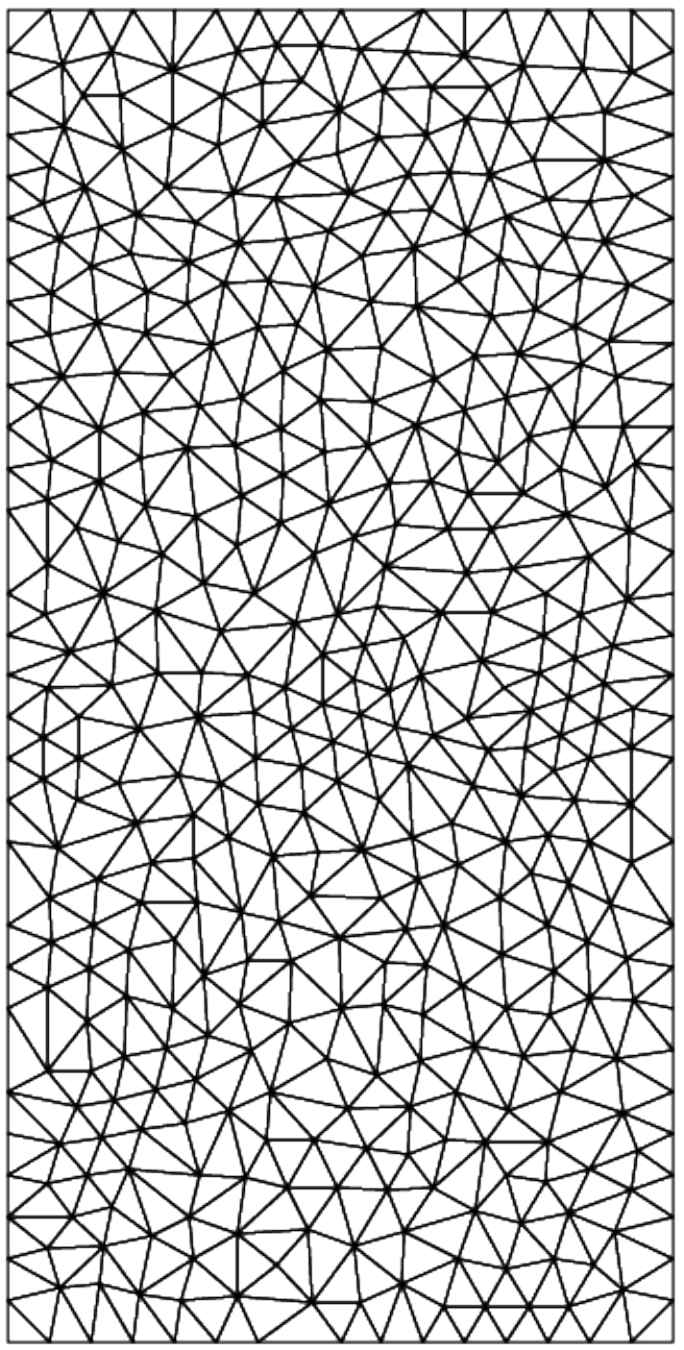}
\subcaption{}\label{fig:poisson_clipped_meshes_b}
\begin{tikzpicture}[overlay]
\draw[red,thick] (0.05,4.81) -- (3.91,4.81);  
\end{tikzpicture}
\end{subfigure} \hfill
\begin{subfigure}{0.24\textwidth}
\includegraphics[width=\textwidth]{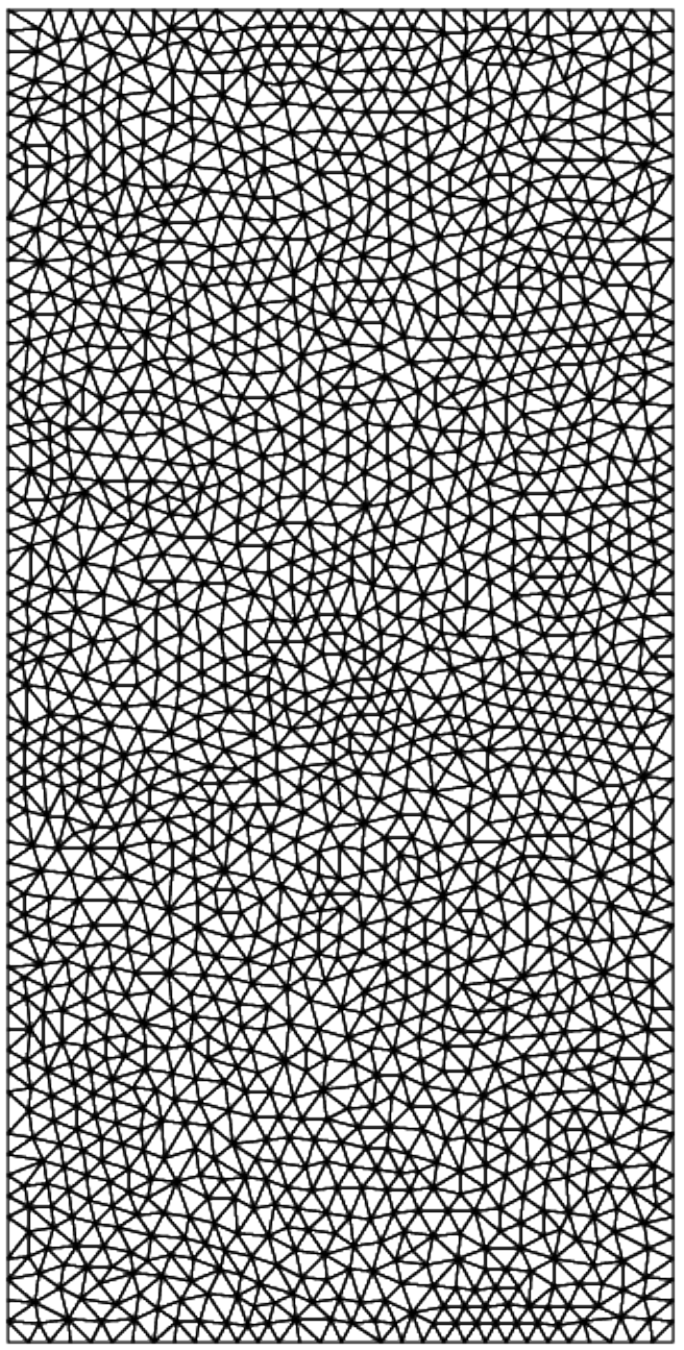}
\subcaption{}\label{fig:poisson_clipped_meshes_c}
\begin{tikzpicture}[overlay]
\draw[red,thick] (0.05,4.81) -- (3.91,4.81);  
\end{tikzpicture}
\end{subfigure} \hfill
\begin{subfigure}{0.24\textwidth}
\includegraphics[width=\textwidth]{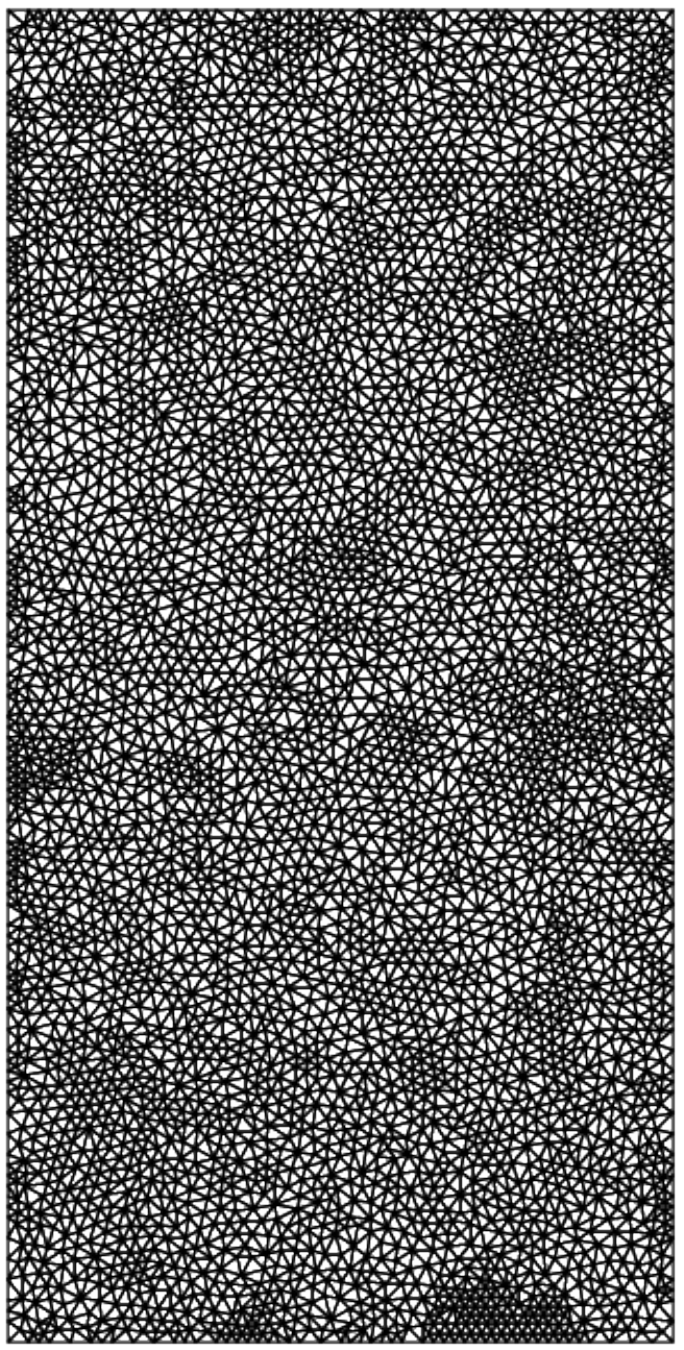}
\subcaption{}\label{fig:poisson_clipped_meshes_d}
\begin{tikzpicture}[overlay]
\draw[red,thick] (0.05,4.81) -- (3.91,4.81);  
\end{tikzpicture}
\end{subfigure}
\begin{subfigure}{0.24\textwidth}
\includegraphics[width=\textwidth]{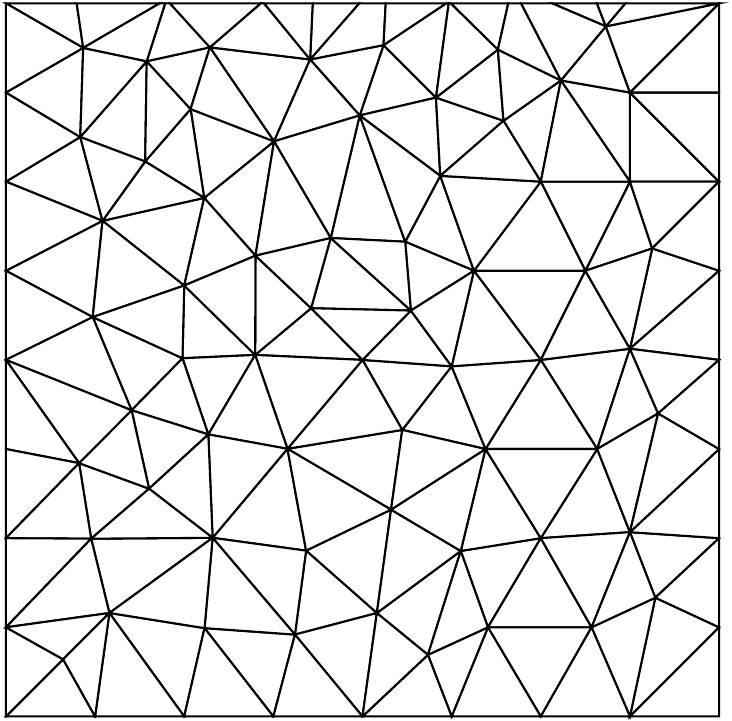}
\subcaption{}\label{fig:poisson_clipped_meshes_e}
\end{subfigure} \hfill
\begin{subfigure}{0.24\textwidth}
\includegraphics[width=\textwidth]{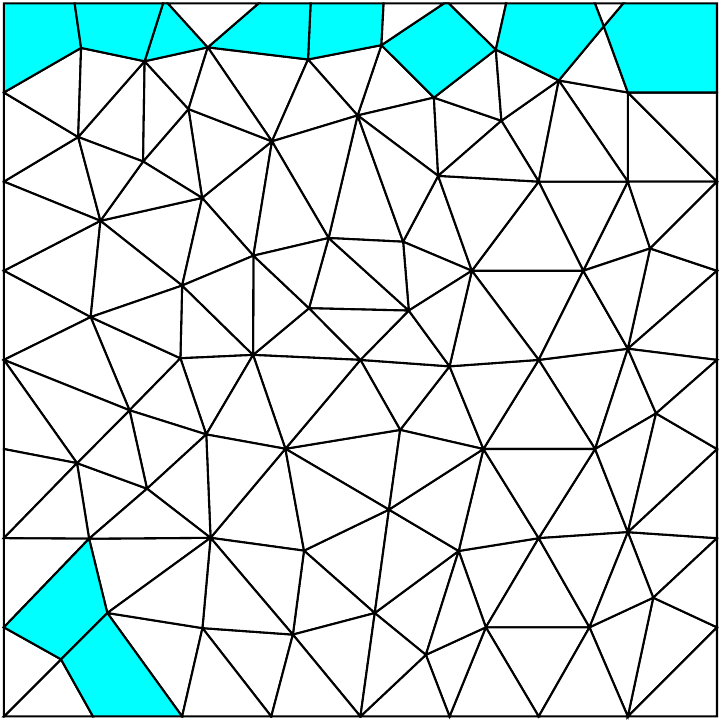}
\subcaption{}\label{fig:poisson_clipped_meshes_f}
\end{subfigure} \hfill
\begin{subfigure}{0.24\textwidth}
\includegraphics[width=\textwidth]{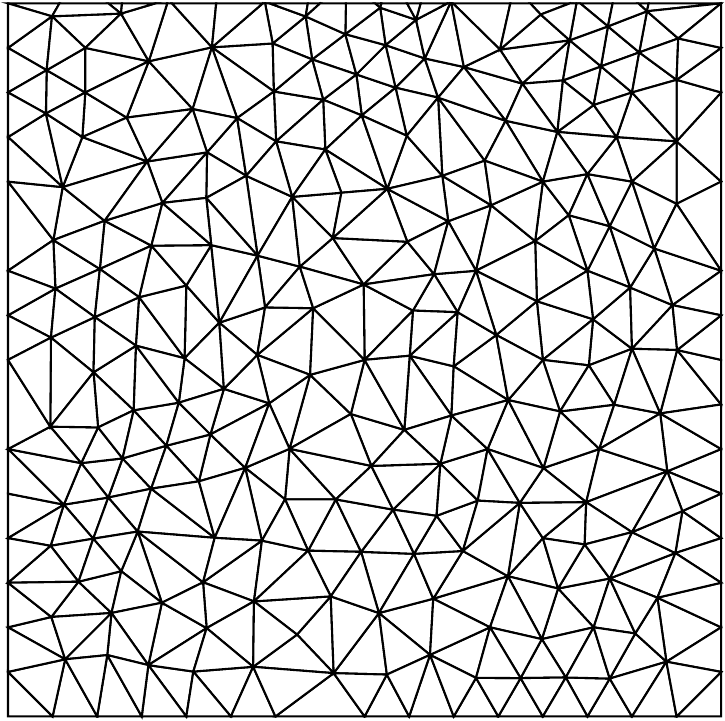}
\subcaption{}\label{fig:poisson_clipped_meshes_g}
\end{subfigure} \hfill
\begin{subfigure}{0.24\textwidth}
\includegraphics[width=\textwidth]{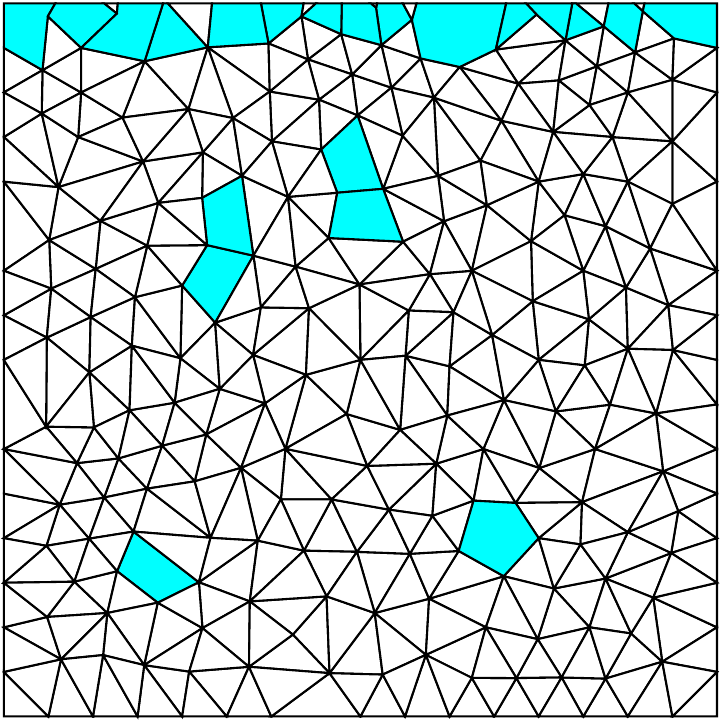}
\subcaption{}\label{fig:poisson_clipped_meshes_h}
\end{subfigure} 
\begin{subfigure}{0.24\textwidth}
\includegraphics[width=\textwidth]{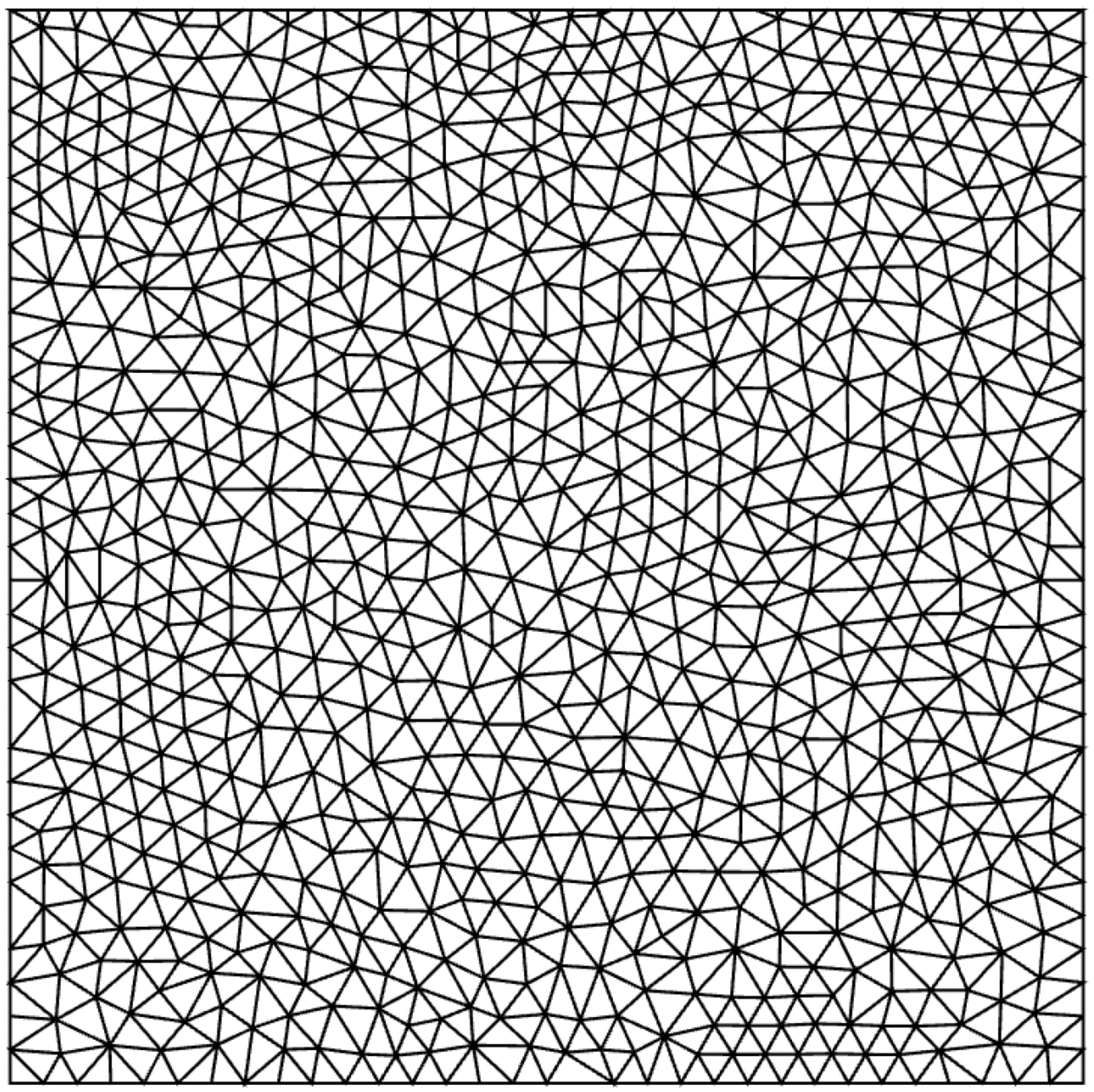}
\subcaption{}\label{fig:poisson_clipped_meshes_i}
\end{subfigure} \hfill
\begin{subfigure}{0.24\textwidth}
\includegraphics[width=\textwidth]{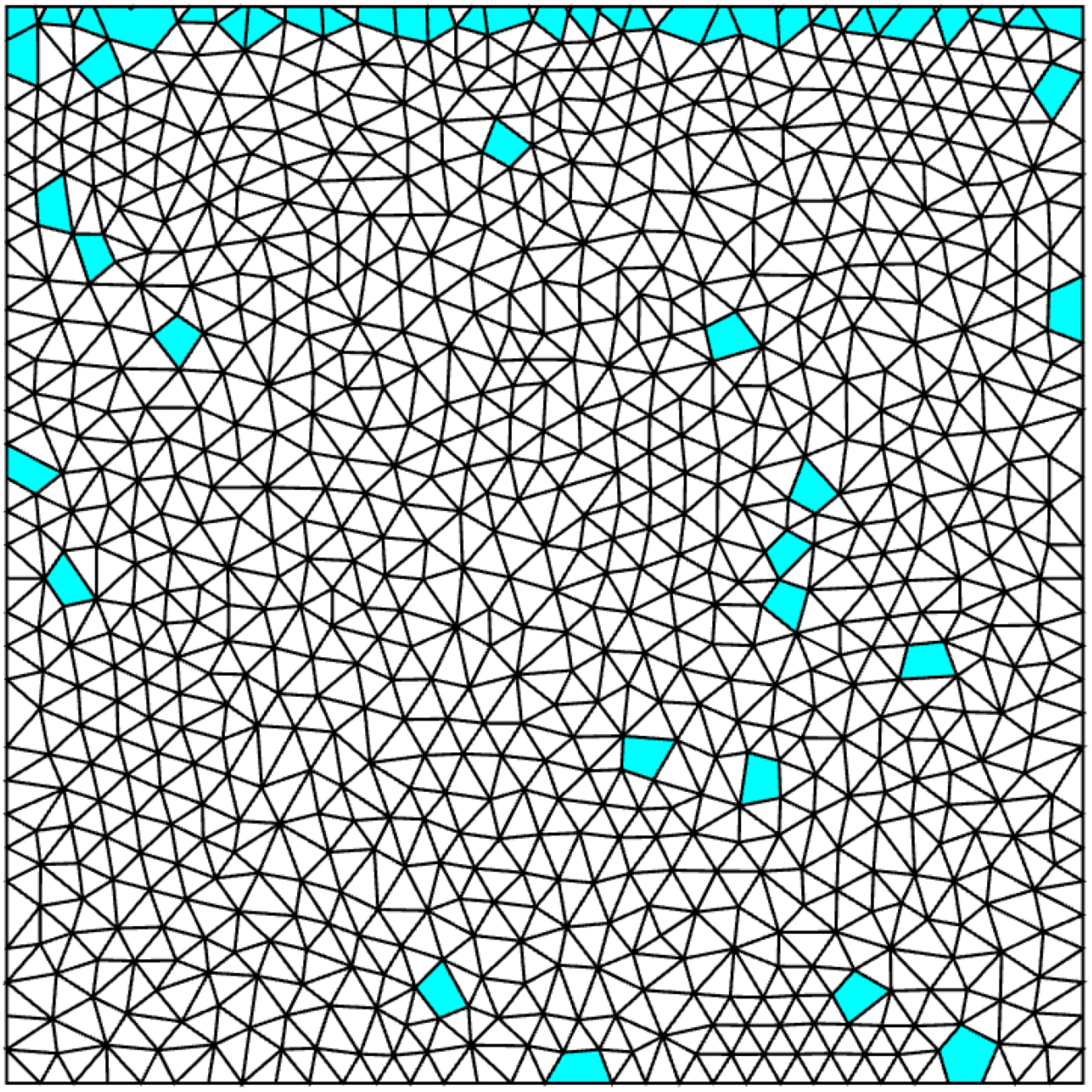}
\subcaption{}\label{fig:poisson_clipped_meshes_j}
\end{subfigure} \hfill
\begin{subfigure}{0.24\textwidth}
\includegraphics[width=\textwidth]{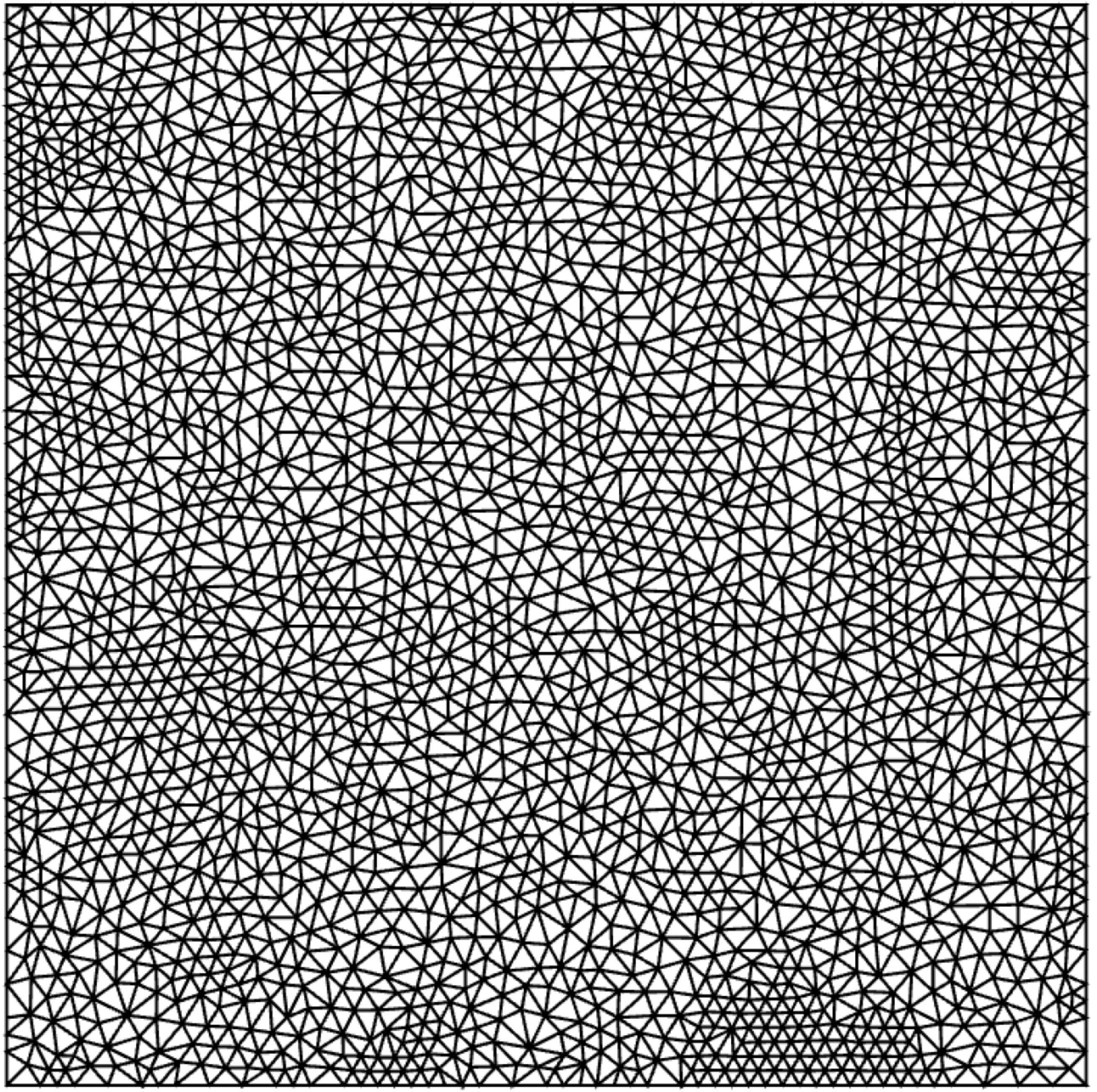}
\subcaption{}\label{fig:poisson_clipped_meshes_k}
\end{subfigure} \hfill
\begin{subfigure}{0.24\textwidth}
\includegraphics[width=\textwidth]{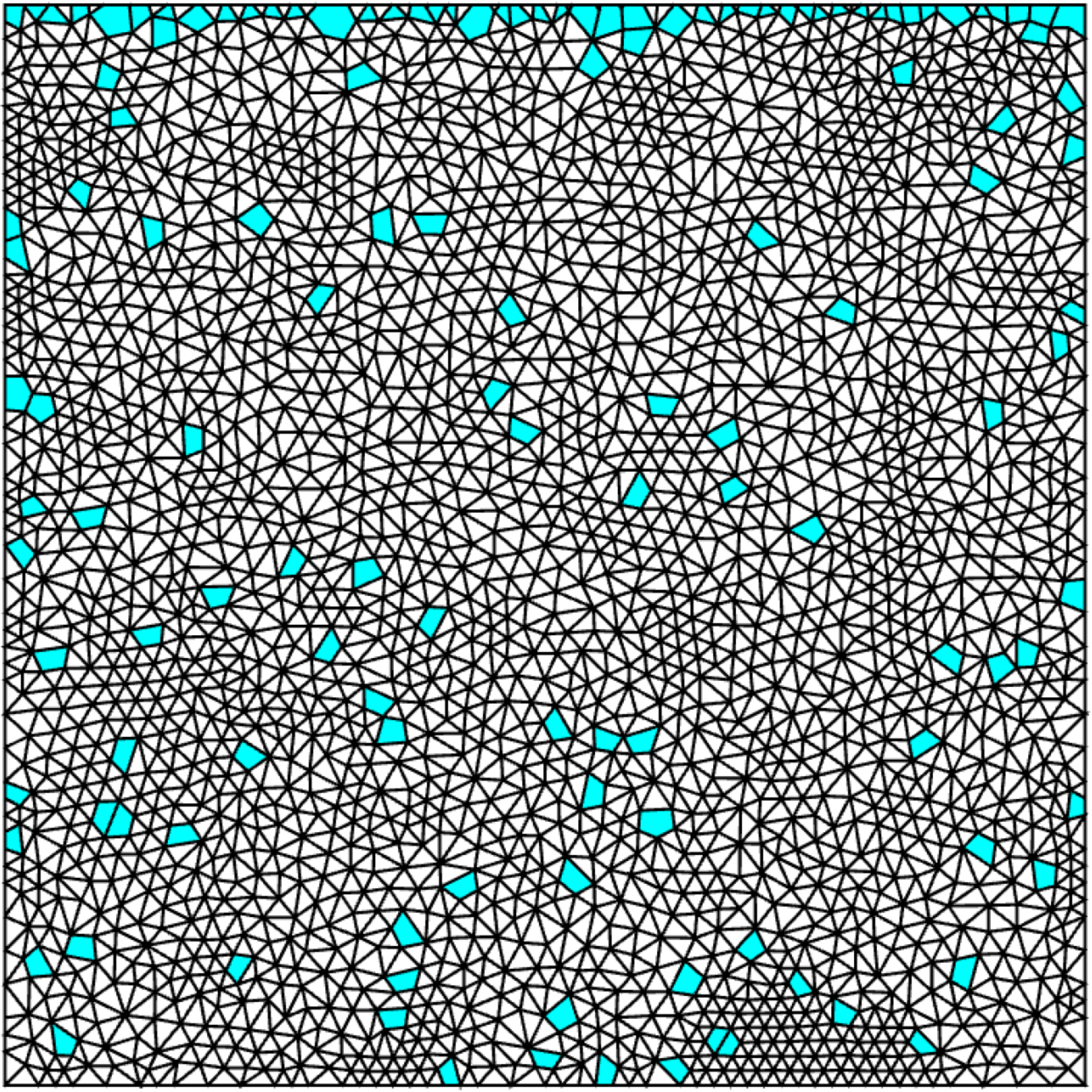}
\subcaption{}\label{fig:poisson_clipped_meshes_l}
\end{subfigure} 
\caption{Clipping finite element meshes to construct agglomerated 
meshes for CutVEM. (a), (b), (c), (d)
Input (four of the six meshes
from the sequence) finite element meshes over the domain $[0,1]\times[0,2]$. (e), (g), (i), (k)  Meshes
formed by clipping the FE mesh in (a), (b), (c) 
and (d)
by the segment from $[0,1]$ to $[1,1]$. (f), (h), (j),
(l) Agglomerated meshes formed from (e), (g), (i) and (k).
Elements with four or more edges (most of them are agglomerated elements) are  shown as filled polygons.}
\label{fig:poisson_clipped_meshes}
\end{figure}

\begin{figure}[!tbh]
\centering
\begin{subfigure}{0.24\textwidth}
\includegraphics[width=\textwidth]{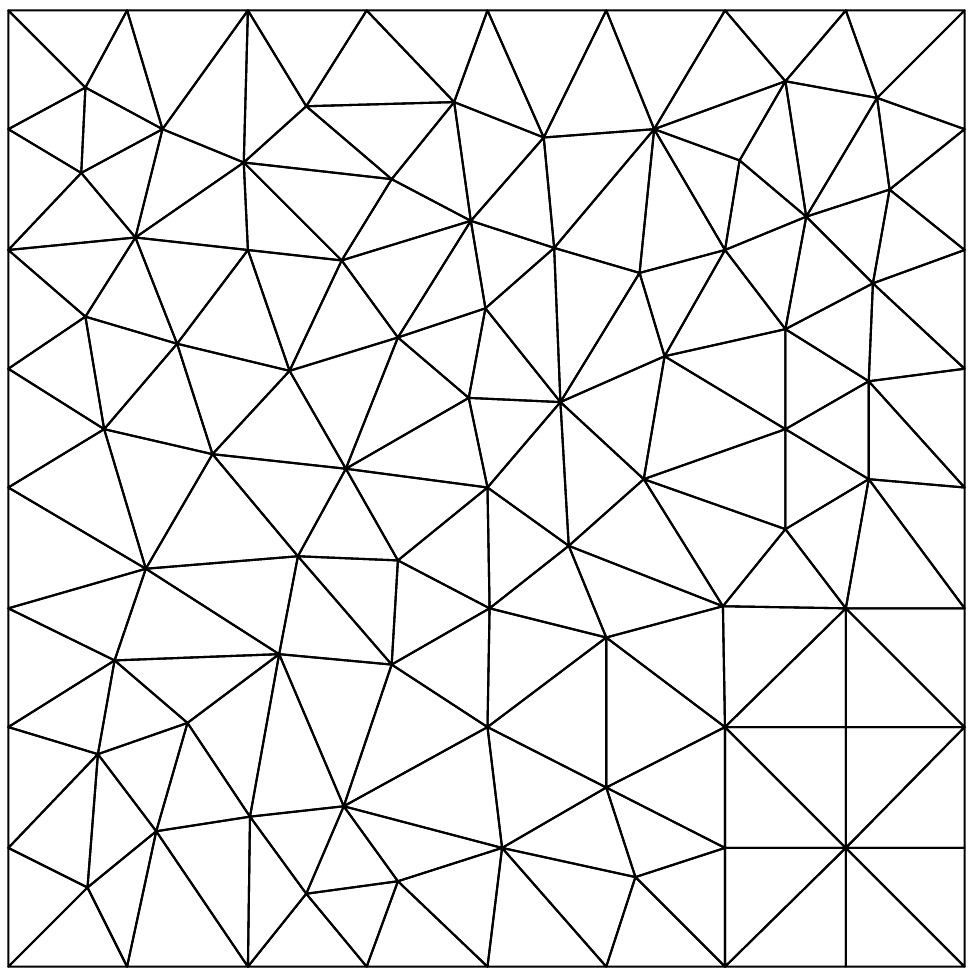}
\subcaption{}\label{fig:poisson_unitsquare_FEM_meshes_a}
\end{subfigure} \hfill
\begin{subfigure}{0.24\textwidth}
\includegraphics[width=\textwidth]{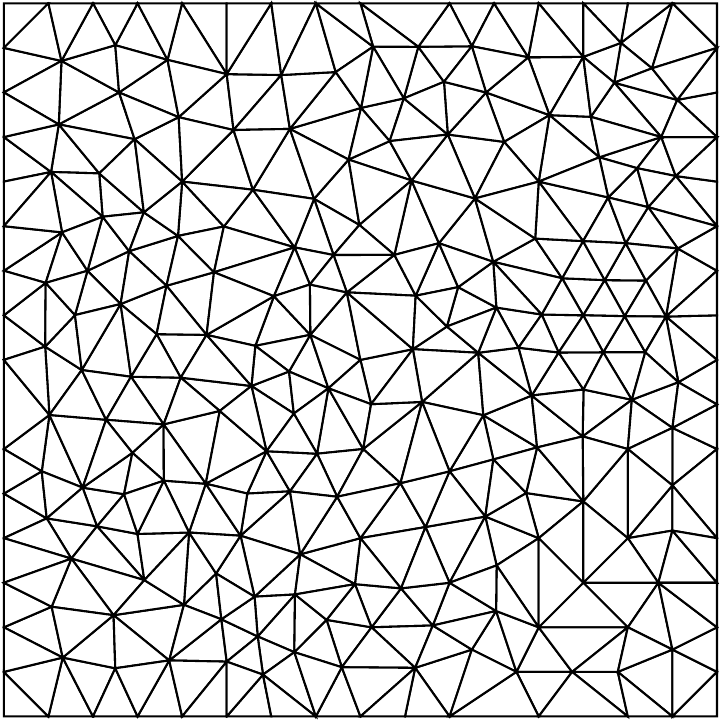}
\subcaption{}\label{fig:poisson_unitsquare_FEM_meshes_b}
\end{subfigure} \hfill
\begin{subfigure}{0.24\textwidth}
\includegraphics[width=\textwidth]{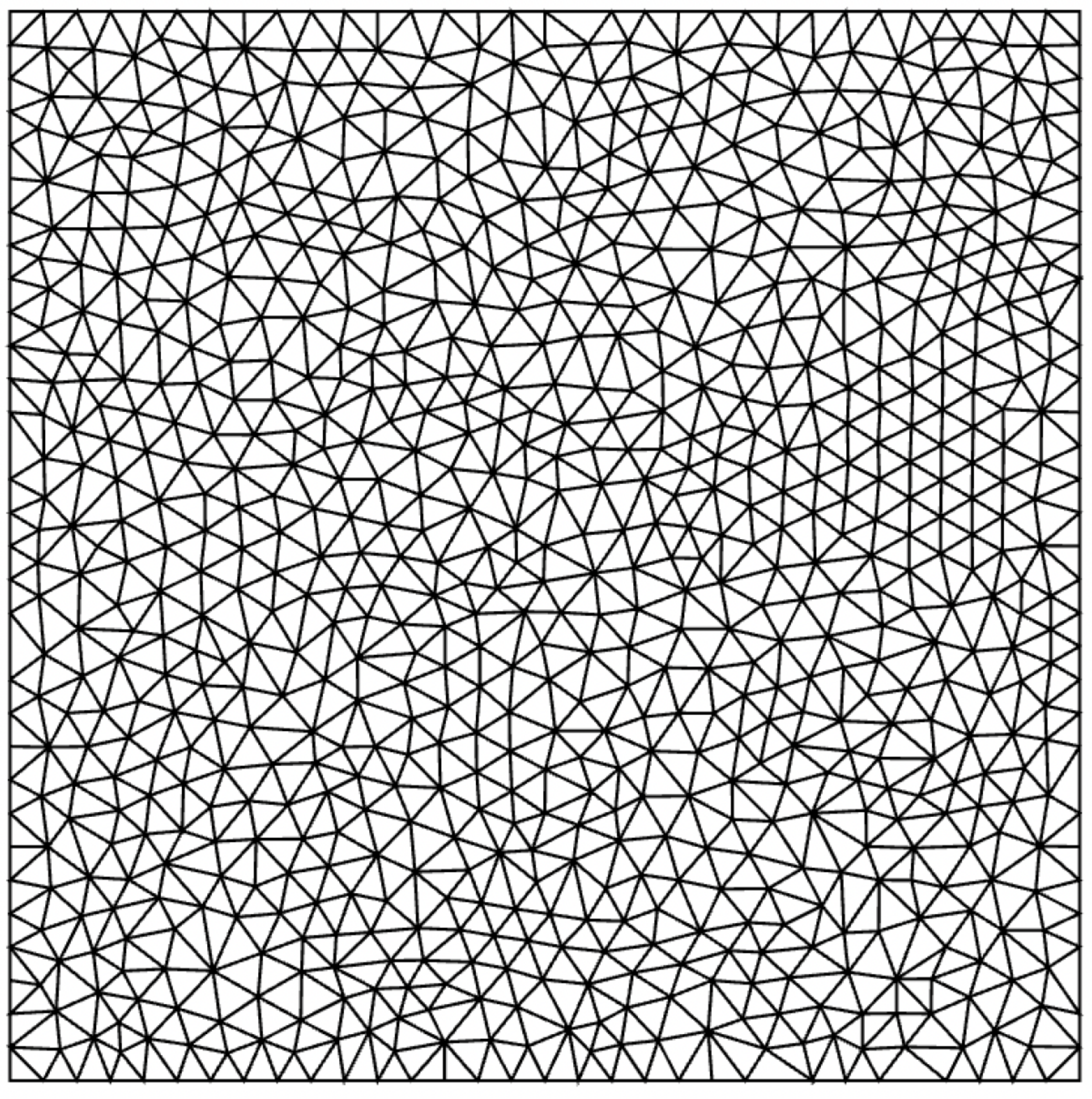}
\subcaption{}\label{fig:poisson_unitsquare_FEM_meshes_c}
\end{subfigure}
\begin{subfigure}{0.24\textwidth}
\includegraphics[width=\textwidth]{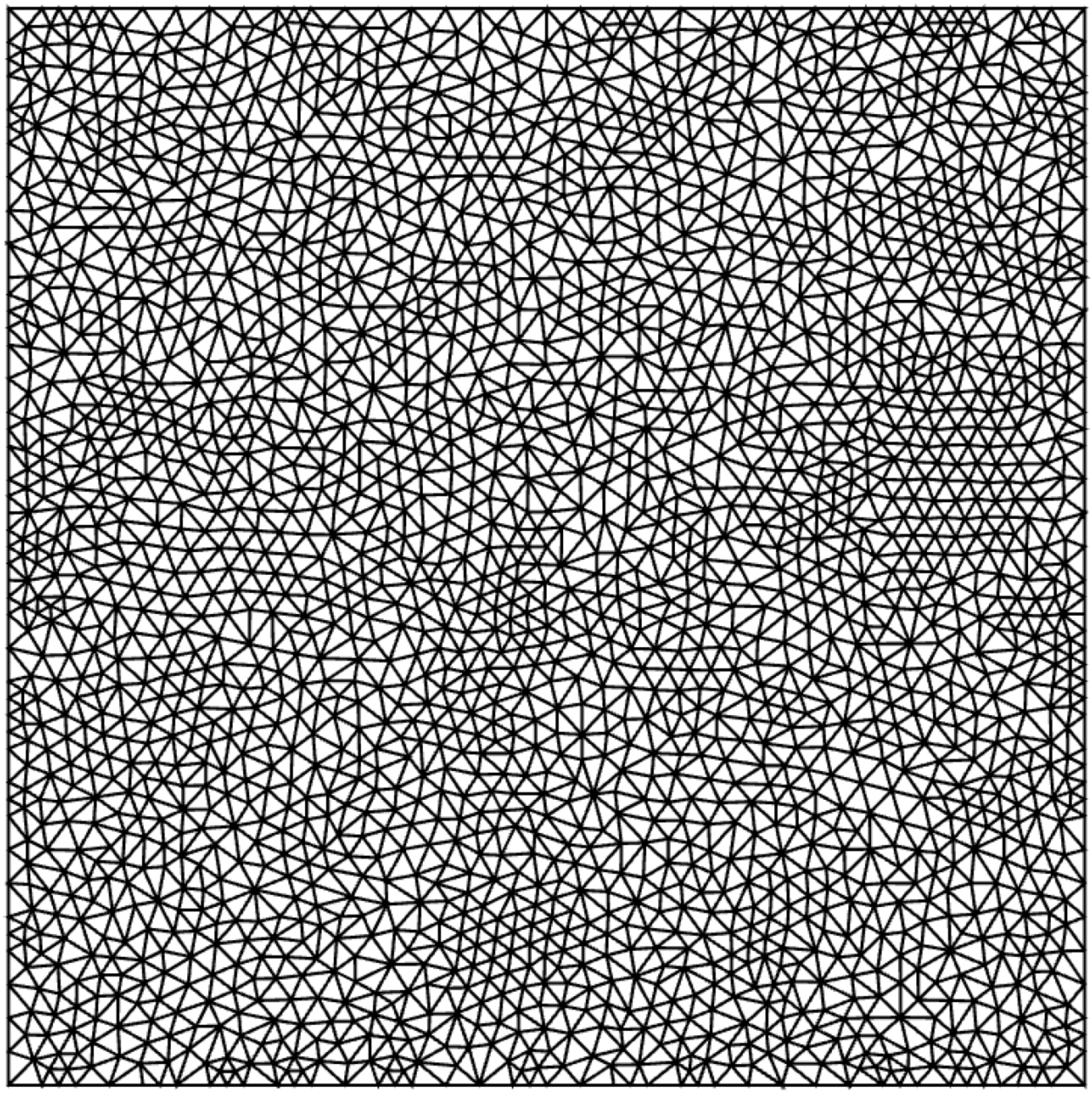}
\subcaption{}\label{fig:poisson_unitsquare_FEM_meshes_d}
\end{subfigure}
\caption{Four finite element
meshes from the sequence over $[0,1]^2$.}
\label{fig:poisson_unitsquare_FEM_meshes}
\end{figure}

We generate six Delaunay triangle meshes over
the domain
$[0,1]\times [0,2]$.  The first 
four meshes in this sequence are reproduced 
in~Figs.~\ref{fig:poisson_clipped_meshes_a}--\ref{fig:poisson_clipped_meshes_d}, where they are clipped along $y = 0$.  The region below the
clipping segment constitutes the problem domain $[0,1]^2$ over which
we solve the Poisson problem.
In Figs.~\ref{fig:poisson_clipped_meshes_e},
\ref{fig:poisson_clipped_meshes_g}, 
\ref{fig:poisson_clipped_meshes_i},
and~\ref{fig:poisson_clipped_meshes_k},
the clipped meshes are depicted.  The agglomeration algorithm discussed
in~\sref{sec:3} is
applied on these meshes, which results in the agglomerated meshes that
are presented 
in Figs.~\ref{fig:poisson_clipped_meshes_f},
\ref{fig:poisson_clipped_meshes_h},
\ref{fig:poisson_clipped_meshes_j},
and~\ref{fig:poisson_clipped_meshes_l},
respectively. 
Elements with four or more edges (most of them are agglomerated elements) are shown as filled polygons.
To assess the accuracy of CutVEM on the agglomerated meshes, we compare
its performance to the FEM. For 
a fair comparison, we use 
the~\texttt{Triangle}~\cite{Shewchuk:1996:TRI}
package to generate a sequence of Delaunay finite element meshes
over $[0,1]^2$. In Figs.~\ref{fig:poisson_unitsquare_FEM_meshes_a}--\ref{fig:poisson_unitsquare_FEM_meshes_d}, the first four finite element meshes in
the sequence are shown, and the corresponding sequence for CutVEM
appear in~\fref{fig:poisson_clipped_meshes}.
We generate FE meshes
so that they have comparable number of nodes (degrees of freedom) as those
used in CutVEM. The number of nodes range from 99 to 28,882 for the finite element analysis, and from 107 to 29,459 for the agglomerated meshes used in CutVEM.

\begin{figure}[!thb]
\centering
\begin{subfigure}{0.48\textwidth}
\includegraphics[width=\textwidth]{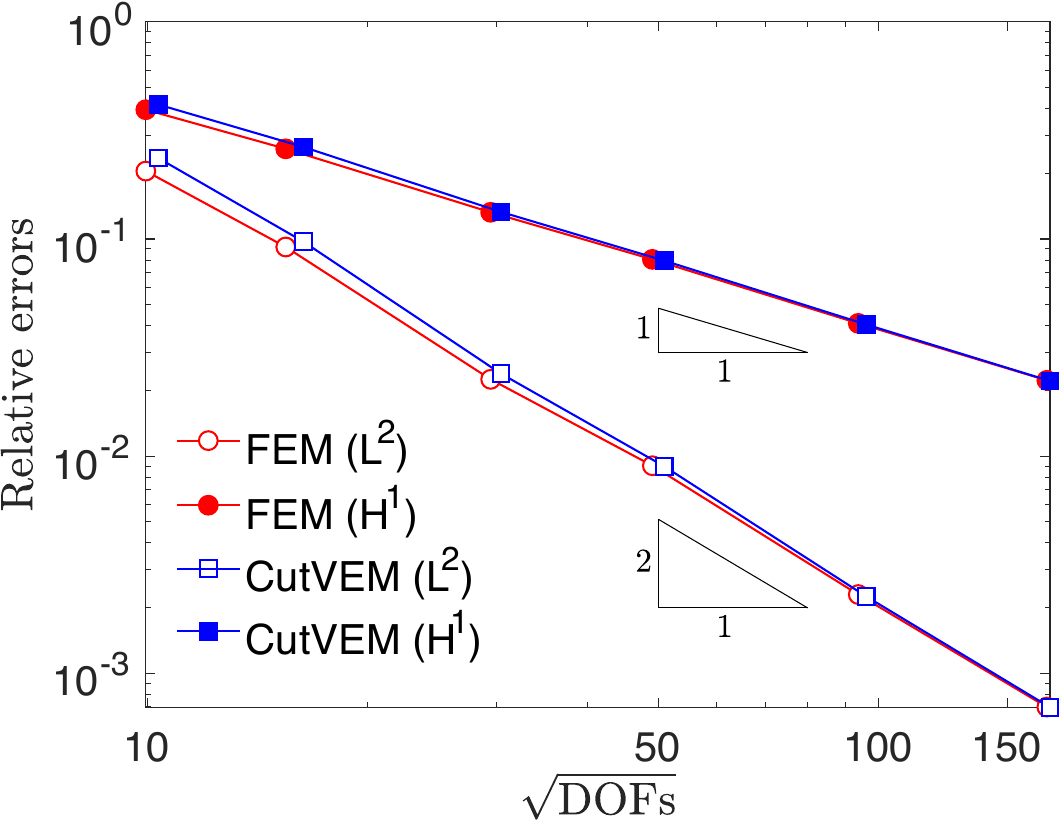}
\subcaption{}\label{fig:poisson_clipped_conv_a}
\end{subfigure} \hfill
\begin{subfigure}{0.48\textwidth}
\includegraphics[width=\textwidth]{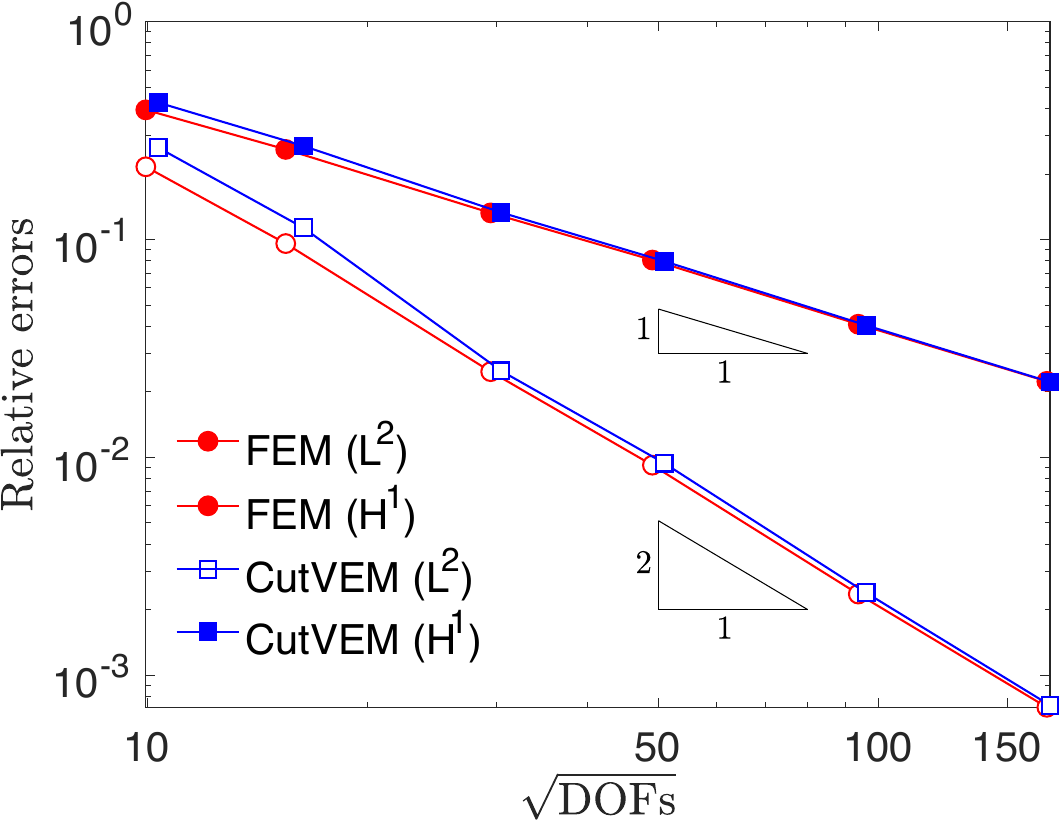}
\subcaption{}\label{fig:poisson_clipped_conv_b}
\end{subfigure}
\caption{Convergence of CutVEM on clipped meshes and  
         Delaunay FEM for the Poisson problem
         in~\protect\eqref{eq:Poisson_clipped}. 
         (a) Dirichlet boundary
         conditions; and
         (b)  mixed (Dirichlet and Neumann) boundary conditions.}\label{fig:poisson_clipped_conv}
\end{figure}

We assess the accuracy and convergence with FEM and CutVEM.  As
indicated earlier, two different boundary conditions are considered:
(i) Dirichlet and (ii) mixed (Dirichlet and Neumann).  Numerical
results are presented in~\fref{fig:poisson_clipped_conv}, with the
plots of the relative errors versus $\sqrt{\textrm{DOFs}}$ for cases
(i) and (ii) shown in Figs.~\ref{fig:poisson_clipped_conv_a}
and~\ref{fig:poisson_clipped_conv_b}, respectively. The accuracy of
CutVEM and FEM are proximal and both methods display optimal
convergence in $u$, with rates of convergence of $2$ and $1$ in the
$L^2$ norm and $H^1$ seminorm, respectively.

\subsection{Heat conduction on immersed geometries}
\label{subsec:inhomogeneneous_heat_conduction}
We consider heat conduction problems in a material with a void (hole) and 
in a two-phase composite disc with inhomogeneous thermal conductivity.

\subsubsection{Material with a hole}\label{subsubsec:hole}
Consider the annular domain $\Omega = \{r : a < r < b\}$ in polar coordinates, 
where $r = a$ and
$r = b$ are the inner and outer radii, respectively. The void
occupies the region $0 < r < a$. For thermal conductivity
$\kappa = 1$ and a unit heat source, the steady-state heat conduction equation is:
\begin{equation}\label{eq:Poisson_hole}
-\nabla^2 u = 1 \ \ \textrm{in } \Omega.
\end{equation}
On imposing zero flux boundary condition on $r = a$ and Dirichlet boundary condition 
$u(b) = 1$ on $r = b$, the exact solution for this problem is:
\begin{equation}\label{eq:exact_Poisson_hole}
u(r) =  1 + \dfrac{a^2}{2} \ln \left( \dfrac{r}{b} \right)
+ \dfrac{b^2 - r^2}{4}.
\end{equation}

\begin{figure}[!thb]
\centering
\begin{subfigure}{0.20\textwidth}
\includegraphics[width=\textwidth]{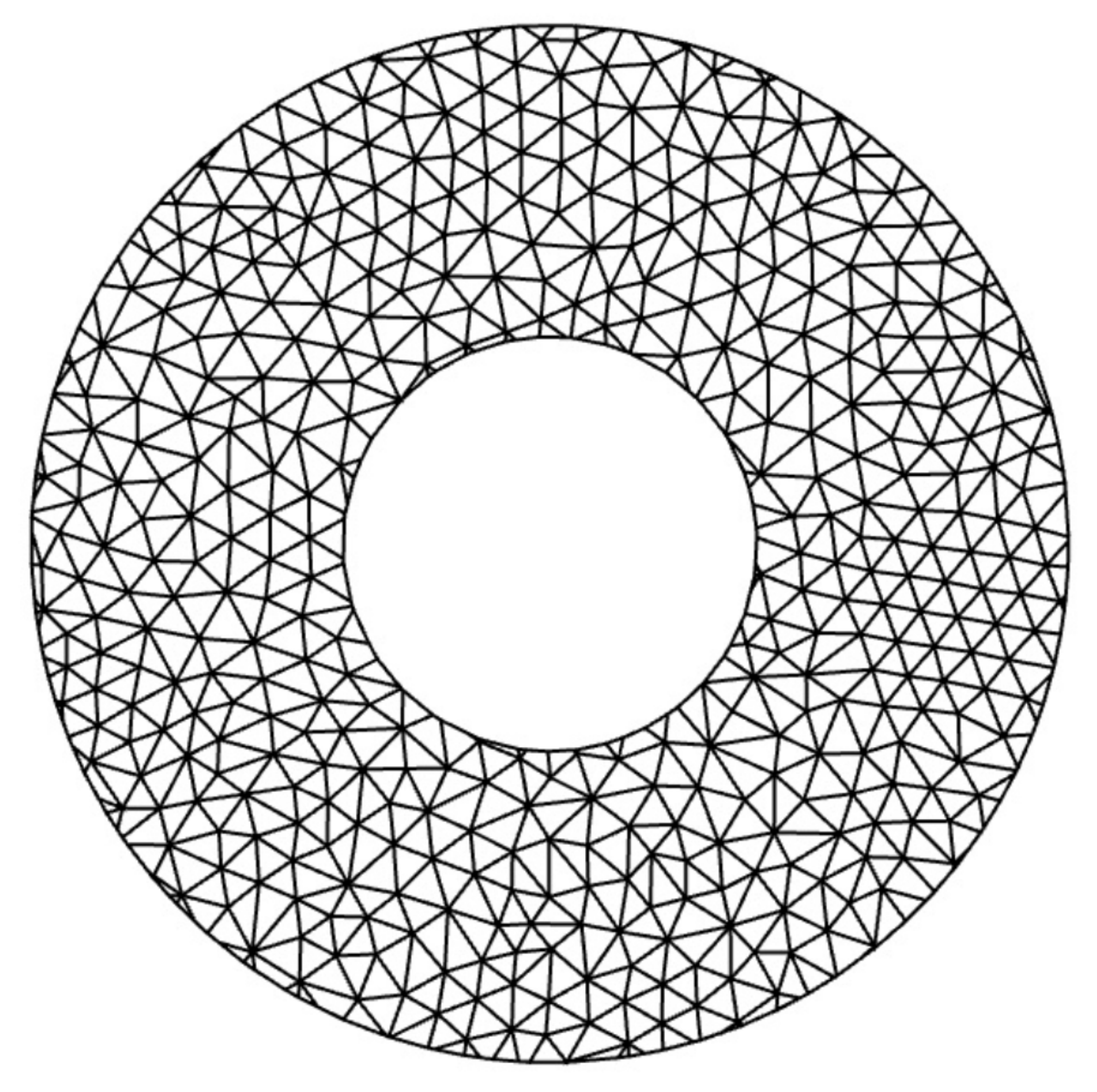}
\subcaption{}\label{fig:meshes_hole_a}
\end{subfigure} \hfill
\begin{subfigure}{0.20\textwidth}
\includegraphics[width=\textwidth]{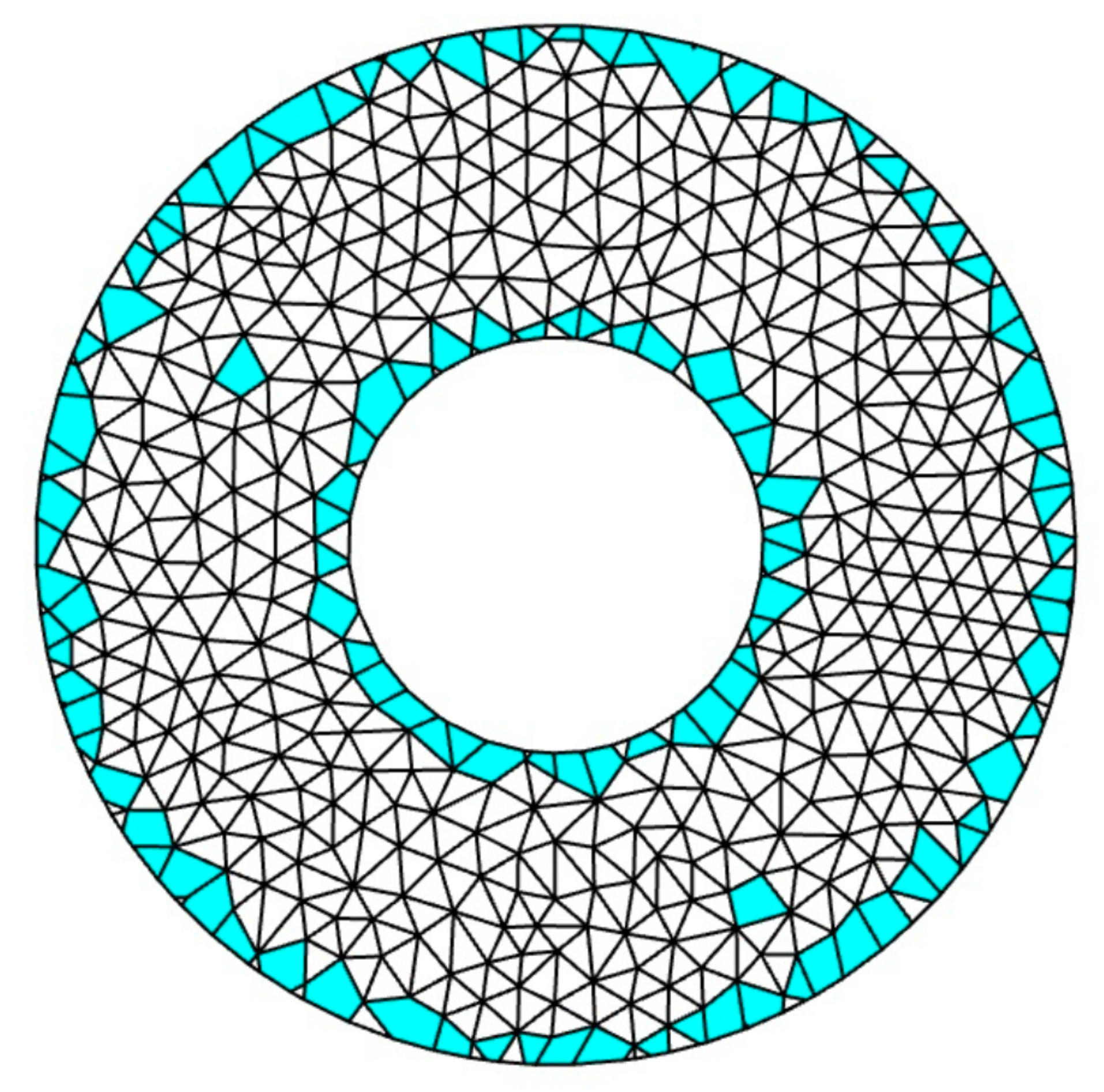}
\subcaption{}\label{fig:meshes_hole_b}
\end{subfigure} \hfill
\begin{subfigure}{0.26\textwidth}
\includegraphics[width=\textwidth]{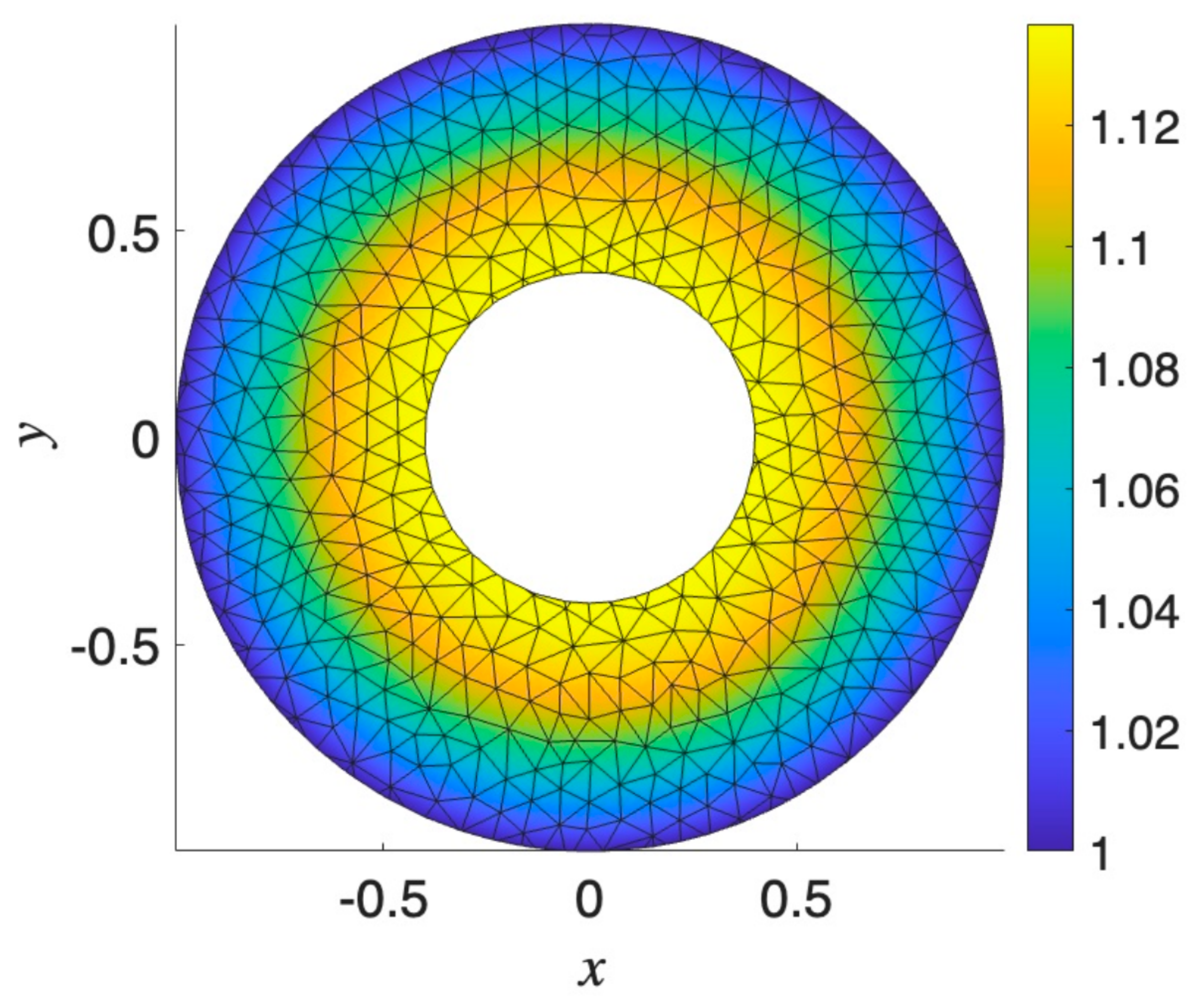}
\subcaption{}\label{fig:meshes_hole_c}
\end{subfigure} \hfill
\begin{subfigure}{0.26\textwidth}
\includegraphics[width=\textwidth]{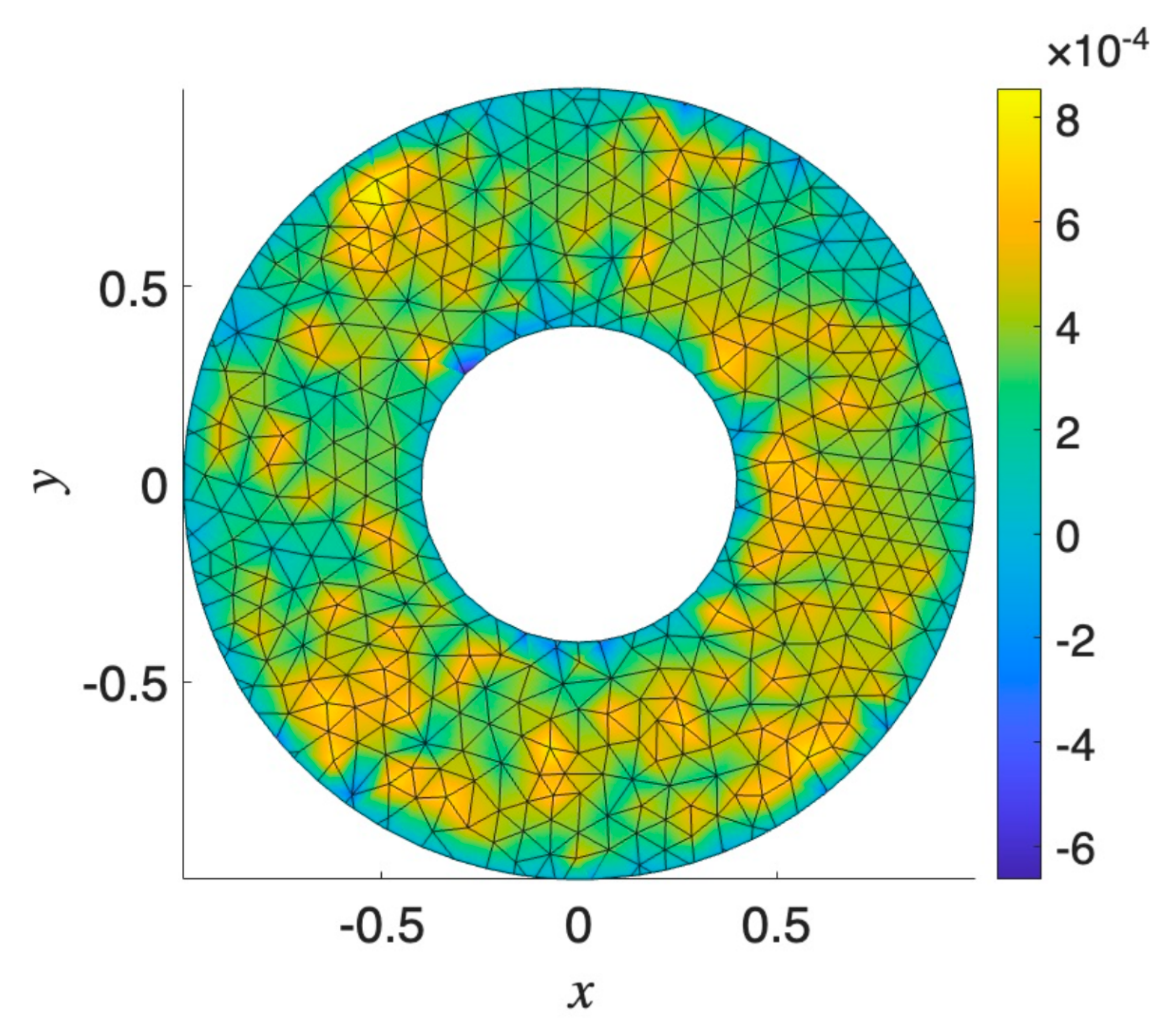}
\subcaption{}\label{fig:meshes_hole_d}
\end{subfigure}
\begin{subfigure}{0.20\textwidth}
\includegraphics[width=\textwidth]{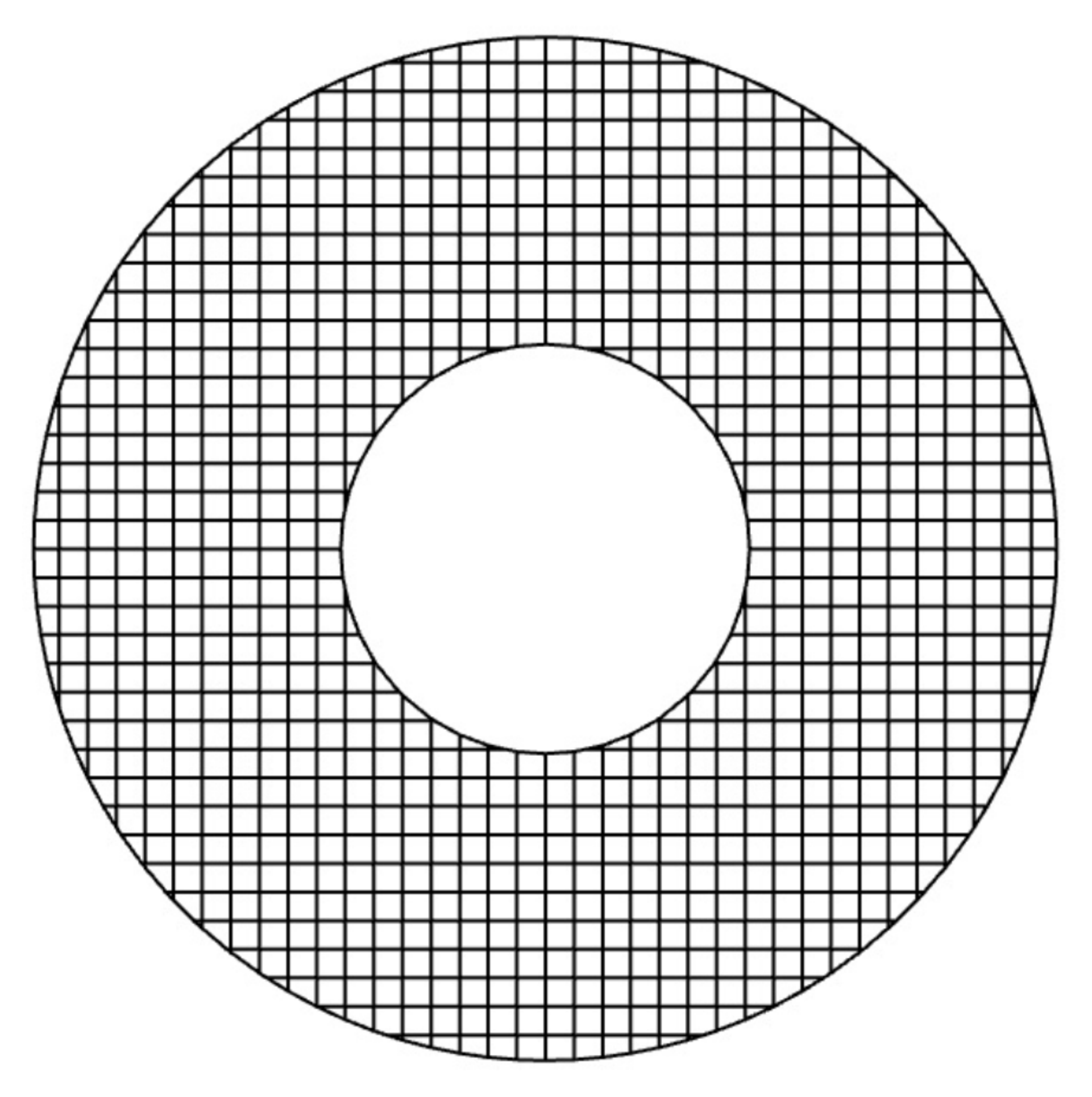}
\subcaption{}\label{fig:meshes_hole_e}
\end{subfigure} \hfill
\begin{subfigure}{0.202\textwidth}
\includegraphics[width=\textwidth]{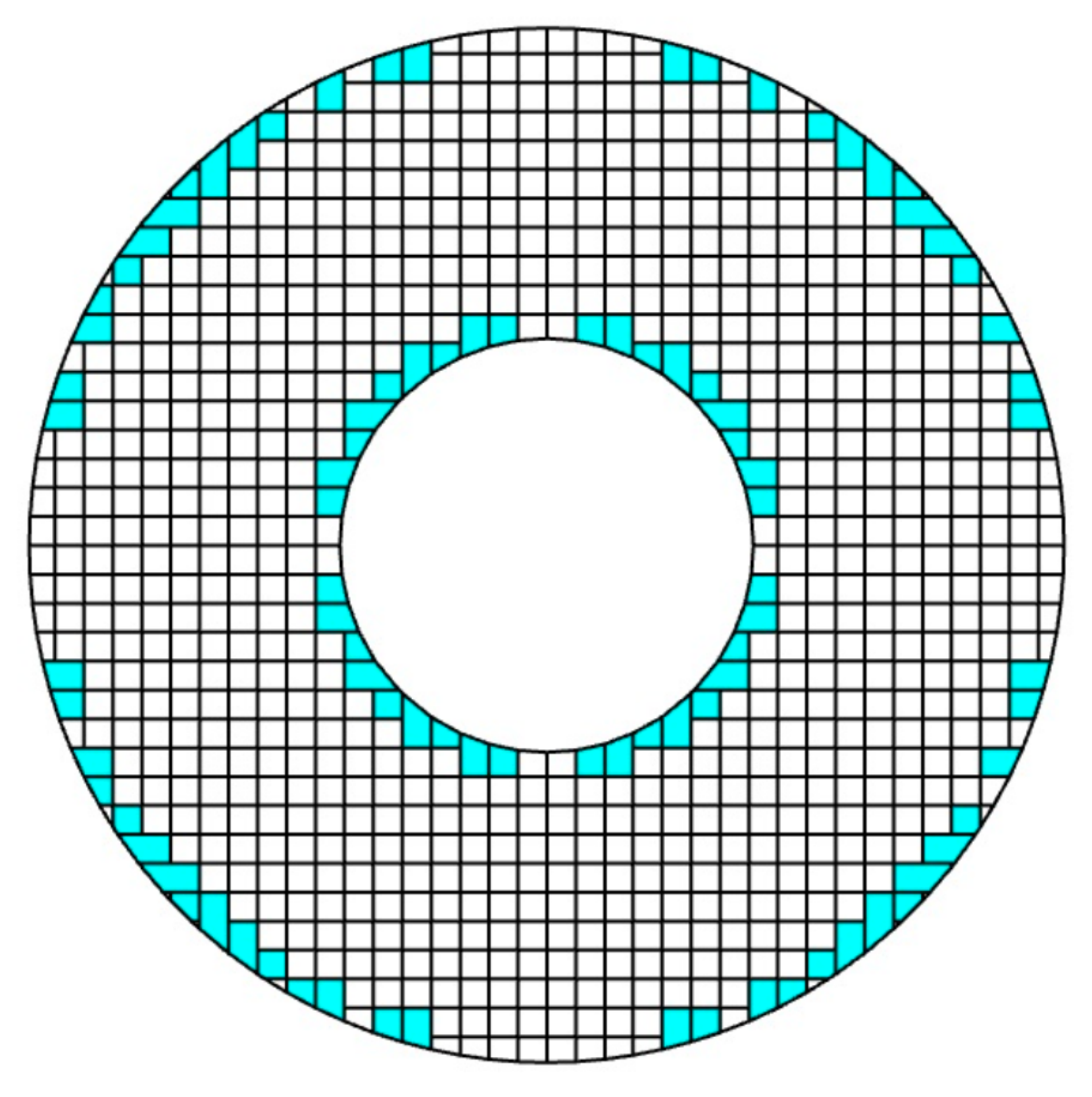}
\subcaption{}\label{fig:meshes_hole_f}
\end{subfigure} \hfill
\begin{subfigure}{0.26\textwidth}
\includegraphics[width=\textwidth]{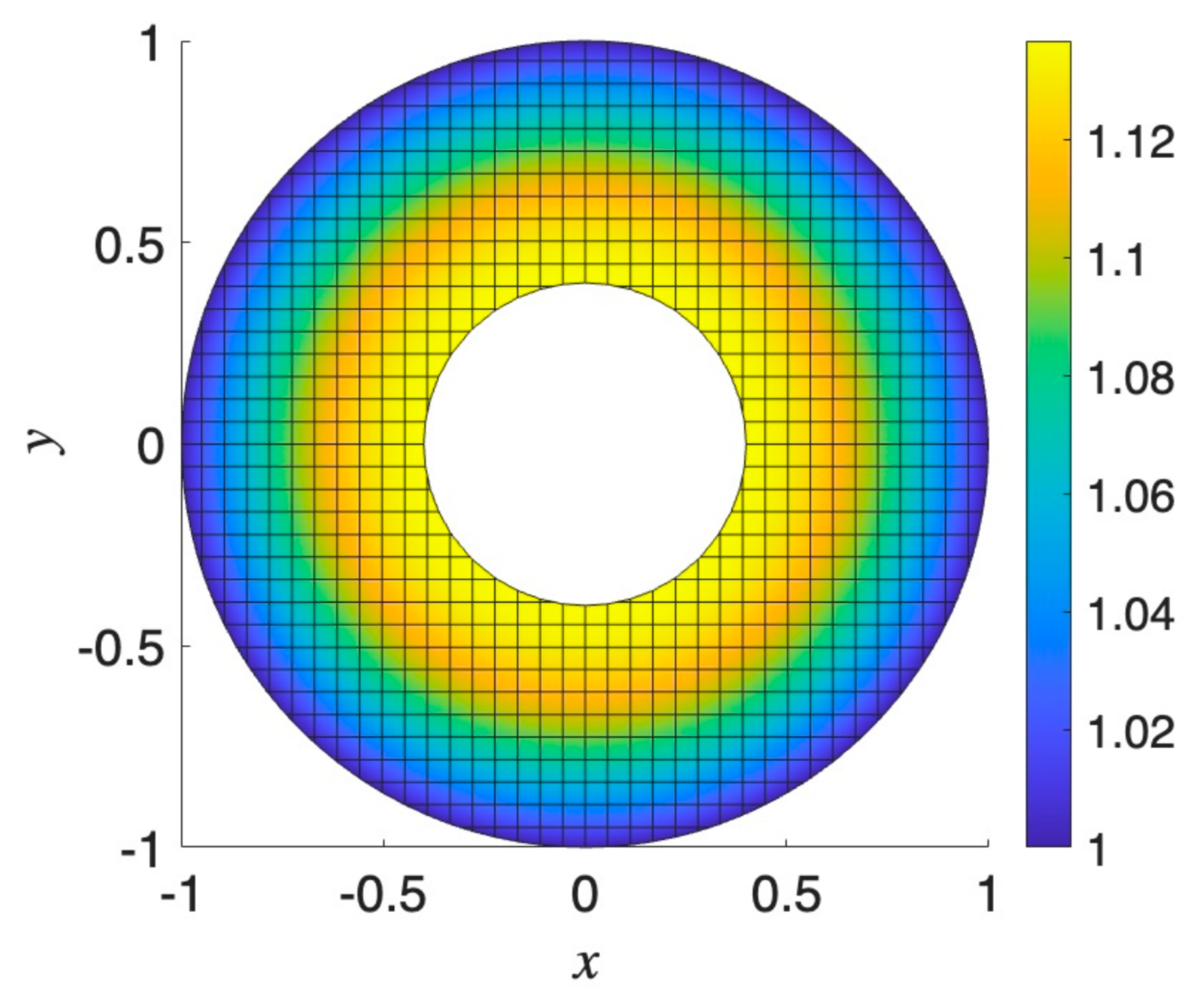}
\subcaption{}\label{fig:meshes_hole_g}
\end{subfigure} \hfill
\begin{subfigure}{0.265\textwidth}
\includegraphics[width=\textwidth]{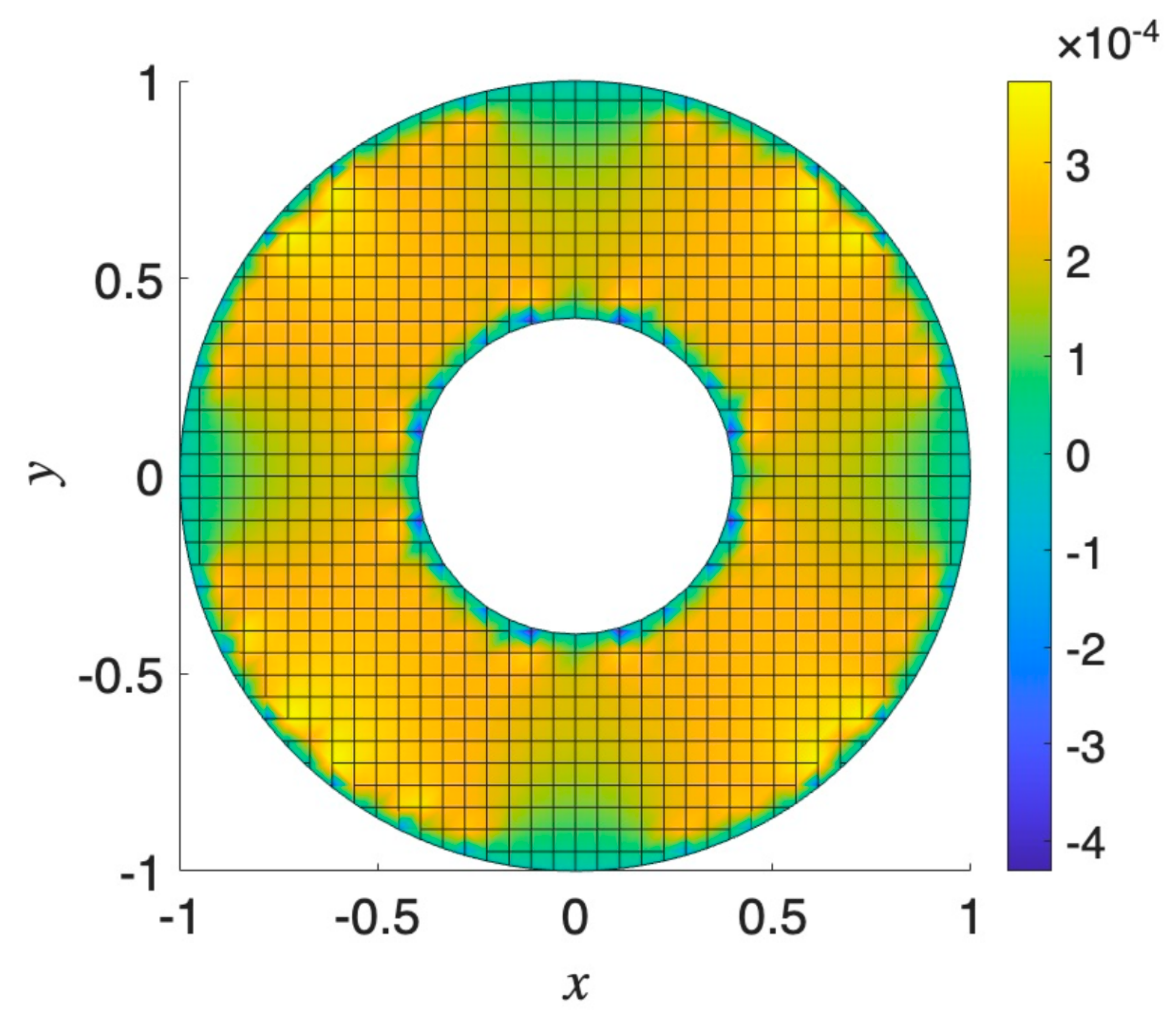}
\subcaption{}\label{fig:meshes_hole_h}
\end{subfigure}
\caption{Numerical solution using CutVEM for the
         heat conduction problem in an annular domain. 
         The annulus is immersed in triangular and uniform rectangular 
        finite element
        meshes. The third mesh in the sequence is shown. For the triangular mesh,
        (a)--(d) are the
        cut-mesh, 
        agglomerated mesh, 
        exact solution and error in CutVEM, and
        the corresponding plots for the
        rectangular mesh are in
        (e)--(h).} 
        \label{fig:meshes_hole}
\end{figure}

\begin{figure}[!thb]
\centering
\begin{subfigure}{0.48\textwidth}
\includegraphics[width=\textwidth]{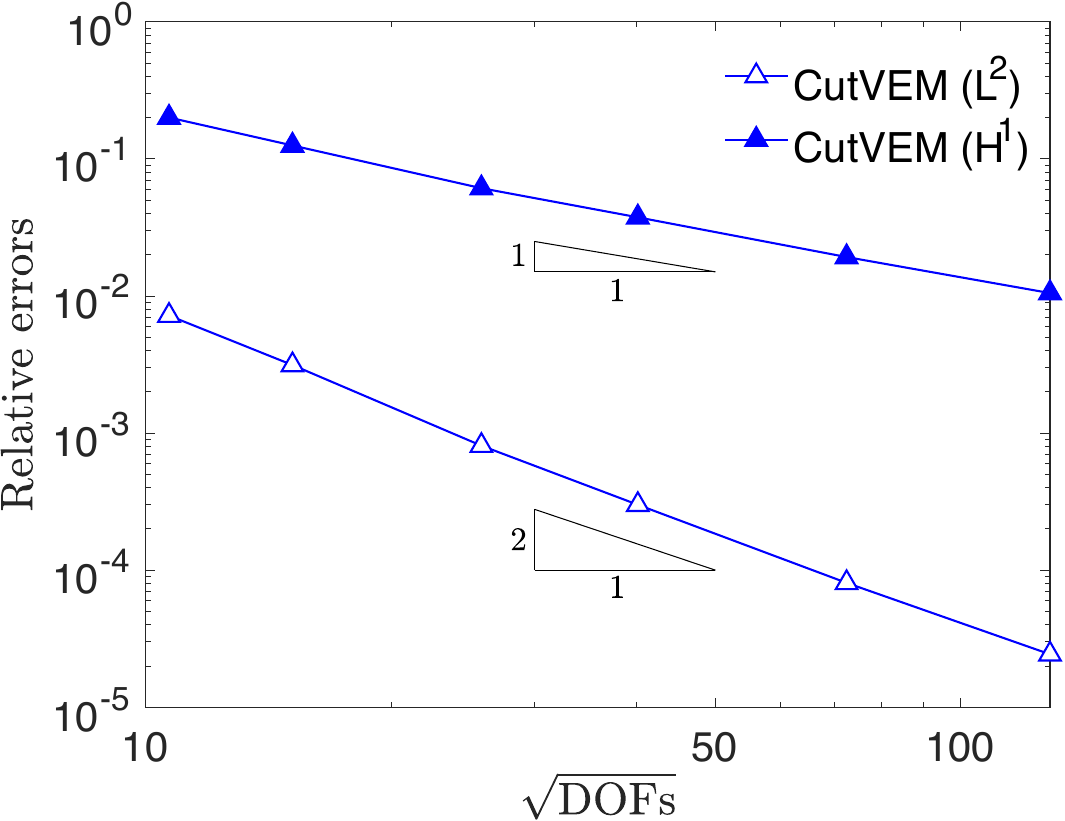}
\subcaption{}\label{fig:convergence_hole_a}
\end{subfigure} \hfill
\begin{subfigure}{0.48\textwidth}
\includegraphics[width=\textwidth]{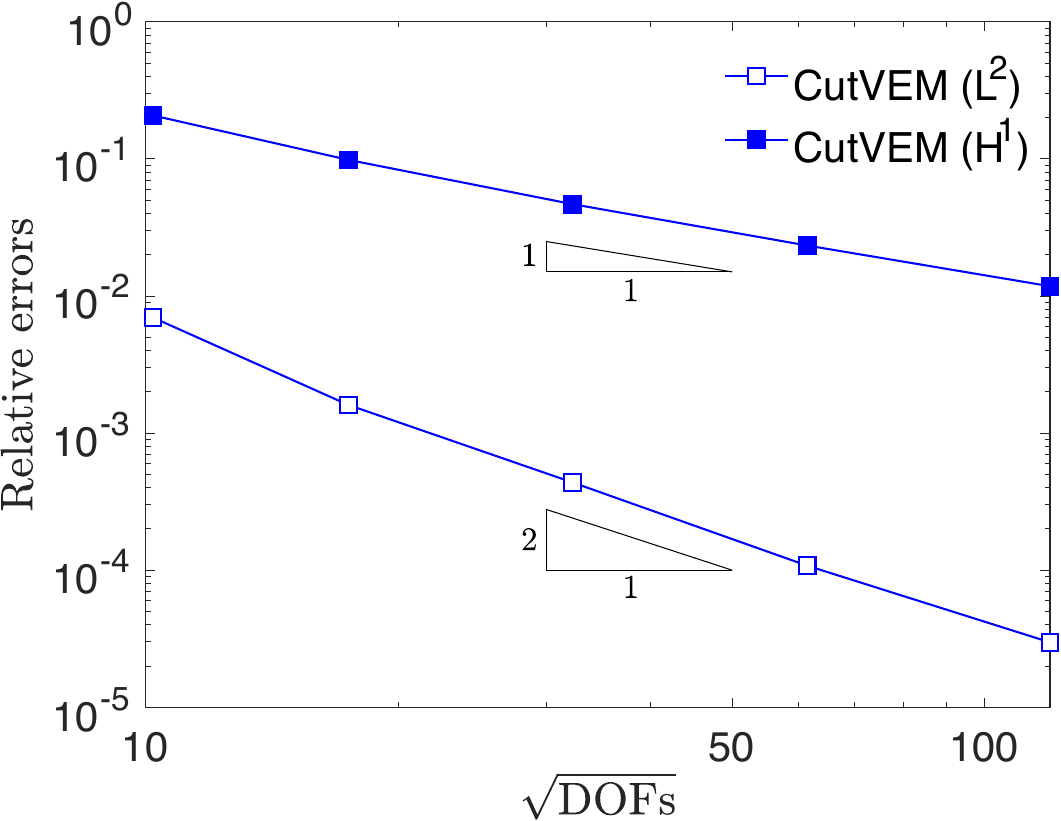}
\subcaption{}\label{fig:convergence_hole_b}
\end{subfigure}
\caption{Convergence study with CutVEM for the heat   
         conduction problem in an   
         annular domain. Agglomerated meshes are formed
         from
         an annulus that is immersed in a (a) triangular
         mesh and in a (b) quadrilateral mesh.}   
        \label{fig:convergence_hole}
\end{figure}

In the numerical computations, we use $a = 0.4$ and
$b = 1$.  The annulus is immersed in the six
triangular meshes considered in~\sref{subsec:homogeneous_heat_conduction}. 
In addition, a
sequence consisting of five 
uniformly refined meshes (square elements) is also considered.  For the third mesh in both sequences,
the cut-cell mesh, agglomerated mesh, exact
solution and the error in the solution
obtained using CutVEM
are presented in~\fref{fig:meshes_hole}.  
\revised{CutVEM} exhibits optimal convergence in
the $L^2$ norm and $H^1$ seminorm (see~\fref{fig:convergence_hole}).

\subsubsection{Two-phase composite}\label{subsubsec:composite}
Using polar coordinates, consider a two-phase composite with an inclusion occupying the open, bounded circular disc $\Omega_1 = \{r : r < a \}$. If $\Omega = \{r : r < b\}$, then the matrix occupies the
region $\Omega_2 = \Omega - \Omega_1$.  We follow the radial solution
derived for the bimaterial elastostatic problem 
in~\cite{Sukumar:2001:MHI} and tailor it for
inhomogeneous heat conduction. The strong form is posed as: 
\begin{subequations}\label{eq:strong_form_bimaterial}
\begin{align}
-\nabla \cdot (\kappa \nabla u)  &= f \ \ \textrm{in } \Omega \subset \Re^2, \\
u &= g \ \ \textrm{on } r = b, \\
\intertext{where $f \in L^2(\Omega)$ is the source function, $g$ is the boundary data and
the scalar thermal conductivity $\kappa$ is discontinuous ($\kappa_1$ and $\kappa_2$ are distinct
constants):}
\kappa &= \begin{cases} \kappa_1 & \textrm{in } \Omega_1 \\
                       \kappa_2 & \textrm{in } \Omega_2
         \end{cases} .
\end{align} 
\end{subequations}
As a verification test, we consider a manufactured problem with an exact solution.
To obtain a purely radial solution
for the steady-state heat conduction problem, we assume a constant heat source
$f = 1$ and Dirichlet boundary data $g = 1$ on $r = b$, i.e., $u(b) = 1$. On using the
matching conditions for the temperature and the normal flux on the interface $r = a$, namely
\begin{equation}\label{eq:bimaterial_ICs}
u(a^-) = u(a^+), \quad
\kappa_1 \left. \frac{du}{dr} \right|_{r = a^-} = 
\kappa_2 \left. \frac{du}{dr} \right|_{r = a^+} ,
\end{equation}
we can write down the exact solution for the strong form of the boundary-value problem 
in~\eqref{eq:strong_form_bimaterial} as:
\begin{equation}
 u(r) = \begin{cases}
     1 + \dfrac{b^2 - a^2}{4 \kappa_2} 
     + \dfrac{a^2 - r^2}{4 \kappa_1} 
     , &0 \leq r < a \\
     1 + \dfrac{b^2 - r^2}{4 \kappa_2} , &a \leq r \leq  b
 \end{cases} .
\end{equation}

\begin{figure}
\centering
\begin{subfigure}{0.25\textwidth}
\includegraphics[width=\textwidth]
{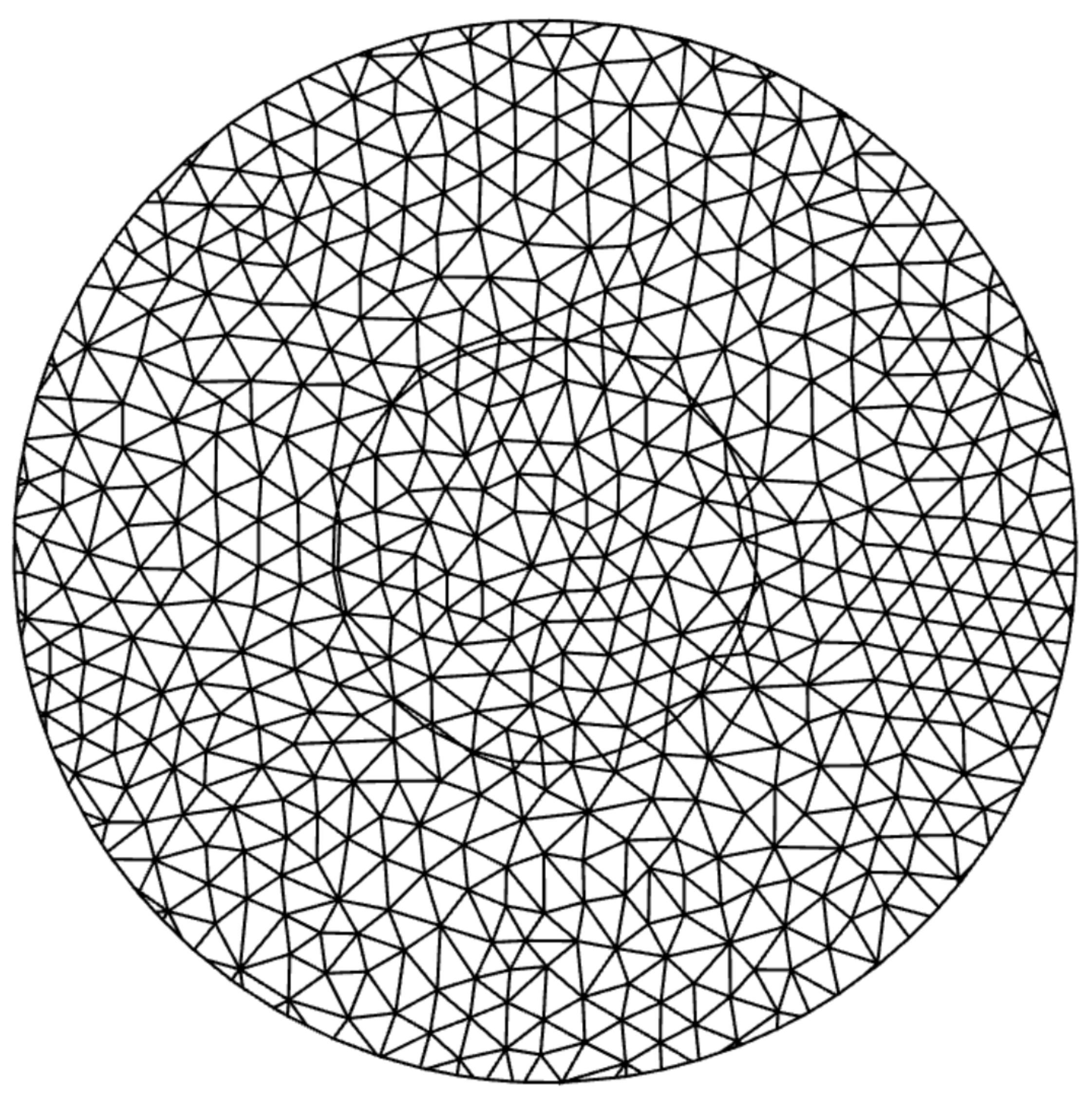}
\subcaption{}\label{fig:meshes_bimaterial_a}
\end{subfigure} \hspace*{0.9in} 
\begin{subfigure}{0.25\textwidth}
\includegraphics[width=\textwidth]
{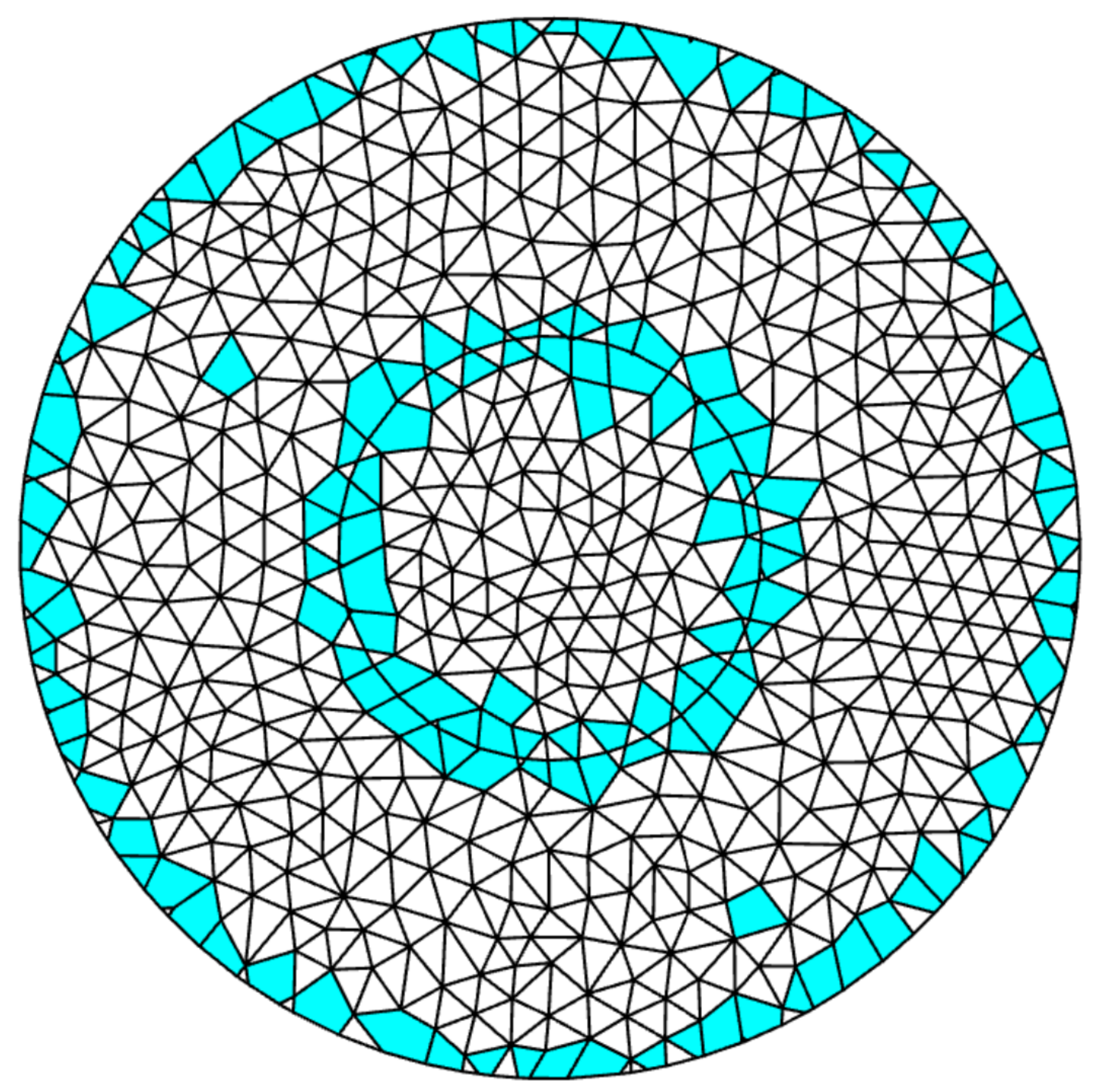}
\subcaption{}\label{fig:meshes_bimaterial_b}
\end{subfigure} \\
\begin{subfigure}{0.24\textwidth}
\includegraphics[width=\textwidth]
{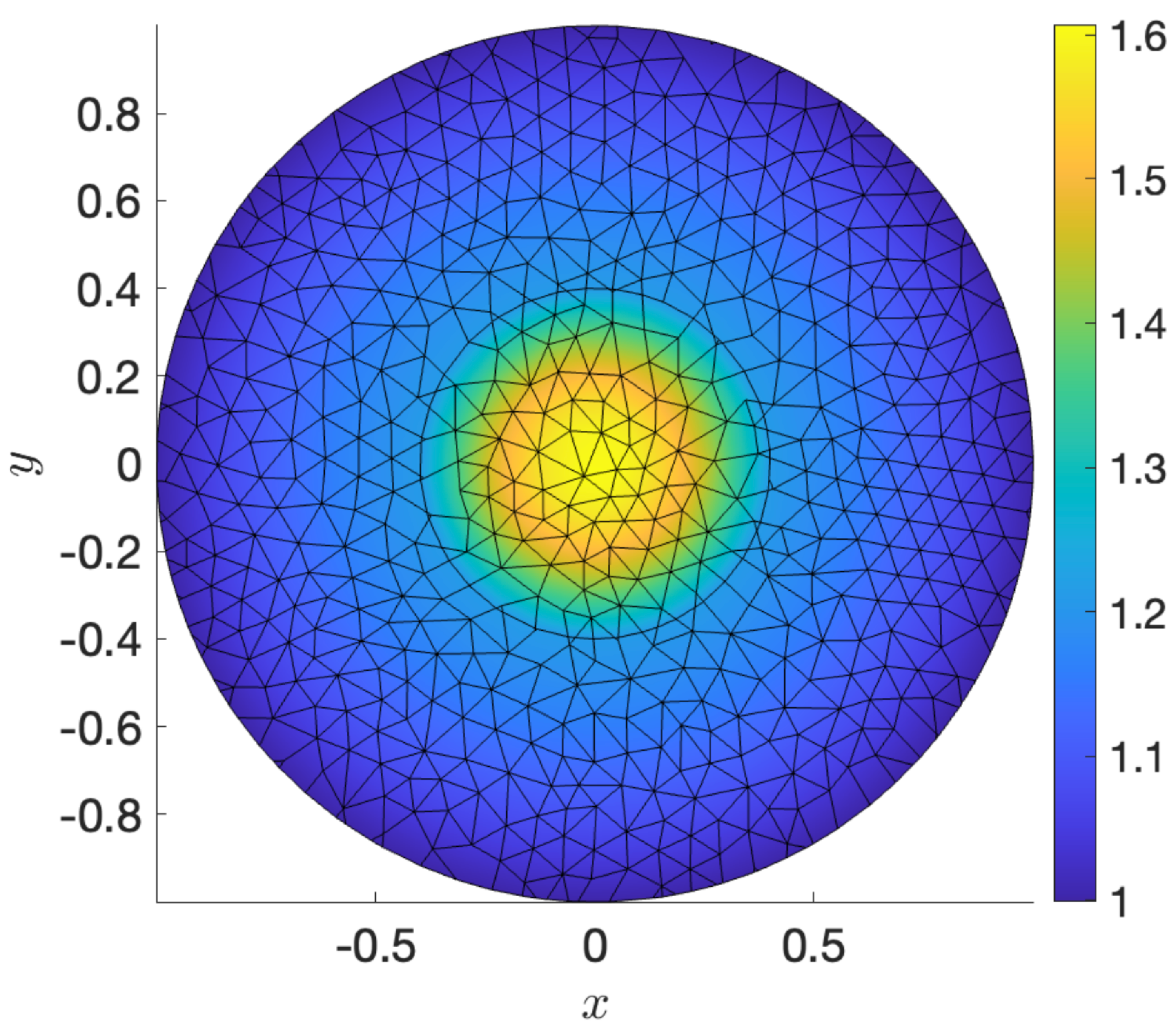}
\subcaption{}\label{fig:meshes_bimaterial_c}
\end{subfigure} \hfill
\begin{subfigure}{0.248\textwidth}
\includegraphics[width=\textwidth]
{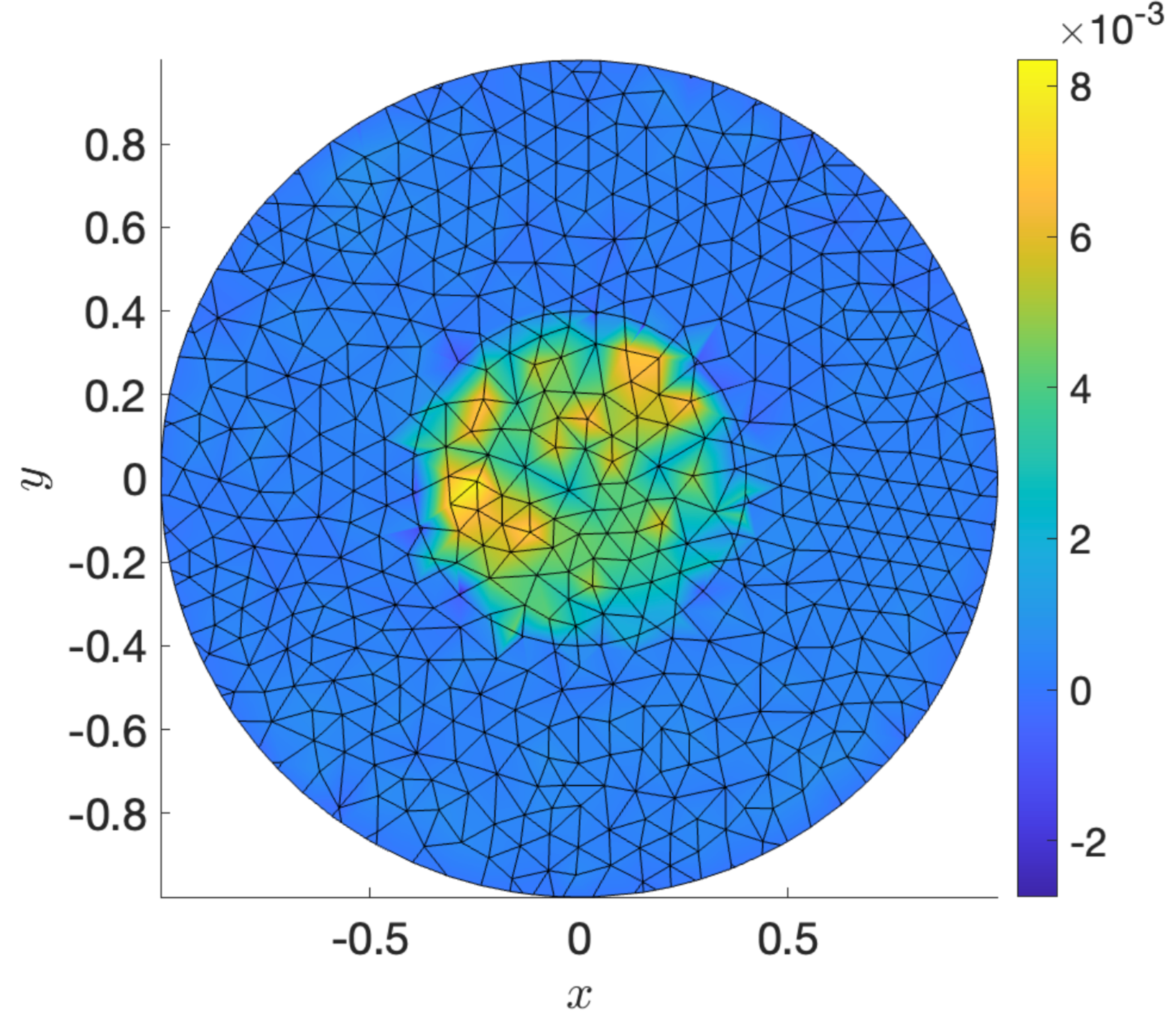}
\subcaption{}\label{fig:meshes_bimaterial_d}
\end{subfigure} \hfill
\begin{subfigure}{0.24\textwidth}
\includegraphics[width=\textwidth]
{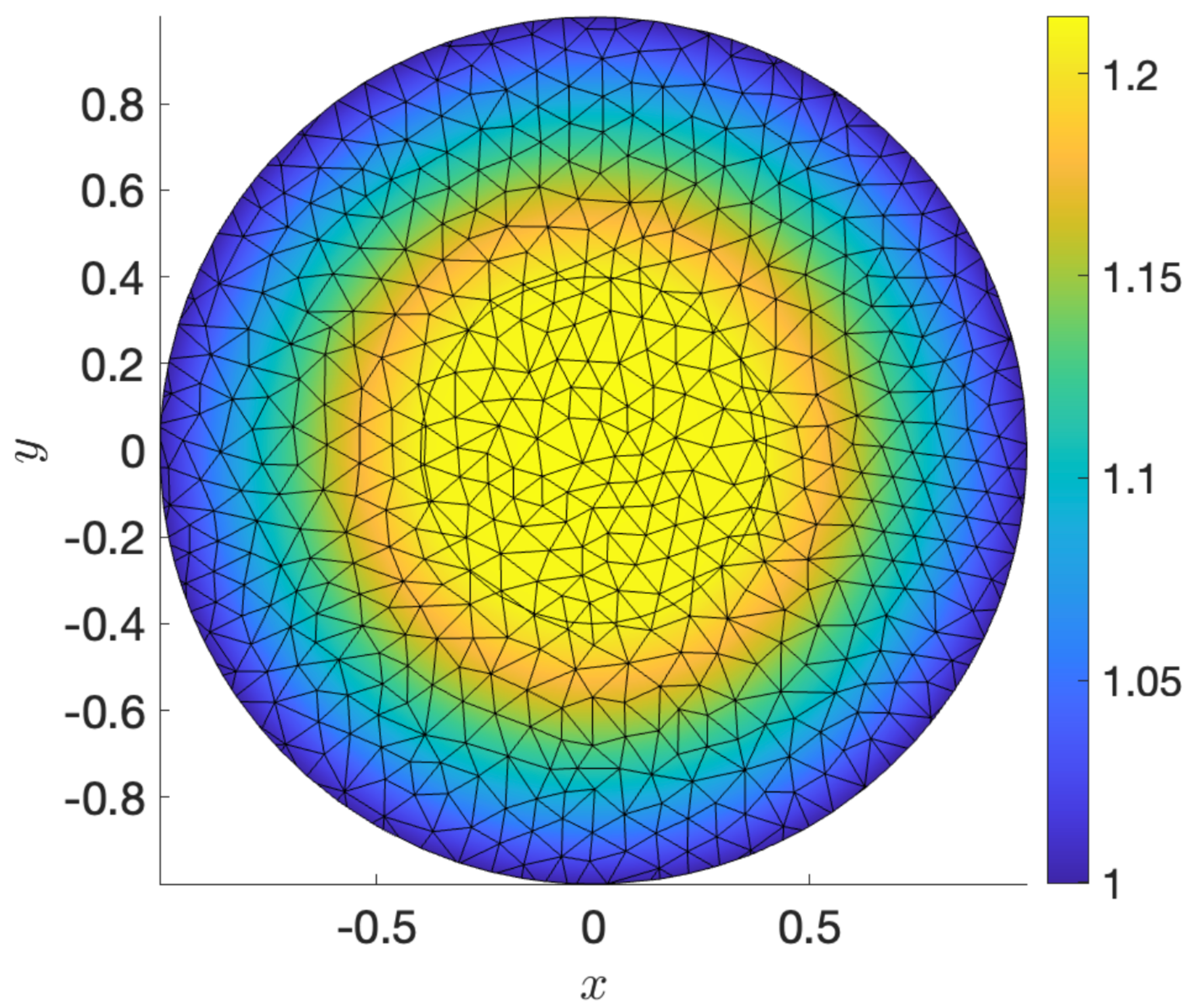}
\subcaption{}\label{fig:meshes_bimaterial_e}
\end{subfigure} \hfill
\begin{subfigure}{0.24\textwidth}
\includegraphics[width=\textwidth]
{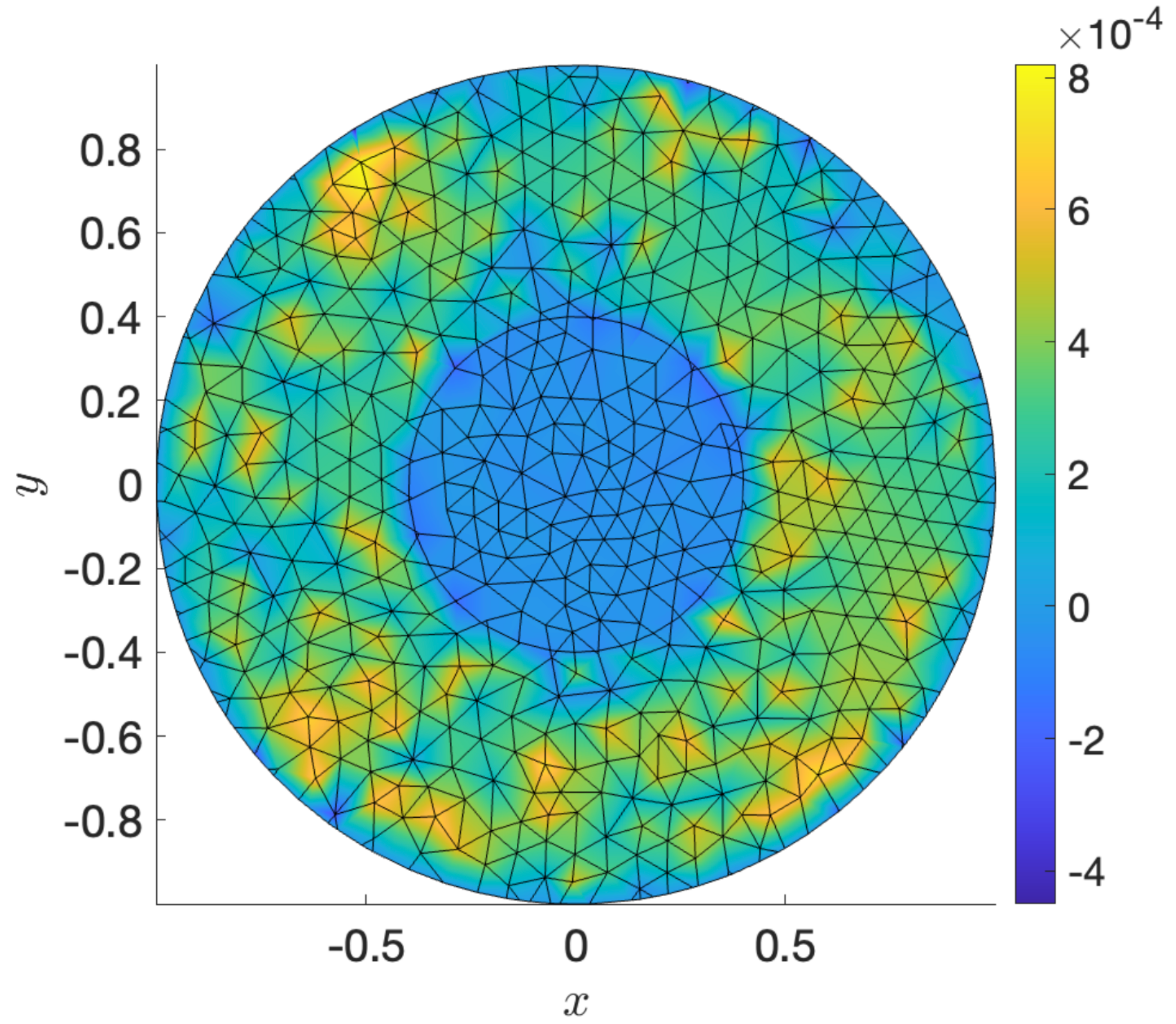}
\subcaption{}\label{fig:meshes_bimaterial_f}
\end{subfigure}
\caption{Numerical solution using CutVEM for the   
         two-phase bimaterial problem.  
         The composite disc is immersed in
         triangular finite element
         meshes. The third mesh in the sequence 
         is shown: (a) cut-mesh and
         (b) agglomerated mesh.  (c), (d) Exact solution and
         error in CutVEM 
         for $\kappa_1 / \kappa_2 = 0.1$.
         (e), (f) Exact solution and
         error in CutVEM 
         for $\kappa_1 / \kappa_2 = 10$.}
        \label{fig:meshes_bimaterial}
\end{figure}

\begin{figure}
\centering
\begin{subfigure}{0.48\textwidth}
\includegraphics[width=\textwidth]
{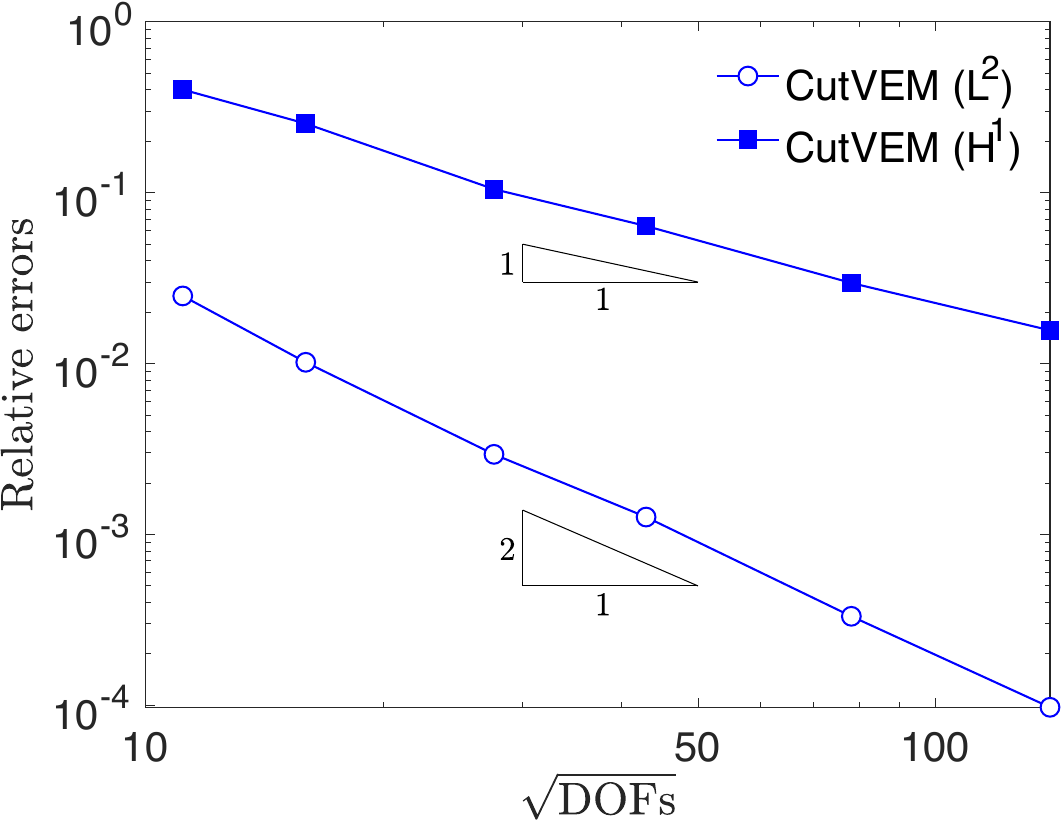}
\subcaption{}\label{fig:convergence_bimaterial_a}
\end{subfigure} \hfill
\begin{subfigure}{0.48\textwidth}
\includegraphics[width=\textwidth]
{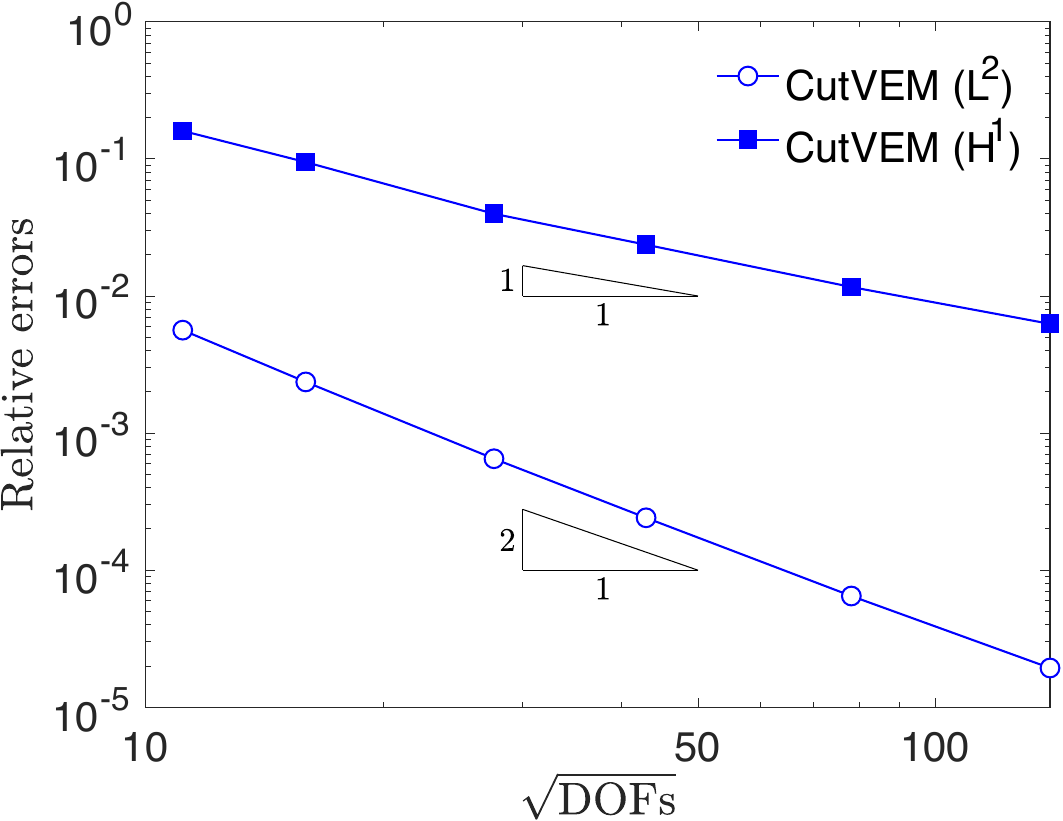}
\subcaption{}\label{fig:convergence_bimaterial_b}
\end{subfigure}
\caption{Convergence study with CutVEM for 
         the two-phase bimaterial heat   
         conduction problem. Agglomerated meshes are formed
         from
         a composite disc that is immersed in triangular
         meshes. (a) 
         $\kappa_1 / \kappa_2 = 0.1$;
         and (b) 
         $\kappa_1 / \kappa_2 = 10$.}   
        \label{fig:convergence_bimaterial}
\end{figure}

As in~\sref{subsubsec:hole}, the domain is immersed in six
triangular meshes. For the third mesh in the sequence,
the cut-cell mesh, agglomerated mesh, exact
solution and the error in the solution obtained
using CutVEM are presented in~\fref{fig:meshes_bimaterial}.  
Two cases are considered with contrasts in thermal
conductivity coefficient: 
$\kappa_1 / \kappa_2 = 0.1$
and $\kappa_1 / \kappa_2 = 10$.
In the numerical computations, we choose $\tau = \kappa_1$ (elements 
that lie within the inclusion) or $\tau = \kappa_2$ (elements that lie in the matrix)
as the stability parameter in the VEM. The
convergence curves for the two cases are
presented in~\fref{fig:convergence_bimaterial}.
We observe that CutVEM converges optimally 
with convergence
rates of $2$ and $1$ in the $L^2$ norm and the
$H^1$ seminorm, respectively.

\section{Conclusions}
\label{sec:conclusions} 
Computational methods that embed
domains and interfaces in a finite element mesh generate cut-cells
that are in general polygons.  Such methods require special treatment
to handle ill-conditioning of basis functions 
caused by cut-cells, satisfying boundary and interfacial conditions, and
numerical integration over polygons, cf.~\cite{dePrenter:2023:SCI}.
It is then appealing to resort to the
VEM~\cite{Beirao:2013:BPV,Beirao:2023:VEM}, which is tailored for
polygonal partitions and is known to be robust to poorly-shaped
polygons (triangles being an exception). Even though this
approach helps avoid intricate remeshing operations typical in
embedded finite element methods and facilitates imposing
boundary conditions and interfacial constraints, 
it is still prone to ill-conditioned stiffness matrices. 
Drawing on the idea of element agglomeration introduced
in~\cite{Sukumar:2022:AVE},
the CutVEM proposed in this work directly addresses this difficulty,
with a robust shape-agnostic element agglomeration technique as the
remedy. The CutVEM rids VEM of ill-conditioning caused by cut-cells
in meshes with embedded geometries while also leveraging it
to accommodate polygonal elements resulting from agglomeration.
Furthermore, numerical integration over polygons (without the need for
subpartitioning) is realized through the scaled boundary cubature
scheme~\cite{Chin:2021:SBC}.  The CutVEM is thus an enhanced
domain-conforming VEM that is 
well-suited for domains with embedded geometries.

Existing approaches to element agglomeration are motivated by
techniques in the meshing literature, relying primarily on polygon
shape-quality metrics.  Our work here is a departure from this
viewpoint and introduces the element stability ratio, \revised{the
  inverse of the element-level condition number,} as an indicator of
element quality.  \revised{The element-agglomeration criterion is not
  computationally expensive to evaluate since it does not require
  accurate eigenvalues; relaxed tolerances and low-precision
  calculations suffice.}  Using examples of triangle and quadrilateral
elements, we highlighted the relationship between an element's shape
and its stability ratio to be tenuous at best.  We showcased element
agglomeration as a simple yet potent operation that improves the
conditioning of the VEM over meshes with embedded geometries. This key
idea, coupled with the robustness of the VEM on poorly-shaped
polygons, is the basis for the CutVEM.  Distinctively, the
agglomeration algorithm underlying the CutVEM:
\begin{itemize}
\item is agnostic to element shapes, seeking to improve element
  spectra rather than element 
  appearances;\footnote{\revised{It is known
that even on extremely degenerate polygons, interpolation errors using generalized barycentric 
can be small, which underscores the ambiguous connection 
between polygon shape quality and its impact on interpolation error~\citep{Gillette:2017:SQG}.}}
\item is method-aware, i.e., it explicitly relies on computing
  problem-specific element stiffness matrices in the VEM;
\item is iterative, and attempts to improve poorer 
elements even if at the cost of better ones;
\item is atomic in nature, requiring inspection of elements
  individually rather than collectively, and is hence trivially parallelizable as well;
\item does not require resolving large systems of 
equations or optimization problems; and
\item is conveniently implemented as a mesh preprocessing step \revised{for stationary embedded geometries}.
\end{itemize}

Through extensive numerical experiments with immersed 
\revised{(stationary and evolving)}
geometries, we provided
convincing evidence that the CutVEM enjoys dramatically
improved condition numbers of global stiffness matrices over the VEM
despite relying on an element-level agglomeration criterion.
Furthermore, we found that the reduced element count due to
agglomeration did not noticeably affect solution accuracy for several
heat conduction (homogeneous and inhomogeneous) problems.  We
consistently found that the CutVEM retains the optimal convergence
rates of the VEM. Though our simulations were limited to prototypical
heat conduction problems, we fully expect the improved conditioning
and optimal convergence properties of the CutVEM
to also hold for coercive elliptic problems in general.
  
As a well-conditioned conforming method for embedded domains, CutVEM
is immediately useful for a challenging class of moving boundary
problems. This includes simulating crack propagation, problems with
multimaterial
interfaces, and shape/topology optimization applications. To aid such
investigations with CutVEM, we hope to share an implementation of the
agglomeration algorithm as a conveniently interfaceable
library. Efforts are also underway to examine the generalization of
the proposed method to simulate three-dimensional problems with immersed geometries.

\section*{Acknowledgment}
RR gratefully acknowledges financial support from the Anusandhan
National Research Foundation (ANRF) through the 
project
CRG/2022/006569.

\end{document}